\documentclass[12pt]{article}
\usepackage{amsmath}
\usepackage{graphicx,psfrag,epsf}
\usepackage{enumerate}
\usepackage{natbib}
\usepackage{url} 

\usepackage{bm}
\usepackage{amsfonts}
\usepackage{subcaption}
\usepackage{multirow}
\usepackage{enumitem}
\usepackage{yhmath}

\newcommand{\blind}{0}

\newtheorem{theorem}{Theorem}

\newtheorem{assumption}[theorem]{Assumption}
\newtheorem{lemma}[theorem]{Lemma}

\newcommand{\eop}{\hfill $\Box$ \\ \\}

\begin{document}

\def\spacingset#1{\renewcommand{\baselinestretch}%
{#1}\small\normalsize} \spacingset{1}


\if0\blind
{
  \title{\bf Grenander--Stone estimator: stacked constrained estimation of a discrete distribution over a general directed acyclic graph}
  \author{Vladimir Pastukhov\\
  Department of Statistics and Operations Research,\\
University of Vienna, Austria
}
  \maketitle
} \fi

\if1\blind
{
  \bigskip
  \bigskip
  \bigskip
  \begin{center}
    {\LARGE\bf Grenander--Stone estimator: stacked constrained estimation of a discrete distribution over a general directed acyclic graph}
\end{center}
  \medskip
} \fi

\bigskip
\begin{abstract}
In this paper we integrate isotonic regression with Stone's cross-validation-based method to estimate a distribution with a general countable support with a partial order relation defined on it. We prove that the estimator is strongly consistent for any underlying distribution, derive its rate of convergence, and in the case of one-dimensional support we obtain Marshal-type inequality for cumulative distribution function of the estimator.  Also, we construct the asymptotically correct conservative global confidence band for the estimator. It is shown that, first, the estimator performs good even for small sized data sets, second, the estimator outperforms in the case of non-isotonic underlying distribution, and, third, it performs almost as good as Grenander estimator when the true distribution is isotonic. Therefore, the new estimator provides a trade-off between goodness-of-fit, monotonicity and quality of probabilistic forecast. We apply the estimator to the time-to-onset data of visceral leishmaniasis in Brazil collected from $2007$ to $2014$.
\end{abstract}

\noindent%
\vfill

\newpage
\spacingset{1.45} 

\section{Introduction}\label{sec:intro}

This work is motivated by the papers in constrained estimation of probability mass functions and by the research in stacked estimators. The first paper about estimation of discrete monotone distributions is \cite{jankowski2009estimation}, where the authors studied constrained maximum likelihood estimator (MLE) and the method of rearrangement. The MLE of a monotone distribution is also known as Grenander estimator. Further, in the papers \cite{balabdaoui2017asymptotics,durot2013least} the authors studied the least squares estimator of a convex discrete distribution. The MLE of log-concave distribution was studied in detail in \cite{balabdaoui2013asymptotics, jankowski2018estimating} and the MLE of a unimodal discrete distribution was considered in \cite{balabdaoui2016maximum}. The estimation of probability mass function (p.m.f.) with unknown labels was discussed in \cite{anevski2017estimating}. Estimation of a completely monotone p.m.f. was addressed in \cite{balabdaoui2020least, balabdaoui2020completely}.

In our paper we combine Grenander estimator with Stone's cross-validation-based concept \citep{stone1974cross} to estimate discrete  infinitely supported distribution which is possibly isotonic with respect to the partial order on its support. 

There are several papers on stacked estimation of discrete distributions, for example, \cite{fienberg1972choice, fienberg1973simultaneous, pastukhov2022stacked,  stone1974cross, trybula1958some}, and the case of continuous support was considered, for example, in \cite{rigollet2007linear, smyth1999linearly}.  To the authors' knowledge stacking of shape constrained estimators has not been studied much. At the same time stacking of estimators is proved to be a powerful method. For example, in the recent paper \cite{hastie2020best} it was shown that in terms of prediction accuracy the stacking of least squares estimator and LASSO performs almost equally to the LASSO in low signal-to-noise ratio regimes, and nearly as well as the best subset selection in high signal-to-noise ratio scenarios. In the paper \cite{smyth1999linearly} it was shown that stacked estimator performs better than selecting the best model by cross-validation. Stacked estimation of one-dimensional constrained estimator of a discrete distribution based on cross-validation estimation of $\ell_{2}$ distance was studied in \cite{pastukhov2022stacked}. 

In most of the papers on constrained discrete density estimation the authors separated the well- and the mis-specified cases and, except for the paper \cite{jankowski2018estimating}, only one-dimensional support was studied. In our work we assume that the true p.m.f. can be non-monotone and prove that the proposed estimator is strongly consistent with $\sqrt{n}$-rate of convergence even if the true p.m.f. is not monotone. Also, we study general case of distribution over a directed acyclic graph.

The idea of the estimator introduced in this paper is in some sense similar to nearly-isotonic regression, cf. \cite{tibshirani2011nearly} and \cite{minami2020estimating, pastukhov2023fused} for multidimensional case. Nearly-isotonic regression provides intermediate less restrictive penalized solution and the isotonic regression is in the solution path. 

The estimator studied in this paper has several important applications. In medicine, a possible application is a monotonic smoothing of bivariate histogram of population based body mass index deciles of mothers and fathers of a group of children with obesity, cf. \cite{hebebrand2000epidemic}. In reliability theory, the application is smoothing the estimator of bivariate distribution, cf. \cite{roy1993reliability}. In machine learning and data mining, an example of discrete isotonic distribution on a tree graph is isotonic data augmentation for knowledge distillation cf. \cite{cui2021isotonic}.

\subsection{Notation}\label{statement_pr}
Assume that $\mathcal{Z}_{n} =  \{\bm{z}_{1}, \bm{z}_{2}, \dots, \bm{z}_{n} \}$ is a sample of $n$ i.i.d. observations with values in some countable set $\mathcal{S} = \{\bm{s}_{0}, \bm{s}_{1}, \bm{s}_{2}, \dots \}$ and generated by p.m.f. $\bm{p}$. For a given data sample $\mathcal{Z}_{n}$ let $\bm{x}$ be the frequency data, i.e. $\bm{x} = (x_{0}, x_{1}, \dots,  x_{t_{n}})$, where $x_{j}= \sum_{i=1}^{n} 1 \{ \bm{z}_{i} = \bm{s}_{j} \}$ and $t_{n} = \sup\{j: x_{j}>0 \}$. Next, $\mathcal{P} = \{ f \in \ell_{1}: \, \sum_{j} f_{j} =1, \, f_{j} \geq 0  \}$ denotes probability simplex. 

Further, for a given data sample $\mathcal{Z}_{n}$ let  
\begin{equation*}\label{}
\hat{p}_{n, \bm{s}_{j}} \equiv \hat{p}_{n, j} = \frac{x_{j}}{n}, \text{ with } j \in \mathbb{N}, 
\end{equation*}
be the empirical estimator of $\bm{p}$, and each element $\bm{z}_{i}$ in the sample $\mathcal{Z}_{n}$ be associated with a multinomial indicator 
\begin{equation*}\label{}
\bm{\delta}^{(i)} \equiv \bm{\delta}^{(\bm{z}_{i})}  = (0, \dots, 0, 1, 0, \dots, 0) \quad (i = 1, \dots, n).
\end{equation*}
The multinomial indicator $\bm{\delta}^{(i)}$ is a vector in $\mathbb{R}^{t_{n}+1}$ with all elements equal to zero, except for the one with the index $l$ such that $\bm{z}_{i} = \bm{s}_{l}$, cf. \citep{bowman1984cross, stone1974cross}. The vector $\bm{\delta}^{(i)}$ is the degenerate distribution associated with the observation $\bm{z}_{i}$ in the observed sample $\mathcal{Z}_{n}$. 

Further, for a frequency vector $\bm{x}$, we associate each component $x_{j}$ of $\bm{x}$ with another multinomial indicator $\bm{\sigma}^{[j]} \in \mathbb{R}^{t_{n}+1}$, given by  
\begin{equation}\label{indgm}
\bm{\sigma}^{[j]} \equiv \bm{\sigma}^{[x_{j}]} = (0, \dots, 0, 1, 0, \dots, 0) \quad (j = 0, \dots, t_{n}).
\end{equation}
All elements of $\bm{\sigma}^{[j]}$ are zeros, except for the one with the index $j$.

Let us define the following binary relation $\preceq$ on the support $\mathcal{S}$ of p.m.f. $\bm{p}$. Assume that for $\preceq$ on $\mathcal{S}$ the following holds:
\begin{enumerate}[label=(\roman*)]
\item it is reflexive, i.e. $\bm{s} \preceq \bm{s}$ for $\bm{s} \in \mathcal{S}$;
\item it is transitive, i.e. $\bm{s}_{i}, \bm{s}_{j}, \bm{s}_{k} \in \mathcal{S}$, $\bm{s}_{i} \preceq \bm{s}_{j}$ and $\bm{s}_{j} \preceq \bm{s}_{k}$ imply $\bm{s}_{i} \preceq \bm{s}_{k}$;
\item it is antisymmetric, i.e. $\bm{s}_{i}, \bm{s}_{j} \in \mathcal{S}$, $\bm{s}_{i} \preceq \bm{s}_{j}$ and $\bm{s}_{j} \preceq \bm{s}_{i}$ imply $\bm{s}_{i} = \bm{s}_{j}$.
\end{enumerate}

The order relation $\preceq$ defined above is called partial order. For any partial order relation $\preceq$ on $\mathcal{S}$ there exists directed acyclic graph $G = (\mathcal{S}, E)$, where $E$ is the set of edges such that
\begin{eqnarray}
E = \{(\bm{s}_{i}, \bm{s}_{j}), \, \text{where} \, (\bm{s}_{i}, \bm{s}_{j}) \,\text{ is the ordered pair of vertices from } \, \mathcal{S}\},
\end{eqnarray}
such that for any $\bm{s}_{k}, \bm{s}_{l} \in \mathcal{S}$ one has $\bm{s}_{k} \preceq \bm{s}_{l}$ iff $E$ contains the chain of edges from $\bm{s}_{k}$ to  $\bm{s}_{l} $. 

A probability vector $\bm{\phi} \in \mathcal{P}$ with the support $\mathcal{S}$ is isotonic with respect to the partial order $\preceq$ on $\mathcal{S}$ if for any $\bm{s}_{i}, \bm{s}_{j} \in \mathcal{S}$ such that $\bm{s}_{i} \preceq \bm{s}_{j}$ we have $\phi_{\bm{s}_{i}} \leq \phi_{\bm{s}_{j}}$. Equivalently, given the directed acyclic graph $G = (\mathcal{S}, E)$, corresponding to the partial order $\preceq$ on $\mathcal{S}$, a probability vector $\bm{\phi}$ with the support $\mathcal{S}$ is isotonic, if for any $\bm{s}_{k}, \bm{s}_{l} \in \mathcal{S}$ we have $\phi_{\bm{s}_{k}} \leq \phi_{\bm{s}_{l}}$ iff $E$ contains the chain of edges from $\bm{s}_{k}$ to  $\bm{s}_{l} $. Let $\mathcal{P}^{is}$ denote the subset of all isotonic probability vectors with the support $\mathcal{S}$. Note that the set $\mathcal{P}^{is}$ is a cone in $\ell_{2}$ space.

For any probability vector $\bm{g} \in \mathcal{P}$ with the support $\mathcal{S}$ the vector $\hat{\bm{g}}^{*}$ is isotonic regression with respect to the partial order $\preceq$, defined on $\mathcal{S}$, if
\begin{eqnarray*}
\bm{g}^{*} = \Pi^{is}(\bm{g}  \mid \mathcal{P}^{is}) = \underset{\bm{f} \in \mathcal{P}^{is}}{\arg \min} \sum_{\bm{s} \in \mathcal{S}}(f_{\bm{s}} - g_{\bm{s}})^{2},
\end{eqnarray*}
where $ \Pi^{is}(\bm{g}  \mid \mathcal{P}^{is})$ denotes the $\ell_{2}$-projection of $\bm{g}$ onto $\mathcal{P}^{is}$.

Equivalently, in terms of directed acyclic graph $G = (\mathcal{S}, E)$, a probability vector $\hat{\bm{g}}^{*} \in \mathcal{P}^{is}$ with the support $\mathcal{S}$ is isotonic regression of a vector $\bm{g}$ with respect to the graph $G = (\mathcal{S}, E)$ if 
\begin{eqnarray*}
\begin{aligned}
&\bm{g}^{*} = \arg \min \sum_{\bm{s} \in \mathcal{S}}(f_{\bm{s}} - g_{\bm{s}})^{2},\\
&\text{s.t.} \quad  f_{\bm{s}_{i}} \leq f_{\bm{s}_{j}} \quad  \text{for any} \quad (\bm{s}_{i}, \bm{s}_{i}) \in E.
\end{aligned}
\end{eqnarray*}

Note, that the MLE of p.m.f. with isotonic constraints, which we denote by $\hat{\bm{g}}_{n}$, is equivalent to the isotonic regression of the empirical estimator \citep{brunk1972statistical, jankowski2009estimation, robertson1980algorithms}
\footnote{A confusing notion of "isotonic regression" is a standard notion in constrained inference, eventhough it is more natural to call it as $\ell_{2}$ projection of a vector onto the isotonic cone.}, i.e.
\begin{equation}\label{Grest}
\hat{\bm{g}}_{n} =\Pi^{is}(\hat{\bm{p}}_{n}  \mid \mathcal{P}^{is}) =  \underset{\bm{f} \in \mathcal{P}^{is}}{\arg \min}\sum_{j=0}^{t_{n}}(\hat{p}_{n, \bm{s}_{j}} - f_{\bm{s}_{j}})^{2}.
\end{equation}
The estimator $\hat{\bm{g}}_{n}$ is called Grenander estimator of isotonic distribution. The solution to (\ref{Grest}) in the case when $\mathcal{S} = \{0, 1, 2, \dots \}$, i.e. one dimensional MLE of increasing p.m.f., is given by one of the following max-min formulas, cf. p. 19 in \citep{brunk1972statistical}:
\begin{equation}\label{minmaxis}
\begin{aligned}
\hat{p}^{MLE}_{n, j}&={} \max_{v\leq j} \min_{w\geq j} \frac{\sum_{s=v}^{w}\hat{p}_{n, s}}{w-v +1} = \min_{w\geq j}\max_{v\leq j} \frac{\sum_{s=v}^{w}\hat{p}_{n, s}}{w-v +1} \\
&=\max_{v\leq j} \min_{w \geq v} \frac{\sum_{s=v}^{w}\hat{p}_{n, s}}{w-v +1} = \min_{w\geq j}\max_{v\leq w} \frac{\sum_{s=v}^{w}\hat{p}_{n, s}}{w-v +1},
\end{aligned}
\end{equation}
where the first equality in (\ref{minmaxis}) can be represented graphically as the slope of the greatest convex majorant of the empirical distribution function for the sample $\mathcal{Z}_{n}$. The solution for the case of a general partial order is more complicated and it is discussed in the supplementary material of the paper.

Further, for an arbitrary vector $\bm{f}\in\ell_{k}$ we define $\ell_{k}$-norm
\begin{equation*}\label{}
  ||\bm{f}||_{k}=\begin{cases}
    \Big(\sum_{j=0}^{\infty}|f_{j}|^{k}\Big)^{1/k}, & \text{if } \, k \in \mathbb{N} \backslash \{0\}, \\
    \sup_{j\in\mathbb{N}}|f_{j}|, & \text{if} \,  k=\infty,
  \end{cases}
\end{equation*}
and for $\bm{v} \in \ell_{2}$ and $\bm{w} \in \ell_{2}$ let $\langle \bm{v}, \bm{w} \rangle = \sum_{j=0}^{\infty}v_{j}w_{j}$ denote the inner product in $\ell_{2}$ space.

For a random sequence $b_{n} \in \mathbb{R}$ we will use notation
\begin{equation*}
b_{n} = O_{p}(n^{q}),
\end{equation*}
if for any $\varepsilon >0$ there exists  a finite $M > 0$ and a finite $N > 0$ such that  
\begin{equation*}
\mathbb{P}(n^{-q}|b_{n}| > M) < \varepsilon,
\end{equation*}
for any $n > N$.

\subsection{Statement of the problem}\label{pr_form}

In our work we construct the estimator of a p.m.f. $\bm{p}$ in the following way. 
\begin{equation}\label{defest}
\hat{\bm{\phi}}_{n} = \hat{\beta}_{n}  \Pi(\hat{\bm{p}}_{n}) + (1-\hat{\beta}_{n})\hat{\bm{p}}_{n},
\end{equation}
where  $\Pi$ is a map, such that $\Pi: \mathcal{P} \to \mathcal{P}$, and $\hat{\beta}_{n}$ is cross-validation selected mixture parameter:
\begin{equation*}
\hat{\beta}_{n} = \underset{\beta \in [0, 1]}{\arg \min } \, CV(\beta, \mathcal{L}).
\end{equation*}
We emphasise here, that, in general, the estimator $ \Pi(\hat{\bm{p}}_{n})$ itself is not assumed to be a consistent estimator of $\bm{p}$.  Analogously to the approach in \cite{bowman1984cross, stone1974cross},  we select the penalisation parameter $\beta$ based on leave-one-out estimation of the following criterions:
\begin{equation}\label{cvs}
CV(\beta, \mathcal{L}) = 
\begin{cases}
\mathbb{E}_{\bm{p}}(L_{1}(\bm{\delta}^{z}, \hat{\bm{\phi}}_{n})\mid\mathcal{Z}_{n}) \equiv S_{1}(\hat{\bm{\phi}}_{n}, \bm{p}), &\text{ for }  \mathcal{L} = L_{1},\\
\mathbb{E}_{\bm{p}}(L_{2}(\bm{\delta}^{z}, \hat{\bm{\phi}}_{n})\mid\mathcal{Z}_{n}) \equiv S_{2}(\hat{\bm{\phi}}_{n}, \bm{p}), &\text{ for }  \mathcal{L} = L_{2},\\
\end{cases}
\end{equation}
where 
\begin{equation*}\label{}
\begin{aligned}
L_{1}(\bm{\delta}, \hat{\bm{\phi}}_{n}) &= \sum_{k=0}^{t_{n}}|\delta_{k} - \hat{\phi}_{n, k}|,\\ 
L_{2}(\bm{\delta}, \hat{\bm{\phi}}_{n}) &= \sum_{k=0}^{t_{n}}(\delta_{k} - \hat{\phi}_{n, k})^{2}.
\end{aligned}
\end{equation*}

Both loss functions $L_{1}$ and $L_{2}$ were studied in \cite{stone1974cross} to perform cross-validation selection of the mixture parameter for the case of stacking the empirical estimator with a fixed probability distribution over a fixed finite domain. The $L_{1}$ loss was also studied in \cite{buja2005loss}, and the $L_{2}$ loss was considered, for example, in \cite{bowman1984cross, haghtalab2019toward}. 

The functions $S_{1}(\hat{\bm{\phi}}_{n}, \bm{p})$ and $S_{2}(\hat{\bm{\phi}}_{n}, \bm{p})$ are also called expected scoring rules under $\bm{p}$ for the probabilistic forecast $\hat{\bm{\phi}}_{n}$ for $L_{1}$ and $L_{2}$ loss, respectively, cf. \citep{buja2005loss, czado2009predictive, haghtalab2019toward}.

Next, we apply this unified approach to estimation of a possibly isotonic p.m.f. with respect to the partial-order relation on its support (or, equivalently, p.m.f. over directed acyclic graph), i.e. in (\ref{defest}) we assume that 
\begin{equation*}
\Pi(\hat{\bm{p}}_{n}) = \Pi^{is}(\hat{\bm{p}}_{n}  \mid \mathcal{P}^{is}).
\end{equation*}
Therefore, we combine Grenander estimator with Stone's cross-validation-based concept in \cite{stone1974cross} to estimate discrete  infinitely supported distribution, and call the resulting estimator 
\begin{equation*}\label{}
\hat{\bm{\phi}}_{n} = \hat{\beta}_{n}\Pi^{is}(\hat{\bm{p}}_{n}  \mid \mathcal{P}^{is}) + (1-\hat{\beta}_{n})\hat{\bm{p}}_{n}
\end{equation*}
as Grenander--Stone estimator.

\subsection{Contribution of the paper}\label{contr_p}

The contribution of the paper is the following. We provide a unified computationally feasible framework for the stacked isotonically constrained estimation of a discrete infinitely supported distribution with a probability forecast control. We obtain sufficient conditions for the stacked estimator to be strongly consistent with parametric $\sqrt{n}$-rate of convergence and construct the estimator of its asymptotically conservative global confidence band. Next, we prove that Grenander-Stone estimator satisfies these conditions and, also, in one-dimensional case we prove Marshal-type result for the cumulative distribution function of the estimator.

\subsection{Organisation of the paper}\label{org_p}
The paper is organised as follows. In Section \ref{alp_sel} we provide derivation of cross-validation based mixture parameters.  Section \ref{constest} is dedicated to the theoretical properties of the estimator, we prove consistency and derive the rate of convergence. Next, in Section \ref{confband} we construct asymptotic global confidence bands for the stacked estimator. Further,  we compare performance of the stacked estimator with empirical, minimax and Grenander estimators in Section \ref{est_sim}. In Section \ref{app} we apply the estimator to the real data sets. The article closes with discussion in Section \ref{discus}. The proofs of all results and R code for the simulations is available upon request. 

\section{Estimation of the mixture parameter}\label{alp_sel}
The cross-validation criterions defined in (\ref{cvs}), which we aim to minimize, depend on unknown $\bm{p}$ and, therefore, we have to estimate them. Let $\hat{\bm{\phi}}^{\backslash (i)}_{n} \ (i=1, \dots, n)$ denote the leave-one-out version of the stacked estimator $\hat{\bm{\phi}}_{n}$ with respect to data sample $\mathcal{Z}_{n}$ i.e.
\begin{equation}\label{cvestimi}
\hat{\bm{\phi}}^{\backslash (i)}_{n} =  \beta \, \hat{\bm{g}}^{\backslash (i)}_{n} + (1-\beta)\hat{\bm{p}}^{\backslash (i)}_{n},
\end{equation}
where
\begin{equation*}
\hat{\bm{p}}^{\backslash (i)}_{n} = \frac{\bm{x} - \bm{\delta}^{(i)}}{n-1}, \quad \hat{\bm{g}}^{\backslash (i)}_{n} = \Pi \big(\hat{\bm{p}}^{\backslash (i)}_{n} \big).
\end{equation*}

Next, let $\hat{\bm{\phi}}^{\backslash [j]}_{n} \ (j=0, \dots, t_{n})$ denote the leave-one-out version of stacked estimator $\hat{\bm{\phi}}_{n}$ for the frequency data $\bm{x}$, i.e. for $j$ such that $x_{j} > 0$ let
\begin{equation}\label{cvestimj}
\hat{\bm{\phi}}^{\backslash [j]}_{n} =  \beta \, \hat{\bm{g}}^{\backslash [j]}_{n} + (1-\beta)\hat{\bm{p}}^{\backslash [j]}_{n},
\end{equation}
where
\begin{equation*}
\hat{\bm{p}}^{\backslash [j]}_{n} = \frac{\bm{x} - \bm{\sigma}^{[j]}}{n-1}, \quad \hat{\bm{g}}^{\backslash [j]}_{n} = \Pi \big(\hat{\bm{p}}^{\backslash [j]}_{n} \big).
\end{equation*}
Note that $\hat{\bm{g}}^{\backslash (i)}_{n}$ denotes the output of the map $\Pi$ of the empirical estimator based on sample $\mathcal{Z}_{n}$ without $z_{i}$, while $\hat{\bm{g}}^{\backslash[j]}_{n}$ is the output of the map $\Pi$ of the empirical estimator, based on corresponding to the sample $\mathcal{Z}_{n}$ the frequency vector $\bm{x}$ without one count in the cell with index $\bm{s}_{j}$. 

Further, following \cite{bowman1984cross, pastukhov2022stacked, stone1974cross}, we estimate cross-validation criterions in the following way:
\begin{equation}\label{ecvs}
\widehat{CV}(\beta, \mathcal{L})  = 
\begin{cases}
\frac{1}{n}\sum_{i=1}^{n}L_{1}(\bm{\delta}^{(i)}, \hat{\bm{\phi}}^{\backslash (i)}_{n}), &\text{ for }  \mathcal{L} = L_{1}, \\
\frac{1}{n}\sum_{i=1}^{n}L_{2}(\bm{\delta}^{(i)}, \hat{\bm{\phi}}^{\backslash (i)}_{n}), &\text{ for }  \mathcal{L} = L_{2}.\\
\end{cases}
\end{equation}

Next, Stone estimator $\hat{\bm{p}}^{S}_{n}$, introduced in \cite{stone1974cross}, is a mixture of the empirical estimator $\hat{\bm{p}}_{n}$ and some fixed prescribed probability vector $\bm{\lambda}$:
\begin{equation*}\label{}
\hat{\bm{p}}^{S}_{n} = \alpha_{n} \bm{\lambda} + (1-\alpha_{n})\hat{\bm{p}}_{n},
\end{equation*}
where $\alpha_{n}$ is the mixture parameter selected by the leave-one-out cross-validation, which for $\widehat{CV}(\beta, L_{1})$ criterion is given by 
\begin{equation*}\label{}
  \alpha^{(L_{1})}_{n} = \begin{cases}
    0, & \text{if } \sum_{j=0}^{t_{n}}\{ (x_{j} - n \lambda_{j})^{2} + (n+1)\lambda_{j}(x_{j} - n \lambda_{j}) + n\lambda_{j}^{2} \} \geq n,\\
    1, & \text{otherwise},
  \end{cases}
\end{equation*}
and for $\widehat{CV}(\beta, L_{2})$ criterion it is
\begin{equation*}\label{}
  \alpha^{(L_{2})}_{n} = \begin{cases}
    \frac{b^{*}_{n}}{a^{*}_{n}}, & \text{if } \, b^{*}_{n} \leq a^{*}_{n}, \\
    1, & \text{otherwise},
  \end{cases}
\end{equation*}
where 
\begin{eqnarray*}
a^{*}_{n} &=& n - \sum_{j=0}^{t_{n}}\{2x_{j}(x_{j} - \lambda_{j}(n-1)) -  n(x_{j} - \lambda_{j}(n-1))^{2}\},\\
b^{*}_{n} &=& n^2 - \sum_{j=0}^{t_{n}}x_{j}^{2}.
\end{eqnarray*}

From the derivation of $\alpha_{n}$ in \cite{stone1974cross} it follows that $a^{*}_{n}$ is always positive. Also, $b^{*}_{n}$ is nonnegative and  $b^{*}_{n} = 0$ only if $x_{j'}=n$ for some $j'$ and, clearly, $x_{j} = 0$ for $j \neq j'$.

The mixture parameters $\beta_{n}$ for the stacked estimator $\hat{\bm{\phi}}_{n}$ based on the estimated cross-validation criterions in (\ref{ecvs}) are given in the next theorem. 

\begin{theorem}\label{thmbetacv}
The leave-one-out cross-validation selected mixture parameter $\beta$ for $\widehat{CV}(\beta, L_{1})$ criterion is
\begin{equation*}\label{}
  \beta^{(L_{1})}_{n}= 
  \begin{cases}
    0, & \text{if }  \sum_{j=0}^{t_{n}}\{ (x_{j} - n \hat{g}^{\backslash[j]}_{n, \bm{s}_{j}})^{2} + (n+1)\hat{g}^{\backslash[j]}_{n, \bm{s}_{j}}(x_{j} - n \hat{g}^{\backslash[j]}_{n, \bm{s}_{j}}) + n(\hat{g}^{\backslash[j]}_{n, \bm{s}_{j}})^{2} \} \geq n, \\
    1, & \text{otherwise},
  \end{cases}
\end{equation*}
for $\widehat{CV}(\beta, L_{2})$ criterion:
\begin{equation*}\label{}
  \beta^{(L_{2})}_{n} = \begin{cases}
    \frac{b_{n}}{a_{n}}, & \text{if } \, a_{n} \neq 0 \text{ and } 0 \leq b_{n} \leq a_{n}, \\
    0, & \text{if } \, a_{n} = 0 \text{ or } b_{n} < 0,\\
    1, & \text{otherwise},
  \end{cases}
\end{equation*}
where 
\begin{equation*}
a_{n} = n - \sum_{j=0}^{t_{n}}x_{j}\big\{2(x_{j} - \hat{g}^{\backslash[j]}_{n, \bm{s}_{j}}(n-1)) - \sum_{k=0}^{t_{n}}(x_{k} - (n-1)\hat{g}^{\backslash[j]}_{n, \bm{s}_{k}})^{2} \big\},
\end{equation*}
\begin{equation*}
b_{n} = n^2 - \sum_{j=0}^{t_{n}}x_{j}^{2} + (n-1) \sum_{j=0}^{t_{n}} x_{j} \sum_{k=0}^{t_{n}} x_{k} ( \hat{g}^{\backslash[k]}_{n, \bm{s}_{k}} - \hat{g}^{\backslash[j]}_{n, \bm{s}_{k}}),
\end{equation*}
\end{theorem}
Note that if $\hat{g}^{\backslash[j]}_{n, \bm{s}_{k}} = \lambda_{k}$, i.e. there is no dependence on $j$, then $a_{n} = a^{*}_{n}$ and $b_{n} = b^{*}_{n}$, and, therefore, in this case we have $\beta^{(L_{1})}_{n} = \alpha^{(L_{1})}_{n}$ and $\beta^{(L_{2})}_{n} = \alpha^{(L_{2})}_{n}$.

\section{Consistency and rate of convergence of Grenander-Stone estimator}\label{constest}
In this section we derive sufficient conditions for the stacked estimator in (\ref{defest}) to be consistent with parametric rate of convergence and prove that it holds for the case of Grenander-Stone estimator for both choices of cross-validation criterions. We make the following assumption on the map $\Pi$. 

\begin{assumption}\label{asmt}
Assume that  the map $\Pi$ in (\ref{defest}) has the following properties:
\begin{enumerate}[label=(\roman*)]
\item It is continues from $\ell_{2}$ to $\ell_{2}$.
\item It is contraction with respect to $\ell_{\infty}$-norm.
\item $\langle \Pi(\bm{f}), (\bm{f} -  \Pi(\bm{f})) \rangle = 0$, and $\langle \bm{f}, (\bm{f} -  \Pi(\bm{f})) \rangle \geq 0$ for any $\bm{f} \in  \mathcal{P}$.
\item There exists a sequence $c_{n} \in \mathbb{R}$, $n = 1, 2, \dots$, such that 
\begin{equation*}\label{}
c_{n} \to 1,
\end{equation*}\label{}
and for any frequency data $\bm{x}$:
\begin{equation*}\label{}
c_{n} \bm{\gamma}_{n} \leq \Pi(\hat{\bm{p}}_{n}),
\end{equation*} 
where 
\begin{equation*}\label{}
\begin{aligned}
  \hat{\gamma}_{n,j}=\begin{cases}
    \Pi(\hat{\bm{p}}^{\backslash [j]}_{n}), & \text{if } \, x_{j} > 0, \\
    0, & \text{otherwise},
      \end{cases} 
\end{aligned}
\end{equation*}
and the notation $\bm{a} \leq \bm{b}$ for some vectors $\bm{a},\bm{b} \in \mathbb{R}^{m}$ means $a_{i} \leq b_{i}$ for all $i = 1, \dots, m$.
\item In the case when the true distribution is such that $\Pi(\bm{p}) = \bm{p}$, then
\begin{equation*}
n^{1/2}||\Pi(\hat{\bm{p}}_{n}) - \bm{p}||_{k} = O_{p}(1), \quad (2 \leq k \leq \infty),
\end{equation*}
and if in addition $\sum_{j=0}^{\infty}\sqrt{p}_{j} < \infty$, then for both $L_{1}$ and $L_{2}$ loss functions we have:
\begin{equation*}
n^{1/2}||\Pi(\hat{\bm{p}}_{n}) - \bm{p}||_{1} = O_{p}(1).
\end{equation*}
\end{enumerate}
\end{assumption}

\begin{theorem}\label{consrtgnr}
Assume that Assumption \ref{asmt} holds. Then, for any underlying distribution $\bm{p}$ (not necessarily for $\bm{p}$ such that $\Pi(\bm{p}) = \bm{p}$) and for both $L_{1}$ and $L_{2}$ loss functions for leave-one-out cross validation step the stacked estimator $\hat{\bm{\phi}}_{n}$ defined in (\ref{defest}) is strongly consistent in $\ell_{k}$-norm $(1 \leq k \leq \infty)$
\begin{equation*}
\hat{\bm{\phi}}_{n} \stackrel{a.s.}{\to} \bm{p},
\end{equation*}
and it has $n^{1/2}$-rate of convergence
\begin{equation*}
n^{1/2}||\hat{\bm{\phi}}_{n} - \bm{p}||_{k} = O_{p}(1), \quad (2 \leq k \leq \infty).
\end{equation*}
Furthermore, if in addition $\sum_{j=0}^{\infty}\sqrt{p}_{j} < \infty$, then for both $L_{1}$ and $L_{2}$ loss functions we have:
\begin{equation*}
n^{1/2}||\hat{\bm{\phi}}_{n} - \bm{p}||_{1} = O_{p}(1).
\end{equation*}
\end{theorem}

In Theorem \ref{consrtgnr} above we proved strong consistency of Grenander--Stone estimator for both $L_{1}$ and $L_{2}$ loss functions. Further, the $L_{2}$ loss is known to be a strictly proper loss function, as proved, for example in \cite{czado2009predictive, haghtalab2019toward}, which means that  for the scoring rule $S_{2}(\bm{f}, \bm{p})$ we have 
\begin{equation*}
\bm{p} = \underset{\bm{f} \in \mathcal{P}}{\arg \min} \, S_{2}(\bm{f}, \bm{p}).
\end{equation*} 
Therefore, for the case of penalization by $L_{2}$ loss the consistency in Theorem \ref{consrtgnr} is the expected result. Next, compare to the case of $L_{2}$ loss,  the $L_{1}$ loss function is not even a proper loss function, meaning that, in general, $S_{1}(\bm{f}, \bm{p})$ is not minimized by $\bm{f} = \bm{p}$, cf. \cite{buja2005loss}.
Nevertheless, below we show that in the certain family of distributions the conditional expected $L_{1}$ loss is minimised by the underlying distribution $\bm{p}$.
\begin{lemma}\label{propL1}
Let $\mathcal{M}_{p}$ denote the following family of discrete distributions
\begin{equation*}
\mathcal{M}_{p} = \{ \bm{f} \in \mathcal{M}_{p}: \bm{f} = \beta \bm{h} + (1-\beta)\bm{p} \, \text{ for some } \, \bm{h} \in \mathcal{P} \, \text{ and } \, 0 \leq \beta \leq 1 \}
\end{equation*} 
for some $\bm{p} \in \mathcal{P}$. Then
\begin{equation*}
\bm{p} = \underset{\bm{f} \in\mathcal{M}_{p}}{\arg \min} \, S_{1}(\bm{f}, \bm{p}).
\end{equation*} 
\end{lemma}
Therefore, the result of Proposition \ref{propL1} explains consistency provided by penalisation with $L_{1}$ loss even though, in general, $L_{1}$ loss is not a proper loss function. Moreover, as shown in the proof of Theorem \ref{consrtgnr}, if $\Pi(\bm{p}) \neq \bm{p}$ for the estimator $\hat{\bm{\phi}}_{n}$ in the case of $L_{1}$ loss function the following holds 
\begin{equation*}
\mathbb{P}( \underset{n\to\infty}{ \lim\inf} \{ \hat{\bm{\phi}}_{n} = \hat{\bm{p}}_{n}\}) = 1.
\end{equation*}

Now, let us assume that $\Pi(\hat{\bm{p}}_{n}) = \Pi^{is}(\hat{\bm{p}}_{n}  \mid \mathcal{P}^{is})$ and consider Grenander-Stone estimator:
\begin{equation}\label{defestGS}
\hat{\bm{\phi}}_{n} = \hat{\beta}_{n}  \Pi^{is}(\hat{\bm{p}}_{n}  \mid \mathcal{P}^{is}) + (1-\hat{\beta}_{n})\hat{\bm{p}}_{n}.
\end{equation}

\begin{theorem}\label{GScons}
For any underlying distribution $\bm{p}$ (not necessarily isotonic) and for both $L_{1}$ and $L_{2}$ loss functions for leave-one-out cross validation step the statements of Assumption \ref{asmt} for Grenander-Stone estimator hold. Therefore, the estimator in (\ref{defestGS}) is strongly consistent in $\ell_{k}$-norm $(1 \leq k \leq \infty)$
\begin{equation*}
\hat{\bm{\phi}}_{n} \stackrel{a.s.}{\to} \bm{p},
\end{equation*}
and it has $n^{1/2}$-rate of convergence
\begin{equation*}
n^{1/2}||\hat{\bm{\phi}}_{n} - \bm{p}||_{k} = O_{p}(1), \quad (2 \leq k \leq \infty).
\end{equation*}
Furthermore, if in addition $\sum_{j=0}^{\infty}\sqrt{p}_{j} < \infty$, then for both $L_{1}$ and $L_{2}$ loss functions we have:
\begin{equation*}
n^{1/2}||\hat{\bm{\phi}}_{n} - \bm{p}||_{1} = O_{p}(1).
\end{equation*}
\end{theorem}

Finally, for Grenander--Stone estimator in one-dimensional case we have the following Marshal-type result (for the case of estimation of p.m.f. with a convex constraint cf. \cite{balabdaoui2015marshall, dumbgen2007marshall}).
\begin{theorem}\label{Marshal}
Let $\hat{F}_{n, j} = \sum_{k=0}^{j}\hat{p}_{n, j} \ (j \in \mathbb{N})$ be the empirical distribution function, and $\hat{F}_{n, j}^{GS} = \sum_{k=0}^{j}\hat{\phi}_{n, j} \ (j \in \mathbb{N})$. Then for any underlying distribution $\bm{p}$ the following holds
\begin{equation*}\label{}
||\hat{F}_{n}^{GS}  - D||_{\infty} \leq ||\hat{F}_{n}  - D||_{\infty},
\end{equation*}
where $D_{j} \ (j \in \mathbb{N})$ is a cumulative distribution function of any decreasing distribution on $\mathbb{N}$.
\end{theorem}

\section{Global confidence band}\label{confband}
In this section we construct the asymptotic global confidence band for the true distribution $\bm{p}$. First, assume that the statements of Assumption \ref{asmt} hold. Second, let $\bm{Y}_{\bm{p}}$ be a Gaussian process in $\ell_{2}$ with mean zero and the covariance operator $C$ such that 
\begin{equation*}\label{}
\langle C \bm{e}_{i}, \bm{e}_{i'}\rangle = p_{i}\delta_{i,i'} - p_{i}p_{i'},
\end{equation*} 
with $\bm{e}_{i} \ (i=0,1, \dots)$ the orthonormal basis in $\ell_{2}$, meaning that in a vector $\bm{e}_{i}$ all elements are equal to zero but the one with the index $i$ is equal to $1$, and $\delta_{i,j}=1$ is the Kronecker delta. 

Further, as shown in \cite{jankowski2009estimation}, the process $\bm{Y}_{\bm{p}}$ is the asymptotic distribution of the empirical estimator. Here we emphasise that the asymptotic distribution of $\hat{\bm{p}}_{n}$ in the case of a finite support is a straightforward result. Nevertheless, to the authors' knowledge, \cite{jankowski2009estimation} is the first work where the limiting distribution of the empirical estimator was obtained for the case of an infinite support.

\begin{figure}[!htbp] 
  \begin{subfigure}{3.9cm}
    \centering\includegraphics[scale=0.22]{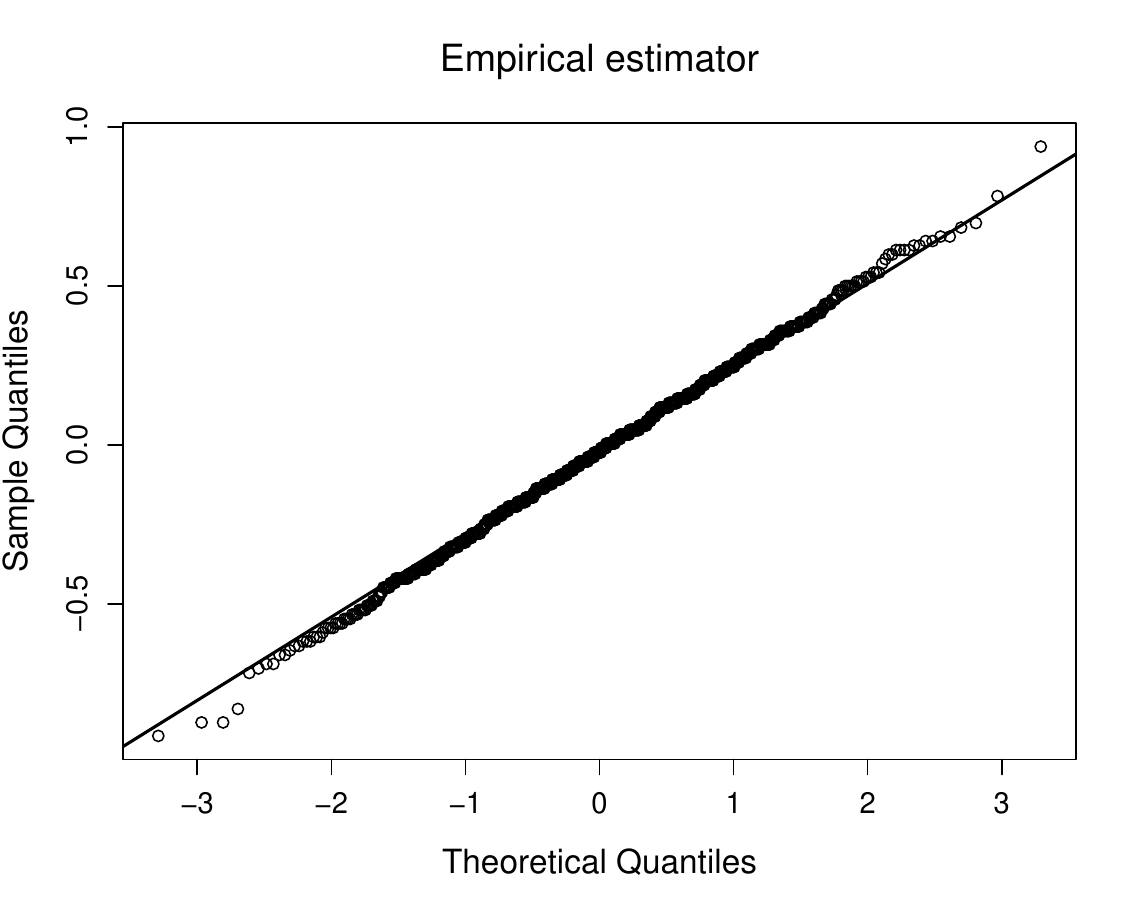}
  \end{subfigure}
  \begin{subfigure}{3.9cm}
    \centering\includegraphics[scale=0.22]{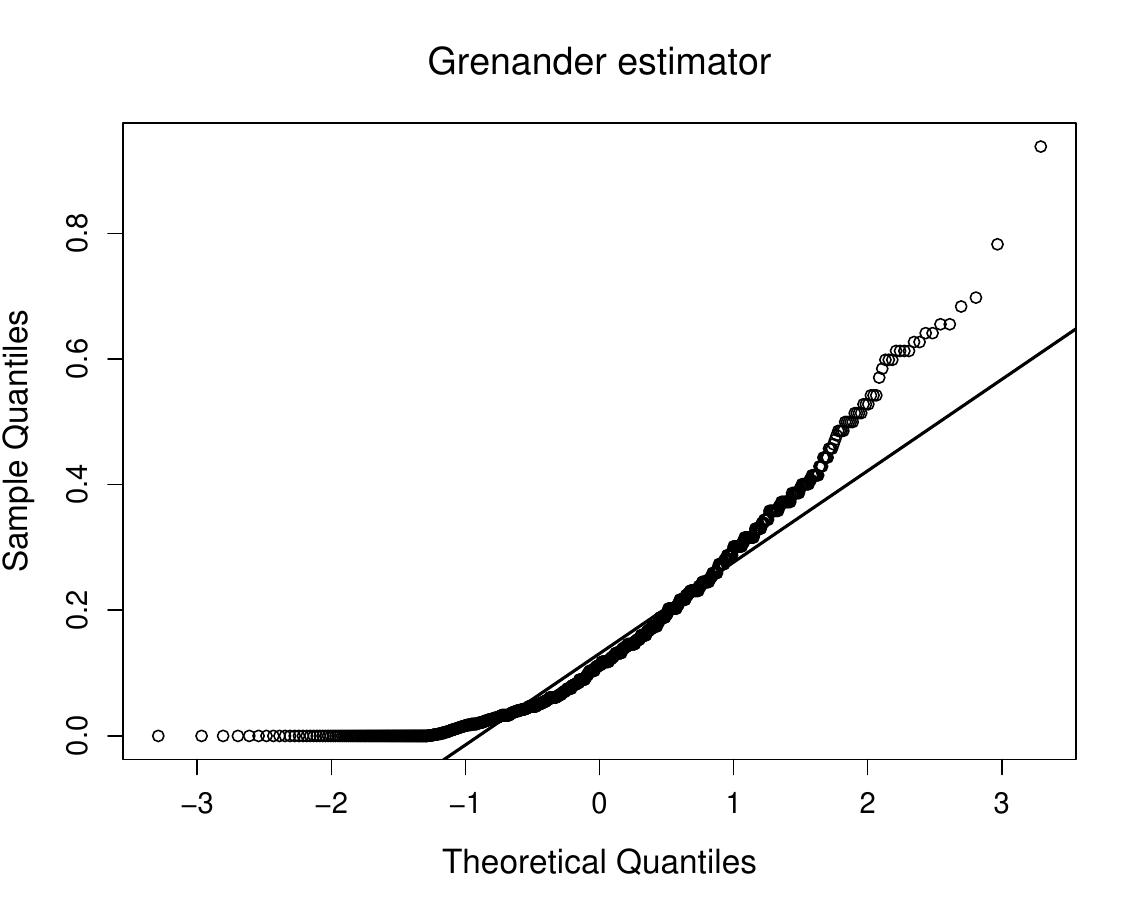}
  \end{subfigure}
    \begin{subfigure}{3.9cm}
    \centering\includegraphics[scale=0.22]{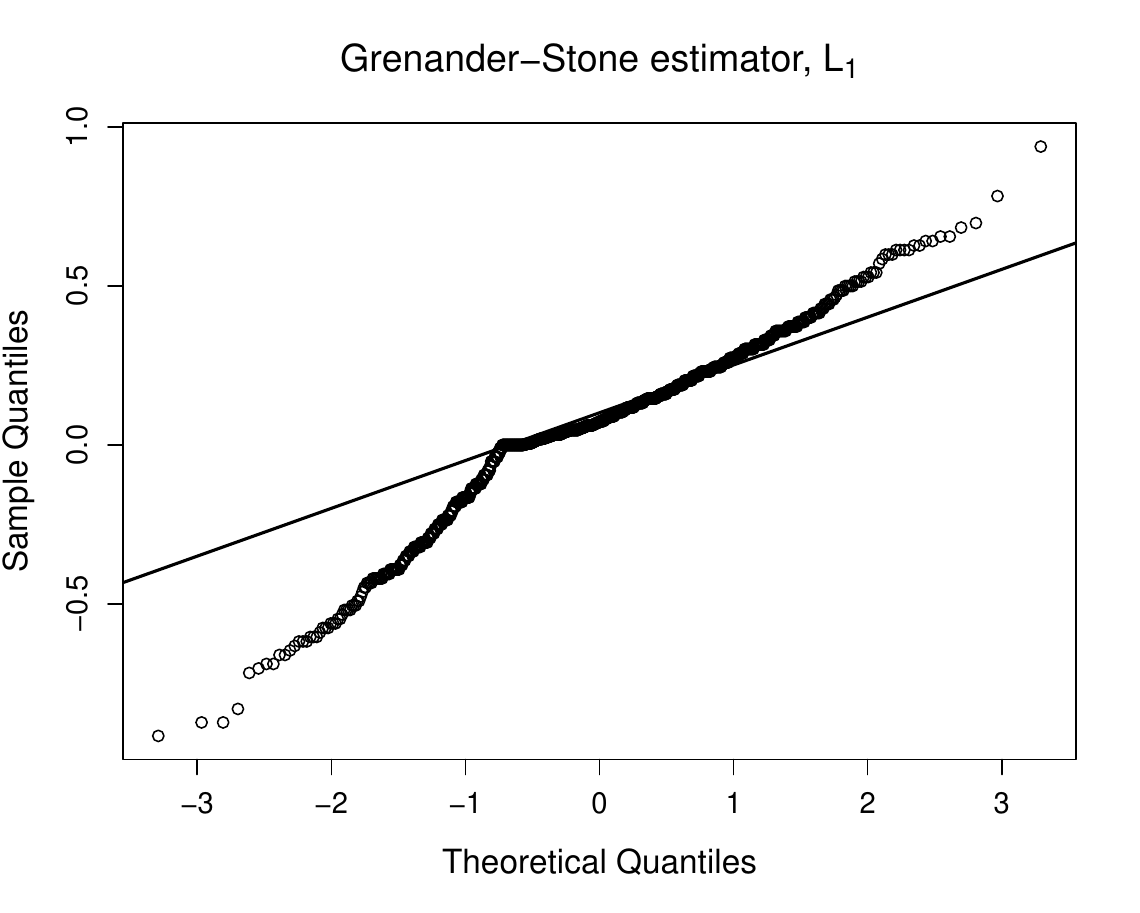}
  \end{subfigure}
  \begin{subfigure}{3.9cm}
    \centering\includegraphics[scale=0.22]{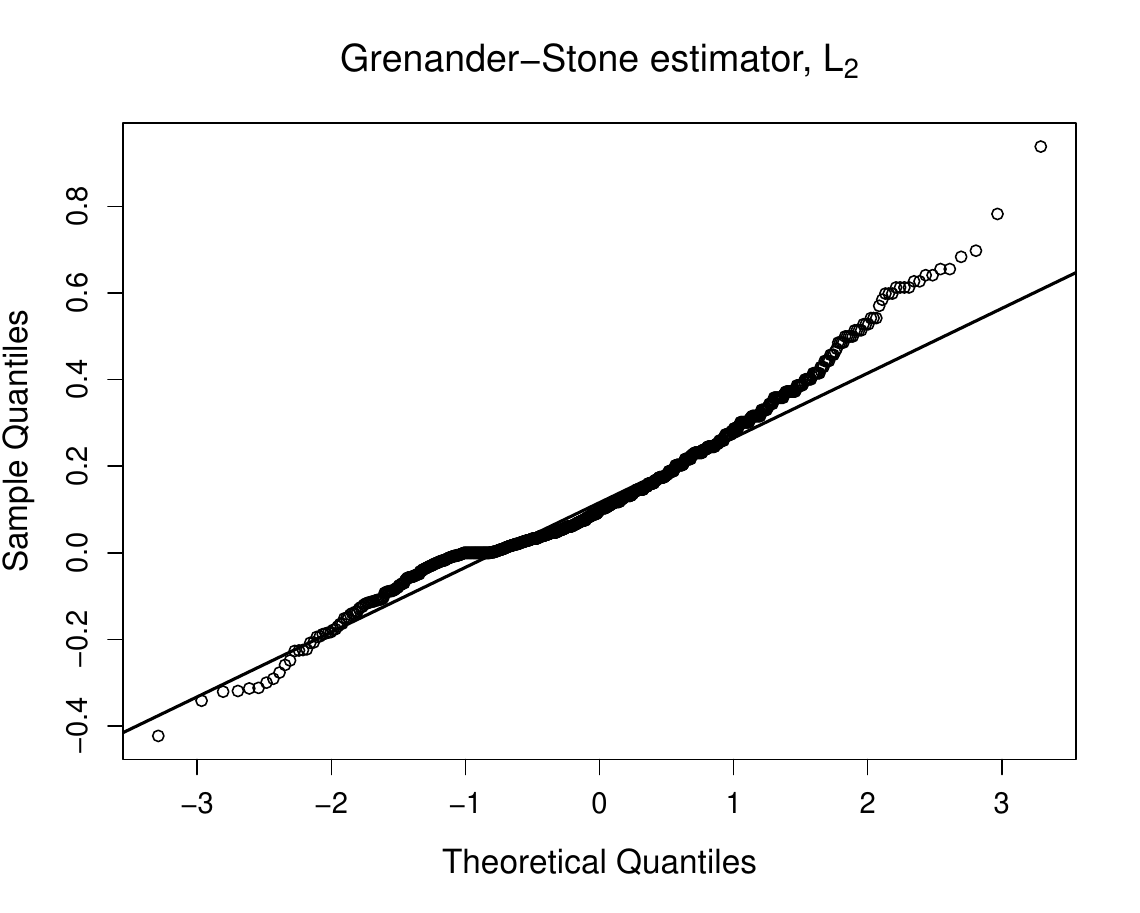}
  \end{subfigure}
   \caption{Standard normal QQ-plots of 1000 samples of the empirical estimator $\sqrt{n}(\hat{p}_{n,1} - p_{1})$, Grenander estimator $\sqrt{n}(\hat{g}_{n,1} - p_{1})$ and Grenander-Stone estimator $\sqrt{n}(\hat{\phi}_{n,1} - p_{1})$ for both $L_{1}$ and $L_{2}$, with $n=1000$ samples for uniform distribution over $\mathcal{S}=\{0,\dots 11\}$.}\label{qq-plots}
\end{figure}

Next, if the underlying distribution $\bm{p}$ is such that $\Pi(\bm{p}) \neq \bm{p}$, then from the proof of Theorem \ref{consrtgnr} it follows that $n^{1/2} \beta^{(L_{d})}_{n} \stackrel{a.s.}{\to} 0$, for $d= 1, 2$. Therefore, in the similar way as for the cases of stacked Grenander and rearrangement estimators in Theorem 5 in \cite{pastukhov2022stacked}, for stacked estimator $\hat{\bm{\phi}}_{n}$ one can show that 
\begin{equation*}
n^{1/2}(\hat{\bm{\phi}}_{n}-\bm{p}) \stackrel{d}{\to} \bm{Y}_{\bm{p}},
\end{equation*} 
for both loss functions $L_{1}$ and $L_{2}$. The asymptotic distribution of stacked estimator in the case, when $\Pi(\bm{p}) = \bm{p}$, is an open problem. On Figure \ref{qq-plots} we illustrate the difference of the asymptotic distributions of the empirical estimator, Grenander estimator and the stacked estimators by plotting standard normal QQ-plots of 1000 samples of the estimators for $n=1000$ samples for uniform distribution over the support $\mathcal{S}=\{0,\dots 11\}$. We can conclude that empirical estimator, stacked estimators and Grenander estimator have different asymptotic distribution.
 
Further, for the process $\bm{Y}_{\bm{p}}$ let $q_{\alpha}$ denote the $\alpha$-quantile of its $\ell_{\infty}$-norm, i.e. 
\begin{equation*}
\mathbb{P}(||\bm{Y}_{\bm{p}}||_{\infty} > q_{\alpha}) = \alpha,
\end{equation*} 
and, analogously, for the stacked estimator let $\hat{q}^{d}_{\alpha, n}$, for $d = 1, 2$, denote the $\alpha$-quantiles for the processes $\bm{Y}_{\hat{\bm{\phi}}_{n}}$, i.e. 
\begin{equation*}
\mathbb{P}(||\bm{Y}_{\hat{\bm{\phi}}_{n}}||_{\infty} > \hat{q}^{d}_{\alpha, n}) = \alpha,
\end{equation*} 
for $d = 1, 2$.

Next, consequently, in the case when the underlying distribution $\bm{p}$ is such that $\Pi(\bm{p}) \neq \bm{p}$ (i.e. $\bm{p}$ is not isotonic when $\hat{\bm{\phi}}_{n}$ is Grenander-Stone estimator) we have
\begin{equation*}
\lim_{n}\mathbb{P}(n^{1/2}||\hat{\bm{\phi}}_{n} - \bm{p}||_{\infty} \leq q^{d}_{\alpha}) =  1- \alpha,
\end{equation*} 
with $d=1,2$ for both penalizing loss functions $L_{1}$ and $L_{2}$, respectively. 

In the case when $\Pi(\bm{p}) = \bm{p}$ (i.e. $\bm{p}$ is isotonic when $\hat{\bm{\phi}}_{n}$ is Grenander-Stone estimator) from the subadditivity of norms and statement $(v)$ of Assumption \ref{asmt} we have 
\begin{equation*}\label{}
\begin{aligned}
||\hat{\bm{\phi}}_{n} - \bm{p}||_{k} {}  \leq \hat{\beta}^{(L_{d})}_{n}||\hat{\bm{g}}_{n} - \bm{p}||_{k} + (1-\hat{\beta}^{(L_{d})}_{n})||\hat{\bm{p}}_{n} - \bm{p}||_{k} \leq ||\hat{\bm{p}}_{n} - \bm{p}||_{k},
\end{aligned}
\end{equation*} 
for $1 \leq k \leq \infty$ and $d = 1, 2$. Therefore, in this case for the stacked estimator $\hat{\bm{\phi}}_{n} $ the following holds 

\begin{equation*}
\liminf_{n}\mathbb{P}(n^{1/2}||\hat{\bm{\phi}}_{n} - \bm{p}||_{\infty} \leq q^{d}_{\alpha}) \geq 1-  \alpha,
\end{equation*} 
for $d = 1, 2$.

In the next theorem, analogously to the unimodal density estimator in \cite{balabdaoui2016maximum}, we prove that $\hat{q}^{d}_{\alpha, n}$, for $d = 1, 2$, converges to $q_{\alpha}$.
\begin{theorem}\label{convqntl}
The $\alpha$-quantiles $\hat{q}^{d}_{\alpha, n}$ of $\ell_{\infty}$-norms of $\bm{Y}_{\hat{\bm{\phi}}_{n}}$, for $d = 1, 2$, are a strongly consistent estimators of 
the $\alpha$-quantile of $\ell_{\infty}$-norm of $\bm{Y}_{\bm{p}}$.
\end{theorem}
 
Finally, the confidence band
\begin{equation*}\label{}
\Big[\max\Big\{(\hat{\phi}_{n,\bm{s}_{j}}  - \frac{\hat{q}^{d}_{\alpha, n}}{n^{1/2}}), 0\Big\}, \,  \hat{\phi}_{n, \bm{s}_{j}} + \frac{\hat{q}^{d}_{\alpha, n}}{n^{1/2}}\Big] \quad (j\in\mathbb{N})
\end{equation*} 
is asymptotically correct global confidence band if $\bm{p}$ is not decreasing and it is asymptotically correct conservative global confidence band if $\bm{p}$ is decreasing for both choices of penalizing loss functions $L_{1}$ and $L_{2}$. We can use Monte-Carlo method, as in \cite{balabdaoui2016maximum}, to estimate the quantiles $\hat{q}^{d}_{\alpha, n}$ of the Gaussian vector $\bm{Y}_{\hat{\bm{\phi}}_{n}}$, $d = 1, 2$.

\section{Computational aspects and simulation study}\label{est_sim}
In this section we compare the performance of Grenander--Stone estimator with the empirical and Grenander estimators and consider both monotone and non-monotone underlying distributions.

In the case of one-dimensional distribution, Grenander estimator, i.e. the optimization problem
\begin{eqnarray*}
\hat{\bm{g}}_{n} = \Pi^{is}(\hat{\bm{p}}_{n}  \mid \mathcal{P}^{is}) =  \underset{f_{1} \geq \dots \geq f_{t_{n}}}{\arg \min}\sum_{j=0}^{t_{n}}(\hat{p}_{n, \bm{s}_{j}} - f_{\bm{s}_{j}})^{2}
\end{eqnarray*}
can be exactly solved with $O(t_{n})$ computational complexity using Pool-Adjacent-Violators algorithm (cf.  \cite{brunk1972statistical}). Further, in the case of the distribution over a general directed acyclic graph, i.e. for a graph $G = (V, E)$ the problem
\begin{eqnarray*}
\begin{aligned}
&\hat{\bm{g}}_{n} = \arg \min \sum_{\bm{s} \in \mathcal{S}}(\hat{p}_{n, \bm{s}} - f_{\bm{s}})^{2},\\
&\text{s.t.} \quad  f_{\bm{s}_{i}} \leq f_{\bm{s}_{j}} \quad  \text{for any} \quad (\bm{s}_{i}, \bm{s}_{i}) \in E
\end{aligned}
\end{eqnarray*}
can be solved with complexity $O(|V|^{2})$ by the algorithm proposed in \cite{burdakov2006n}. In the case when, for example, $G$ is a directed tree, the problem can be solved with complexity $O(|V|\log |V|)$, cf. \cite{pardalos1999algorithms}. The approximate fast numerical solution for a general graph can be obtained by solving a dual problem to the nearly-isotonic regression problem with a big penalisation parameter, cf. \cite{pastukhov2023fused}. 

Next, we compare the performance of stacked estimators with the empirical, Grenander, and minimax estimators. For the distribution with a finite support, for example, $\mathcal{S} = \{0, \dots, s \}$, and the sample size $n$ the minimax estimator of $\bm{p}$ with respect to $\ell_{2}$-loss is given by
\begin{equation}\label{estMM}
\hat{\bm{p}}^{mm}_{n} = \alpha^{mm}_{n} \bm{\lambda} + (1-\alpha^{mm}_{n})\hat{\bm{p}}_{n},
\end{equation}
with $\bm{\lambda} = (\frac{1}{s+1}, \dots, \frac{1}{s+1})$ and $\alpha^{mm}_{n} = \frac{\sqrt{n}}{n + \sqrt{n}}$, cf. \cite{trybula1958some}. To the authors' knowledge, the minimax estimation with respect to $\ell_{2}$-loss for infinitely supported discrete distribution is an open problem.  With abuse of notation, for infinitely supported distributions we still refer the estimator defined in (\ref{estMM}) with $s = t_{n}$ as minimax.

\subsection{Finite sample performance}
\subsubsection{True distribution is isotonic}
We consider the following isotonic probability mass functions to assess the performance of Grenander-Stone estimator:
\begin{eqnarray*}
\textbf{M1}: \bm{p} &=& 0.15 U(3) + 0.1 U(7) + 0.75 U(11),\\
\textbf{M2}: \bm{p} &=& Geom(0.25), \\
\textbf{M3}: \bm{p} &=& q_{1}U_{2d}(1) + \dots + q_{5}U_{2d}(5),\\
\end{eqnarray*}
where  $\bm{q} = (0.1, 0.2, 0.3, 0.2, 0.2)$, $U(k)$ denotes the uniform distribution over $\{0, \dots, k\}$,  $U_{2d}(k)$ denotes the uniform distribution over grid $\{0, \dots, k\}\times \{0, \dots, k\}$,  and $Geom(\phi)$ is Geometric distribution, i.e. $p_{j} = (1-\phi)\phi^{j} \ (j \in \mathbb{N})$, with $ 0< \phi < 1$.

\begin{figure}[!htbp] 
  \begin{subfigure}{5.2cm}
    \centering\includegraphics[scale=0.17]{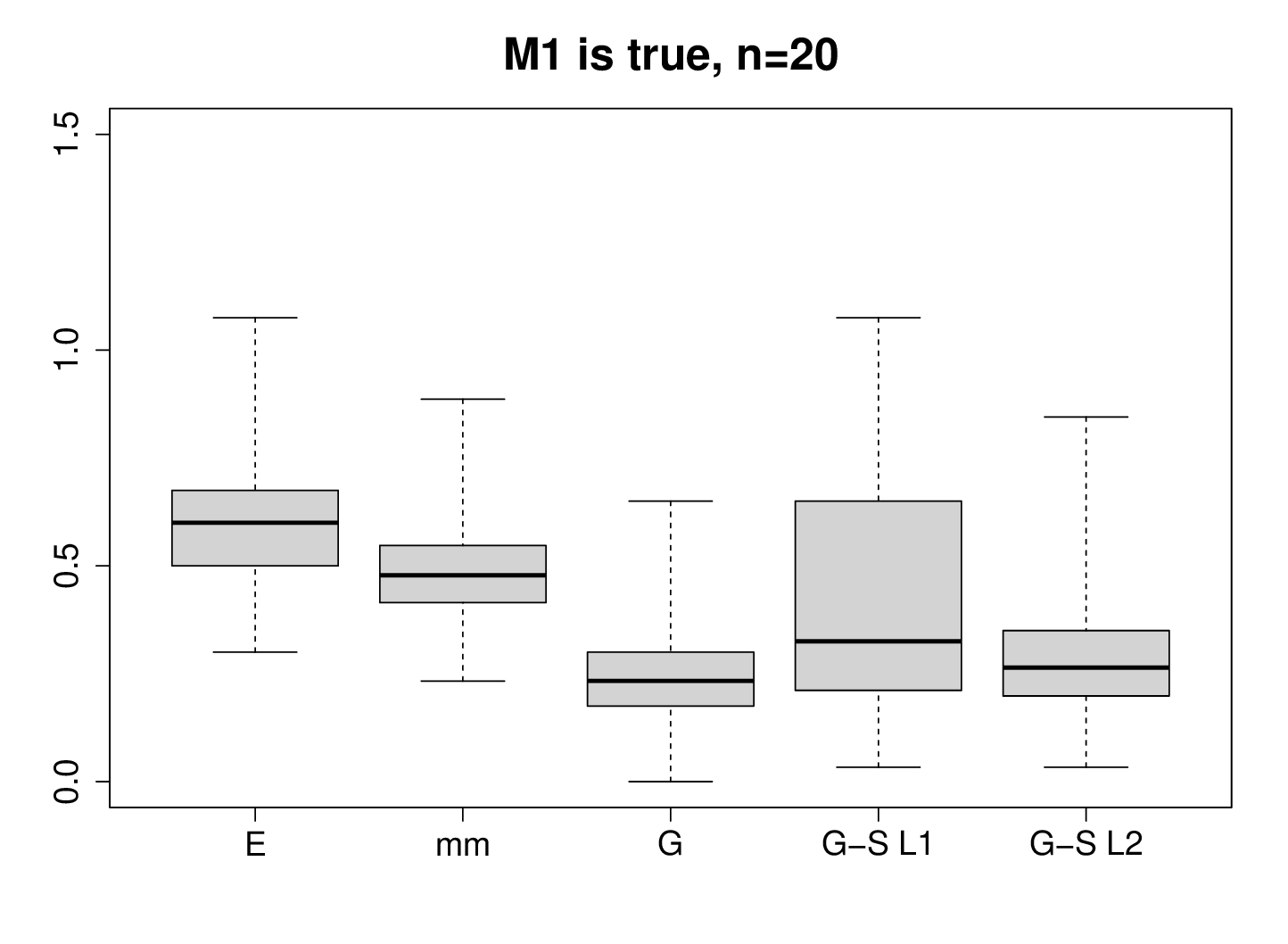}
  \end{subfigure}
  \begin{subfigure}{5.2cm}
    \centering\includegraphics[scale=0.17]{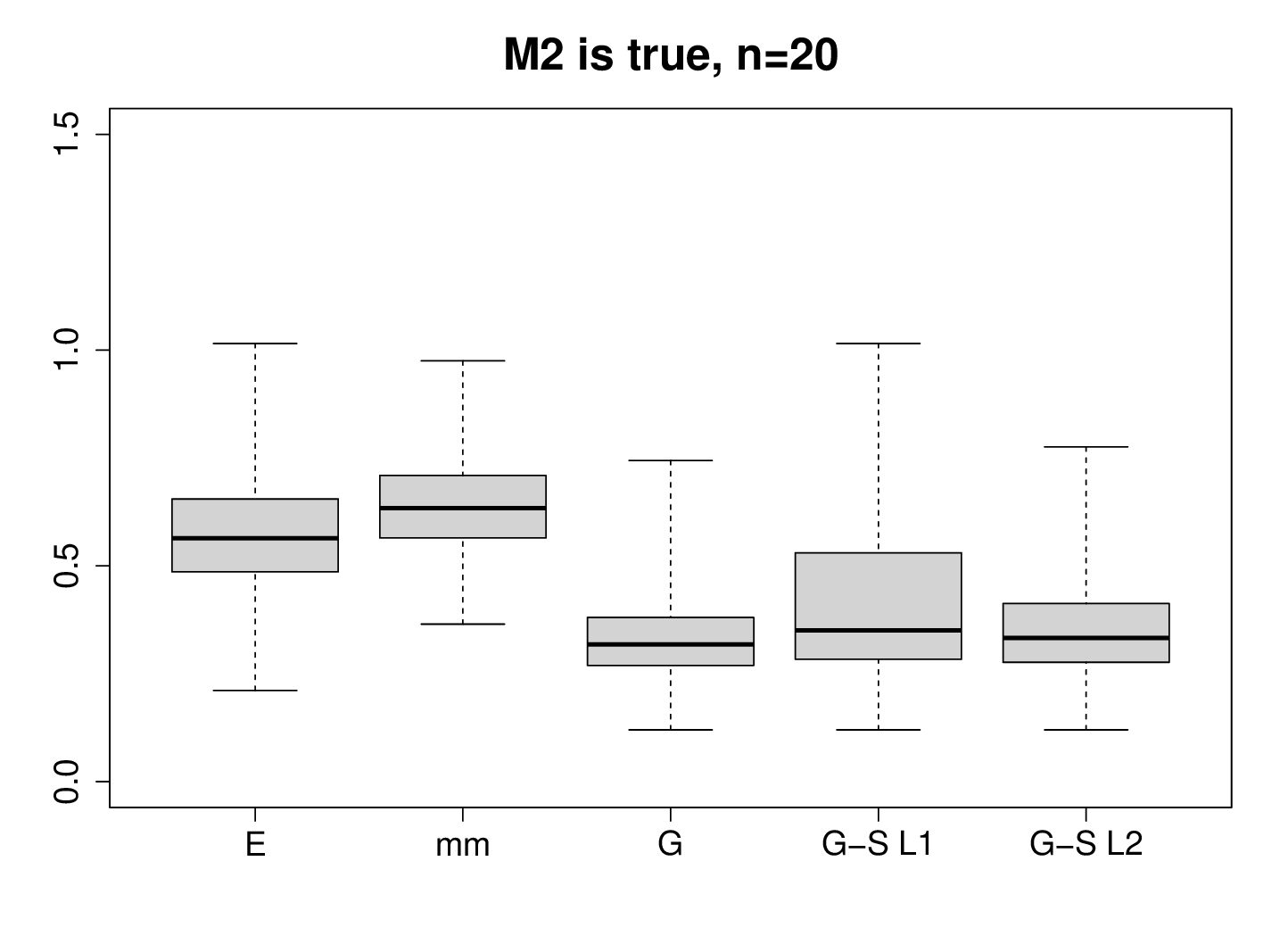}
  \end{subfigure}
  \begin{subfigure}{5.2cm}
    \centering\includegraphics[scale=0.17]{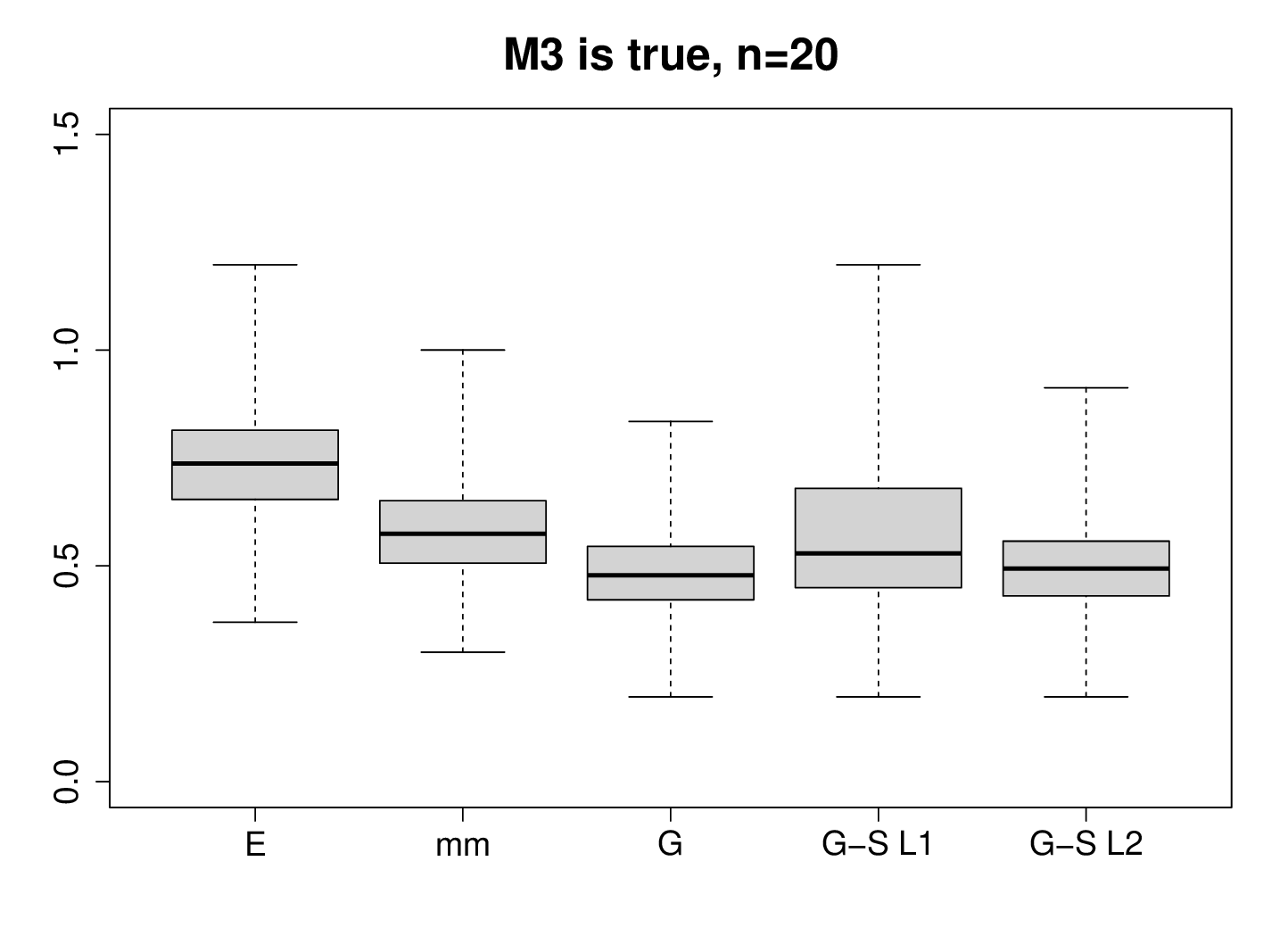}
  \end{subfigure}
  
    \begin{subfigure}{5.2cm}
    \centering\includegraphics[scale=0.17]{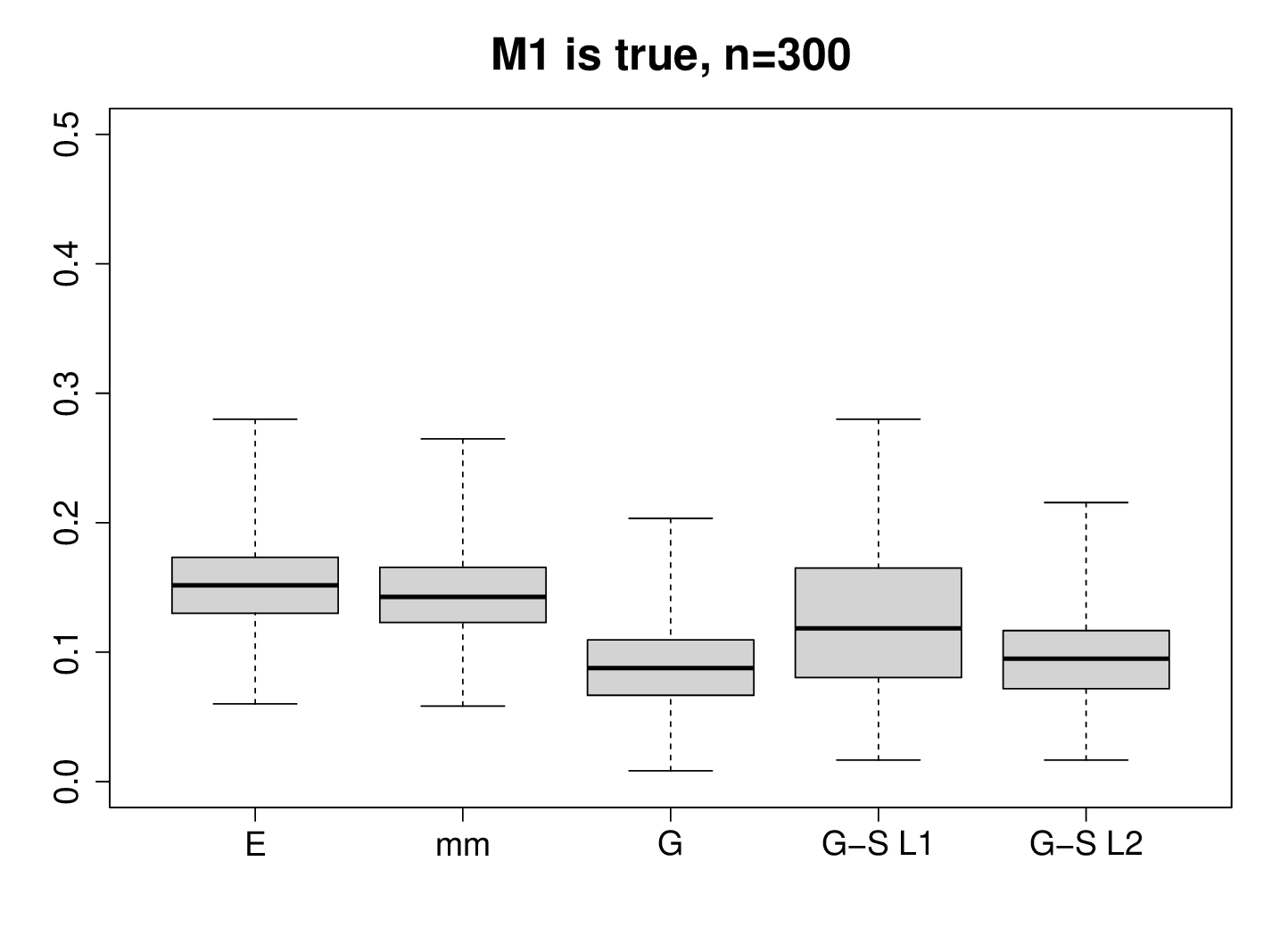}
  \end{subfigure}
  \begin{subfigure}{5.2cm}
    \centering\includegraphics[scale=0.17]{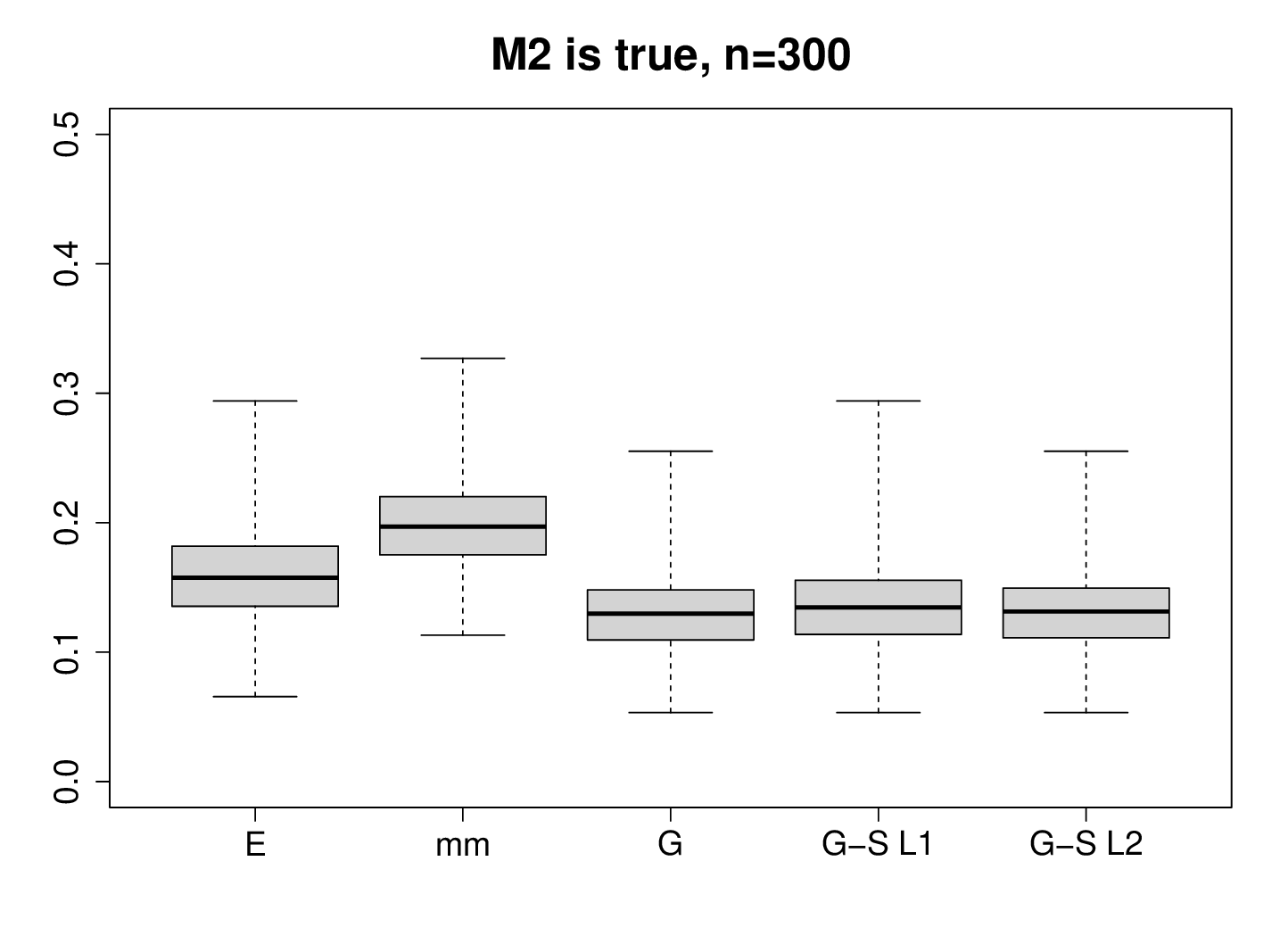}
  \end{subfigure}
  \begin{subfigure}{5.2cm}
    \centering\includegraphics[scale=0.17]{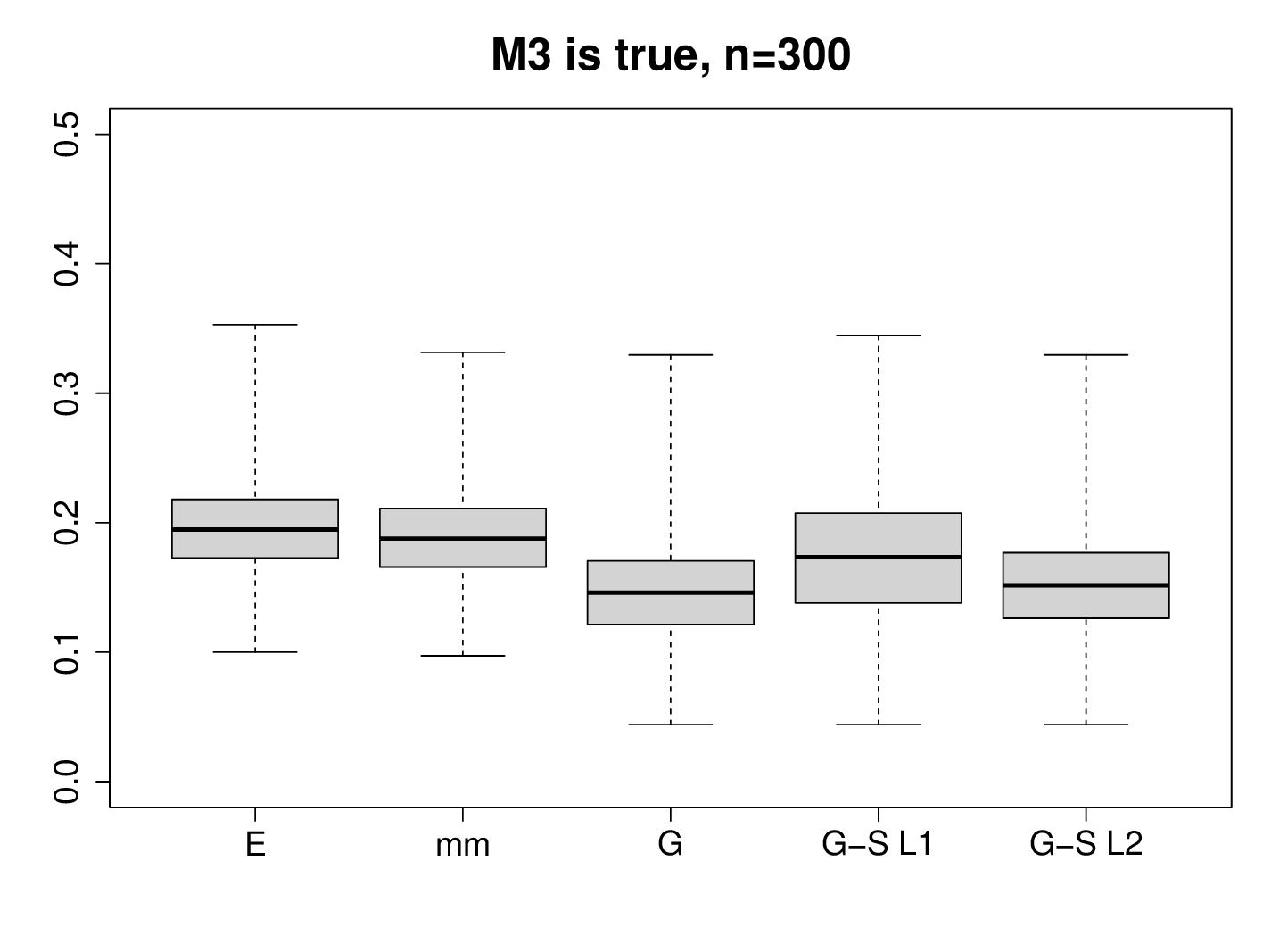}
  \end{subfigure}
  
    \begin{subfigure}{5.2cm}
    \centering\includegraphics[scale=0.17]{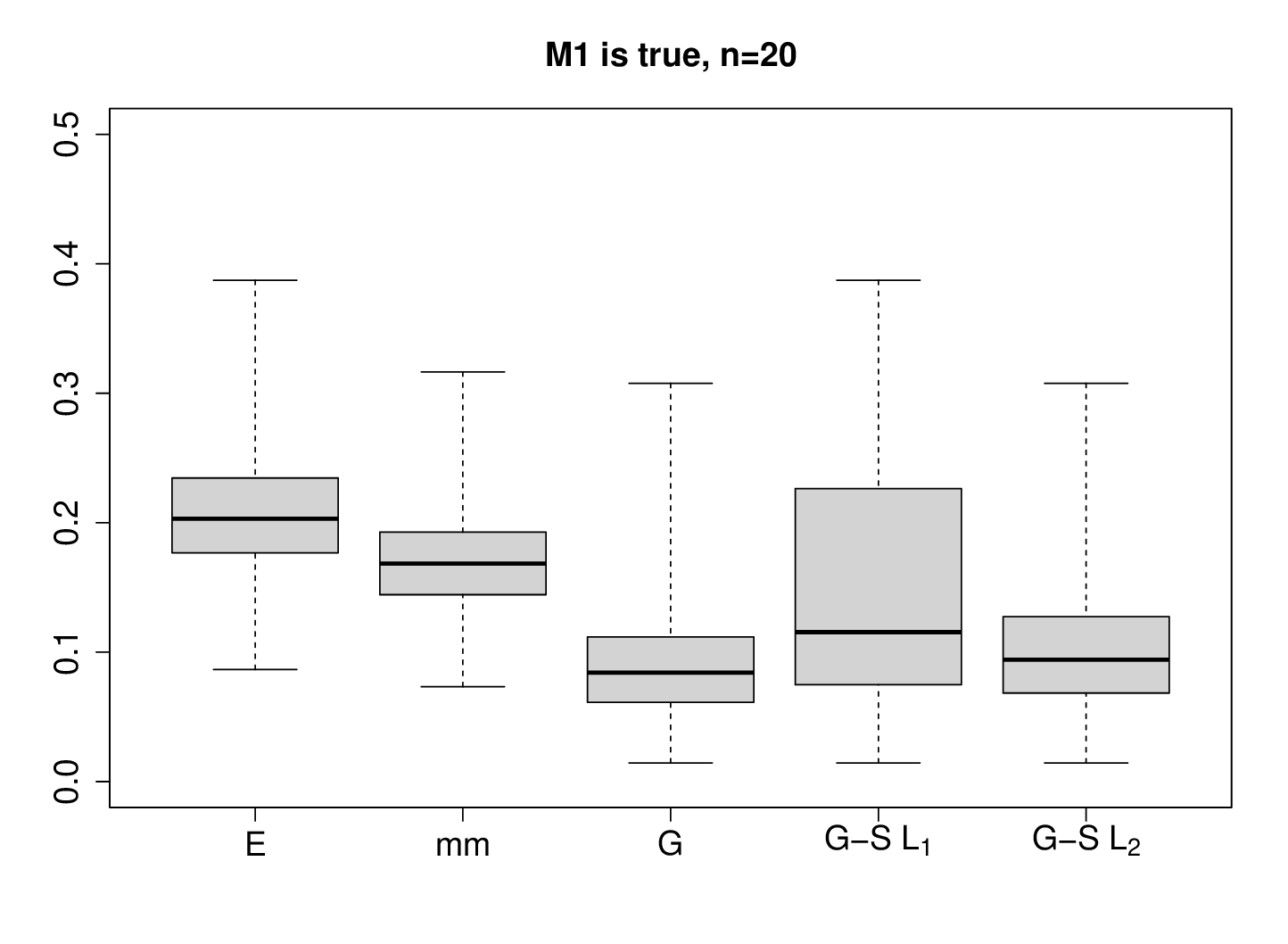}
  \end{subfigure}
  \begin{subfigure}{5.2cm}
    \centering\includegraphics[scale=0.17]{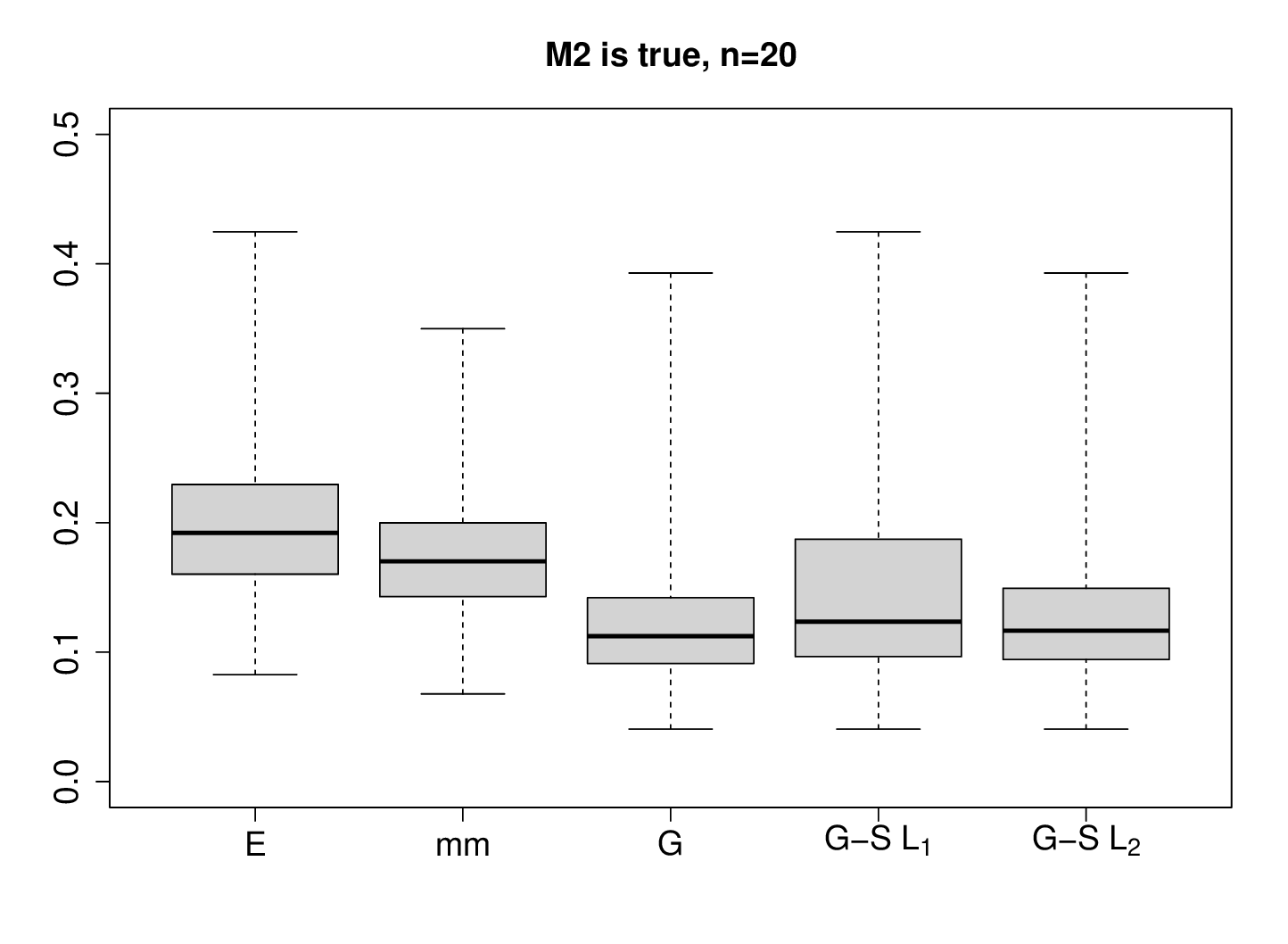}
  \end{subfigure}
  \begin{subfigure}{5.2cm}
    \centering\includegraphics[scale=0.17]{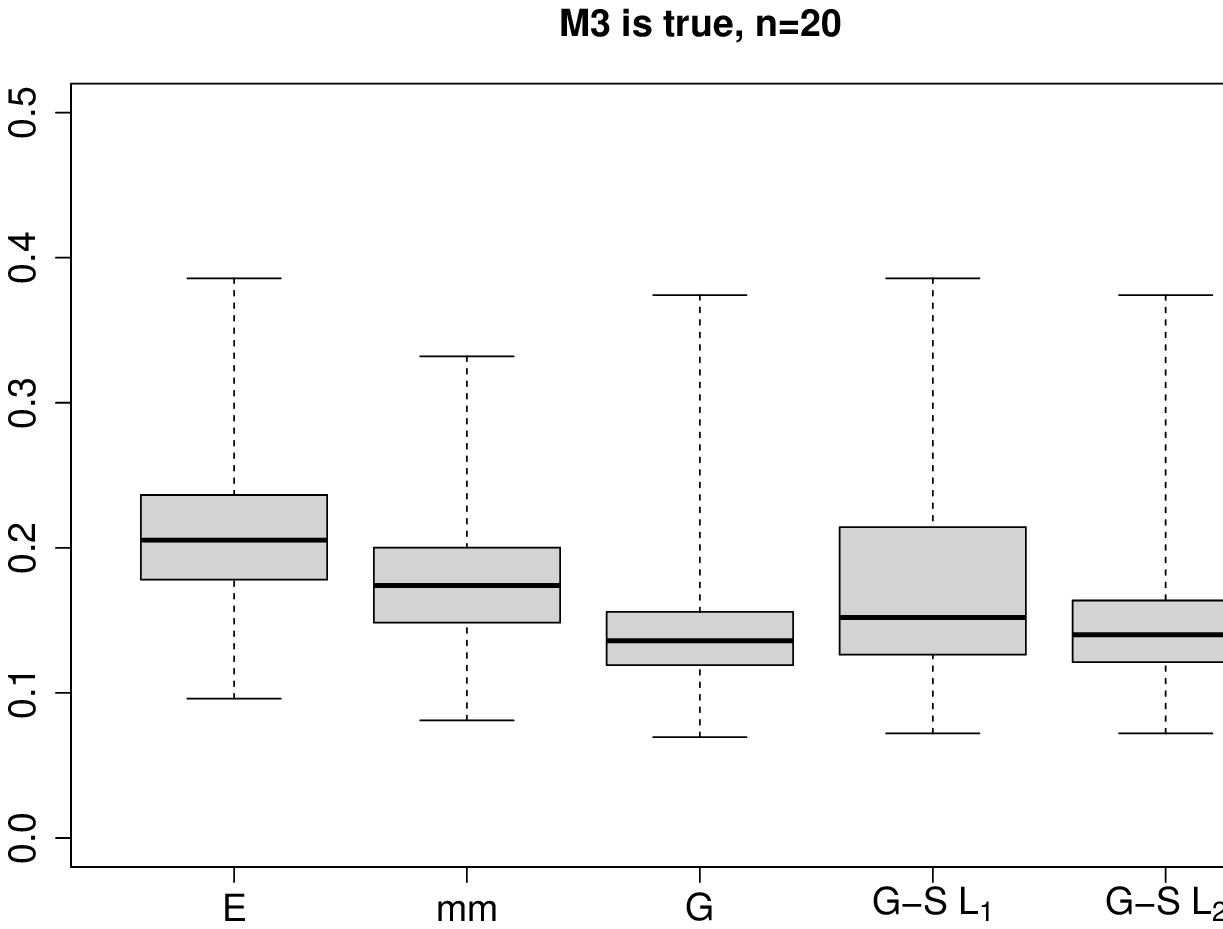}
  \end{subfigure}
  
    \begin{subfigure}{5.2cm}
    \centering\includegraphics[scale=0.17]{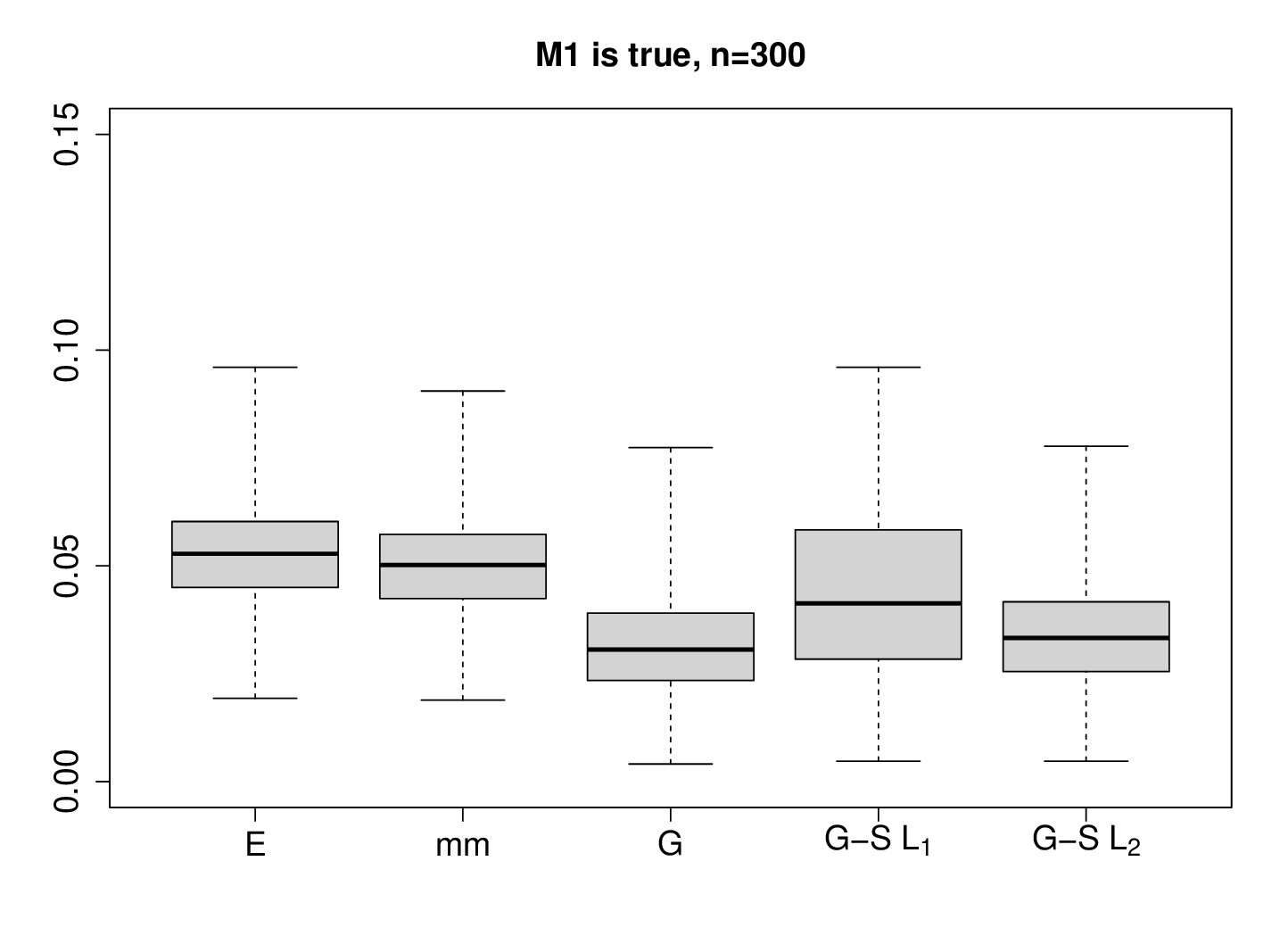}
  \end{subfigure}
  \begin{subfigure}{5.2cm}
    \centering\includegraphics[scale=0.17]{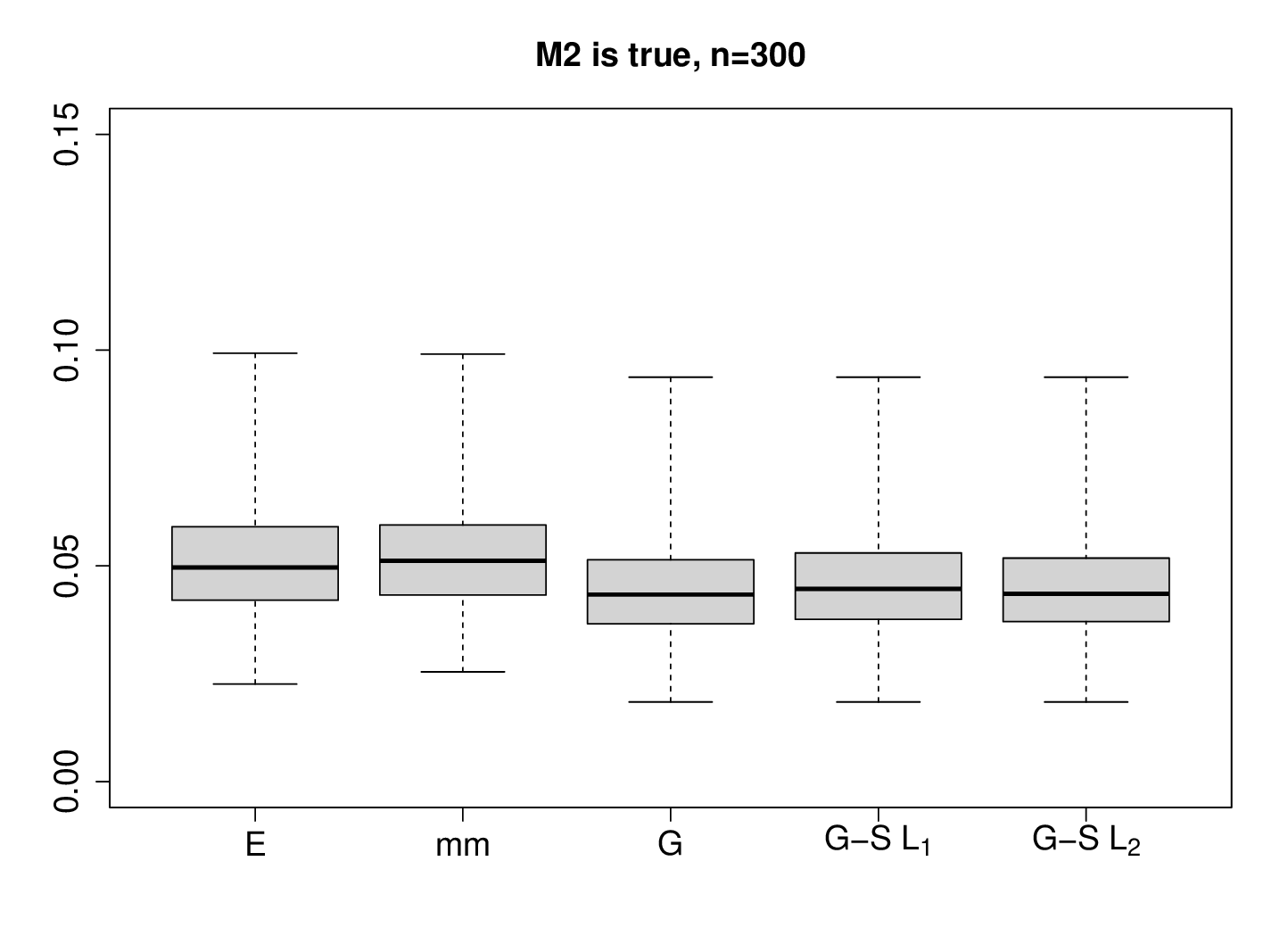}
  \end{subfigure}
  \begin{subfigure}{5.2cm}
    \centering\includegraphics[scale=0.17]{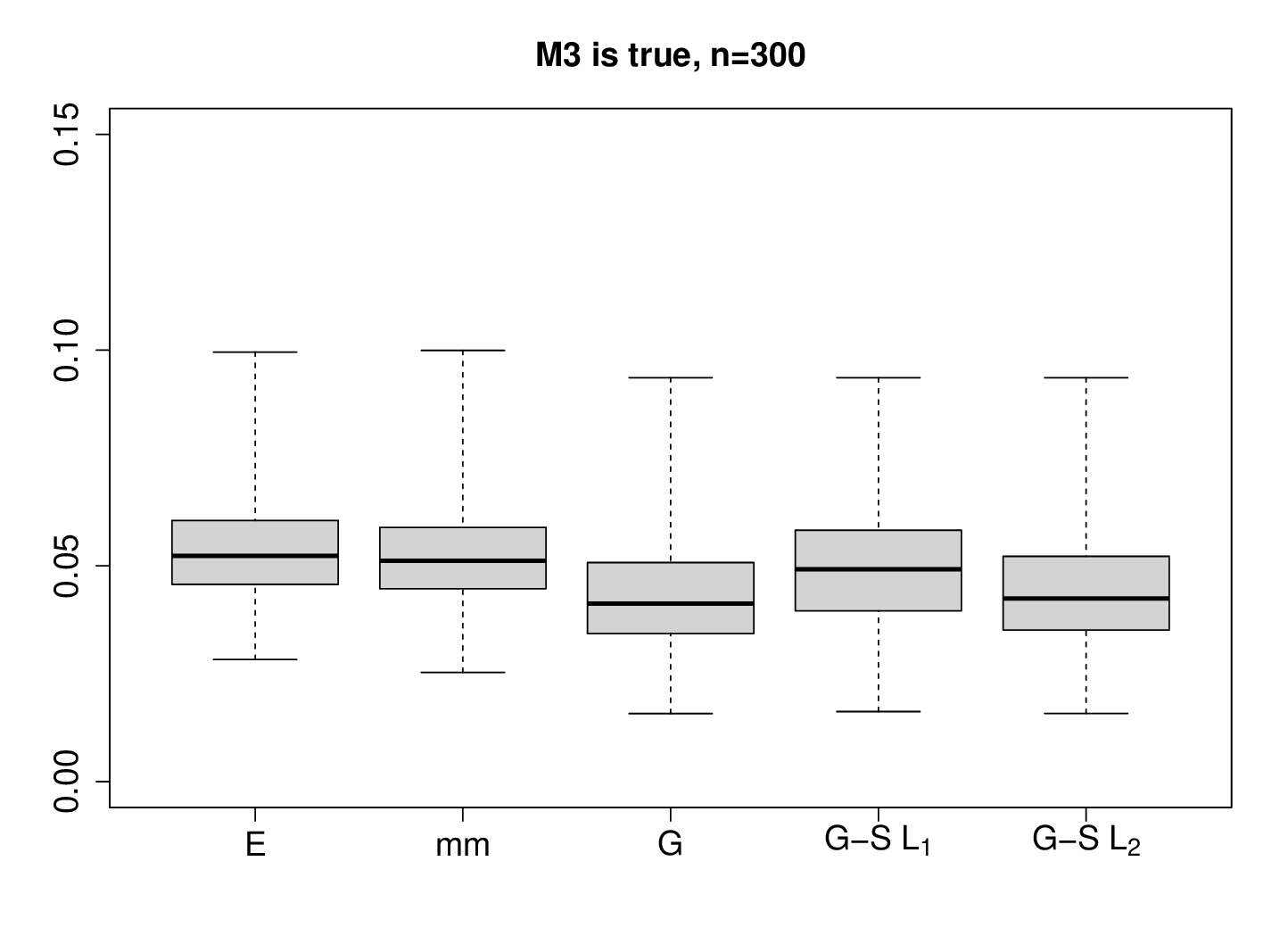}
  \end{subfigure}  
   \caption{The boxplots for $\ell_1$-distances (first and second rows from the top) and $\ell_2$-distances (third and fourth rows)  of the estimators: the empirical estimator $(E)$, minimax esimator $(mm)$, Grenander estimator $(G)$, Grenander-Stone estimator with $L_{1}$ penalization (G--S $L_{1}$) and Grenander-Stone estimator with $L_{2}$ penalization (G--S  $L_{2})$ for the models \textbf{M1}, \textbf{M2} and \textbf{M3}.}\label{decr_l1}
\end{figure} 

  

\begin{figure}[!htbp] 
  \begin{subfigure}{5.2cm}
    \centering\includegraphics[scale=0.27]{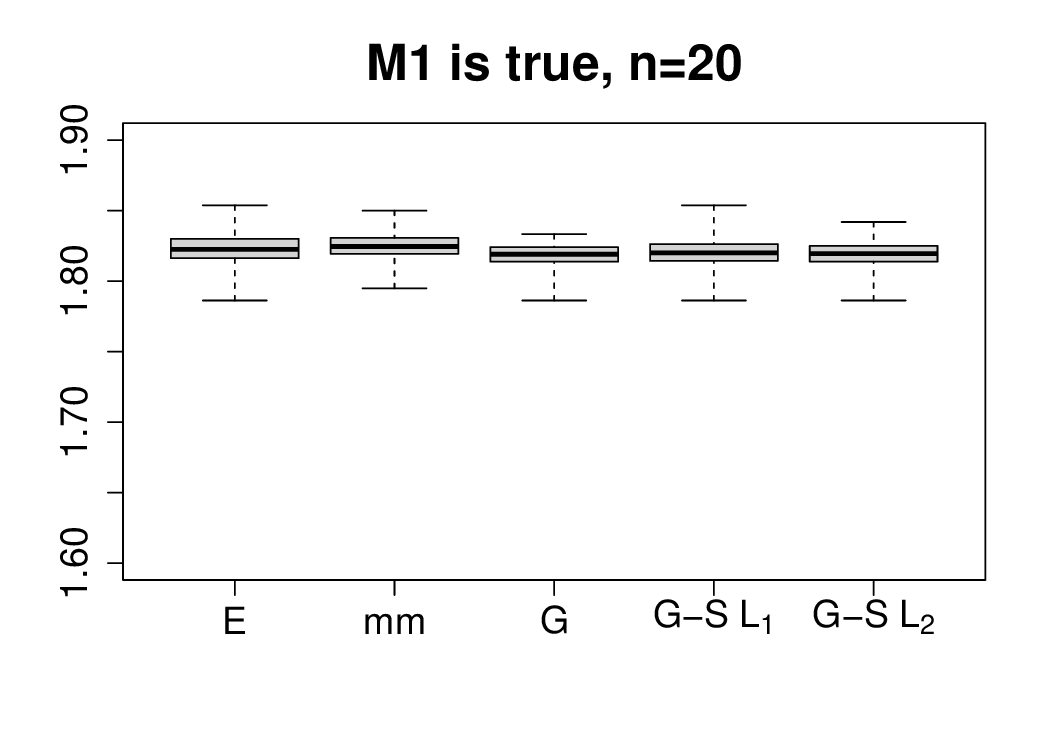}
  \end{subfigure}
  \begin{subfigure}{5.2cm}
    \centering\includegraphics[scale=0.27]{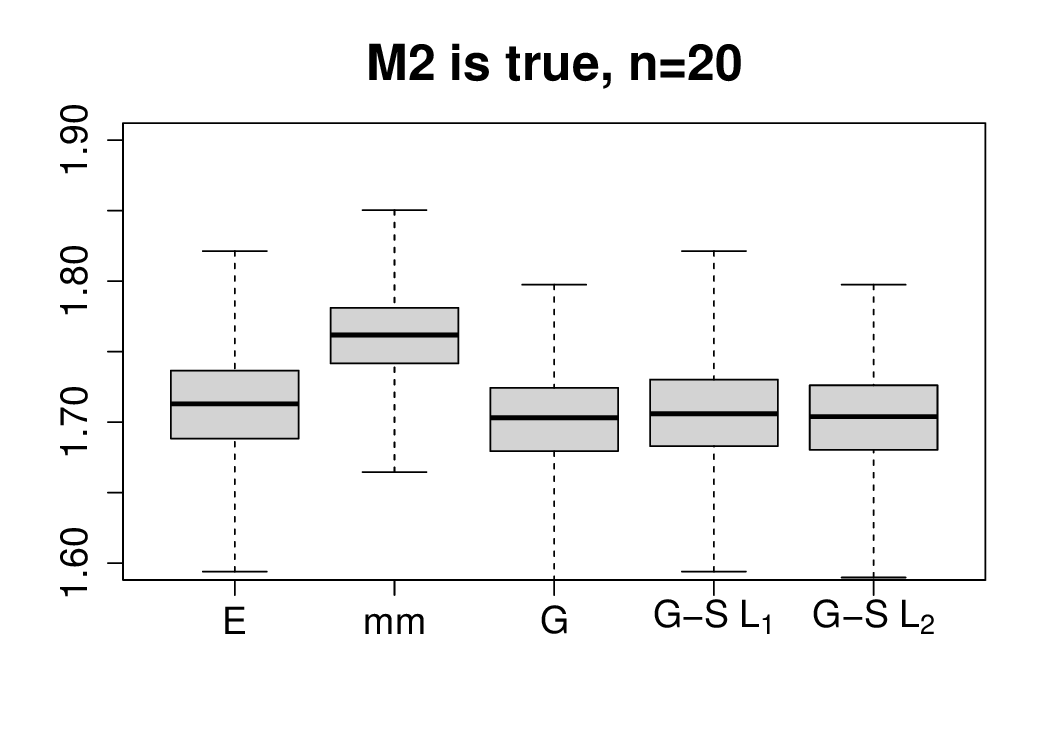}
  \end{subfigure}
  \begin{subfigure}{5.2cm}
    \centering\includegraphics[scale=0.27]{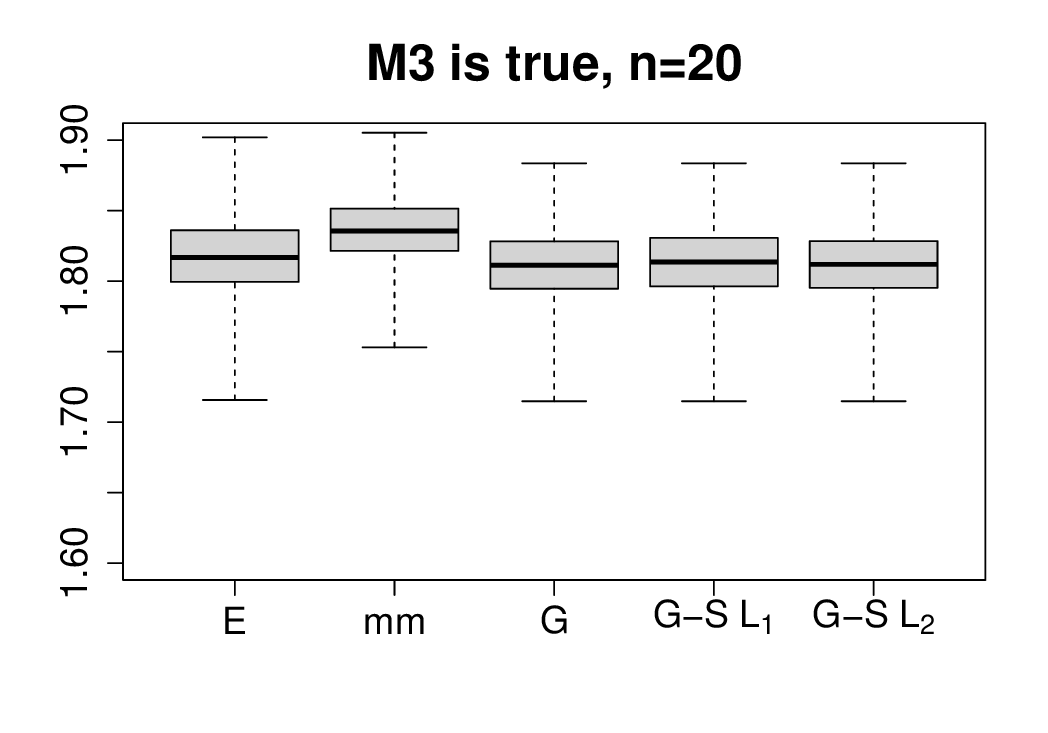}
  \end{subfigure}
  
    \begin{subfigure}{5.2cm}
    \centering\includegraphics[scale=0.27]{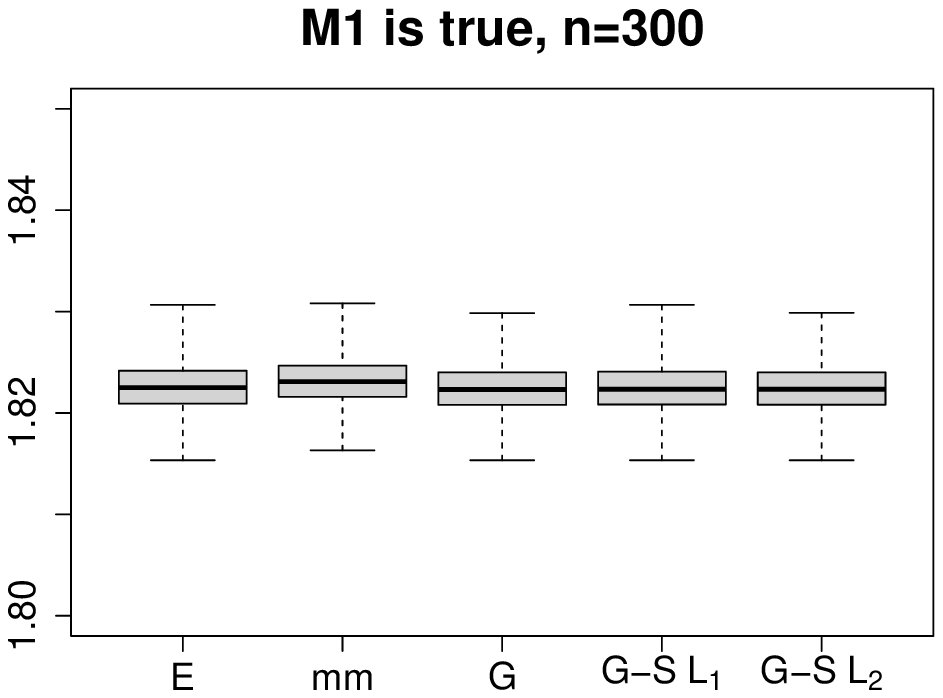}
  \end{subfigure}
  \begin{subfigure}{5.2cm}
    \centering\includegraphics[scale=0.27]{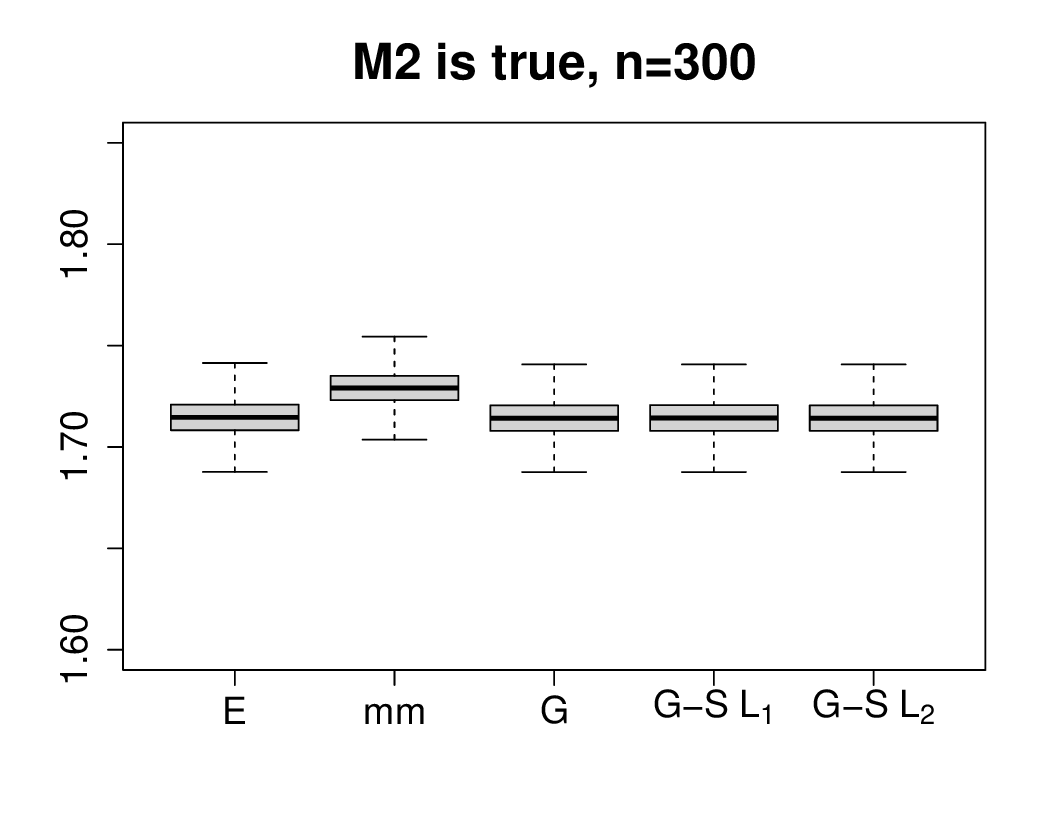}
  \end{subfigure}
  \begin{subfigure}{5.2cm}
    \centering\includegraphics[scale=0.27]{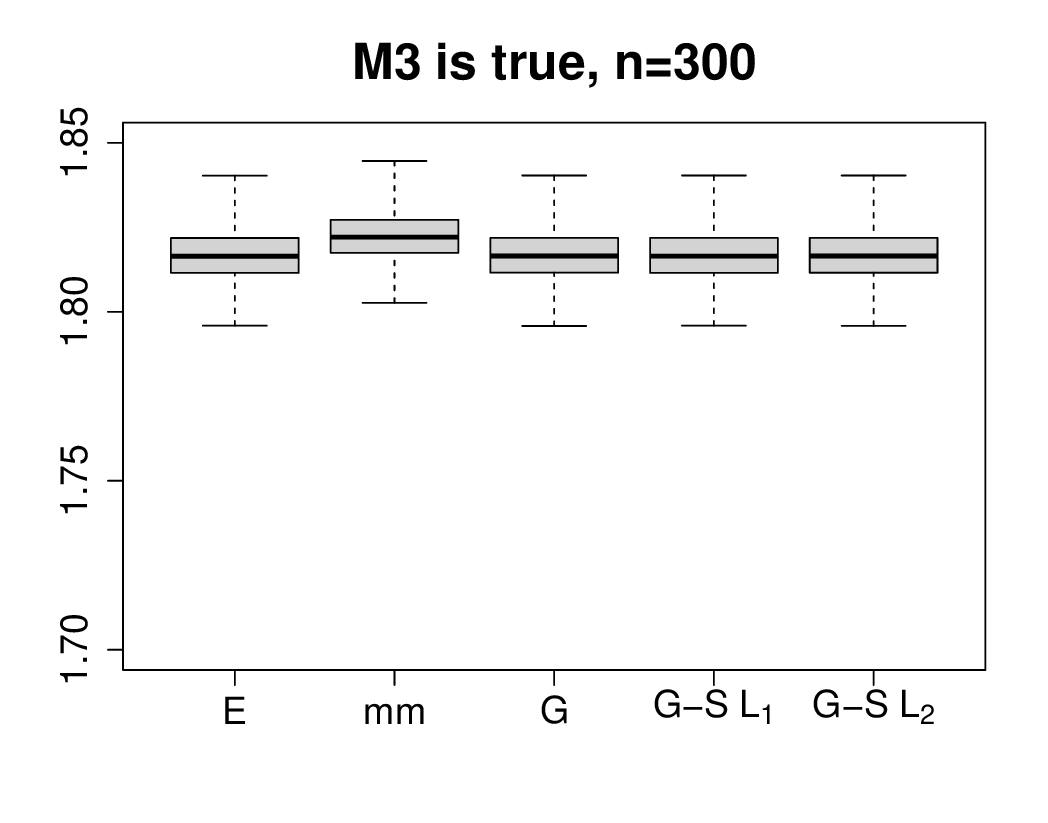}
  \end{subfigure}
  
    \begin{subfigure}{5.2cm}
    \centering\includegraphics[scale=0.27]{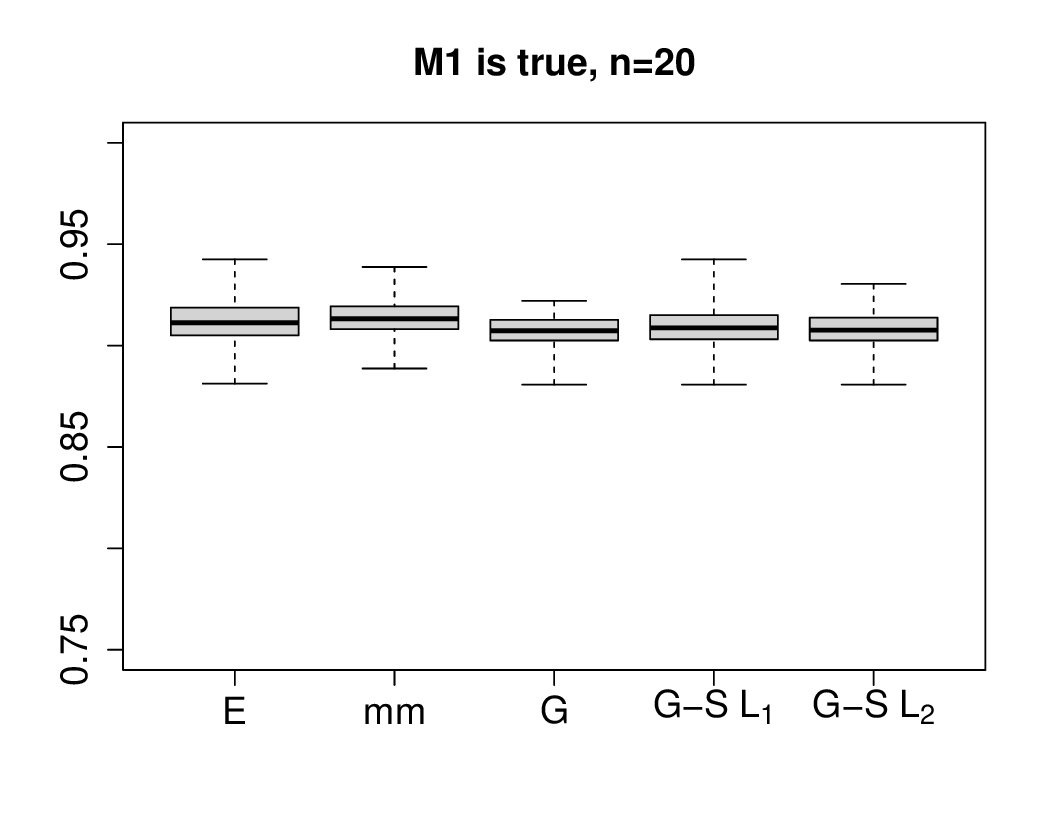}
  \end{subfigure}
  \begin{subfigure}{5.2cm}
    \centering\includegraphics[scale=0.27]{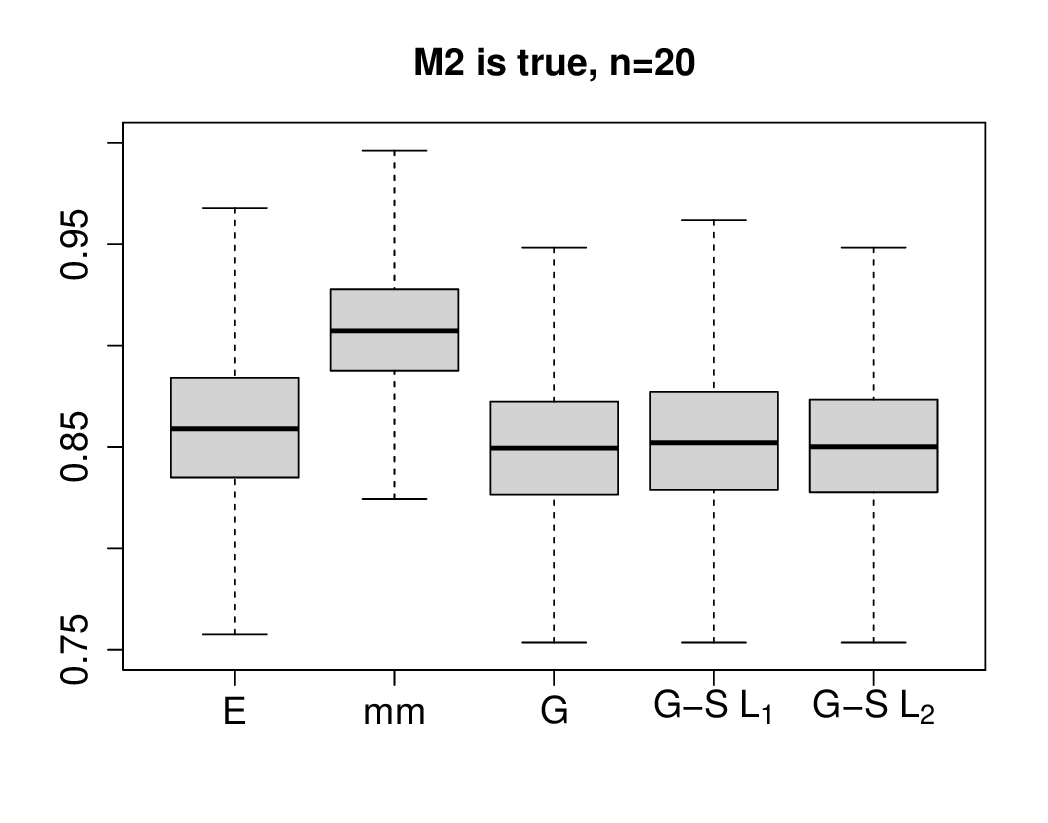}
  \end{subfigure}
  \begin{subfigure}{5.2cm}
    \centering\includegraphics[scale=0.27]{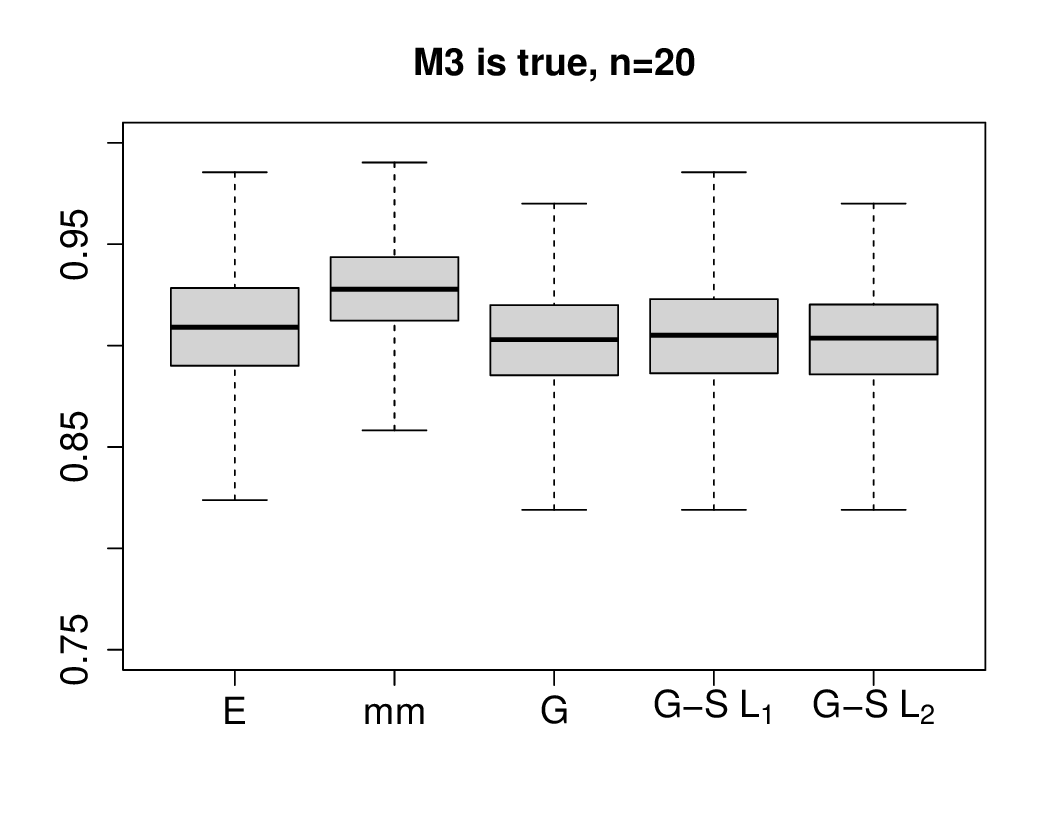}
  \end{subfigure}
  
    \begin{subfigure}{5.2cm}
    \centering\includegraphics[scale=0.27]{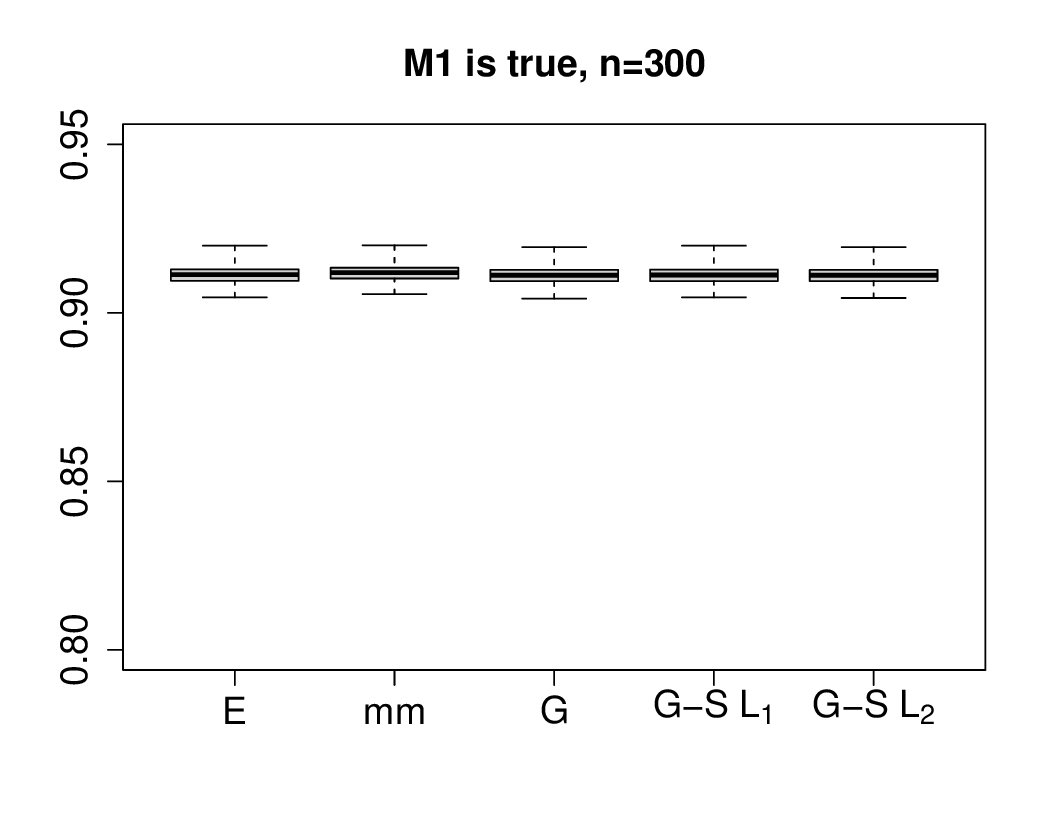}
  \end{subfigure}
  \begin{subfigure}{5.2cm}
    \centering\includegraphics[scale=0.27]{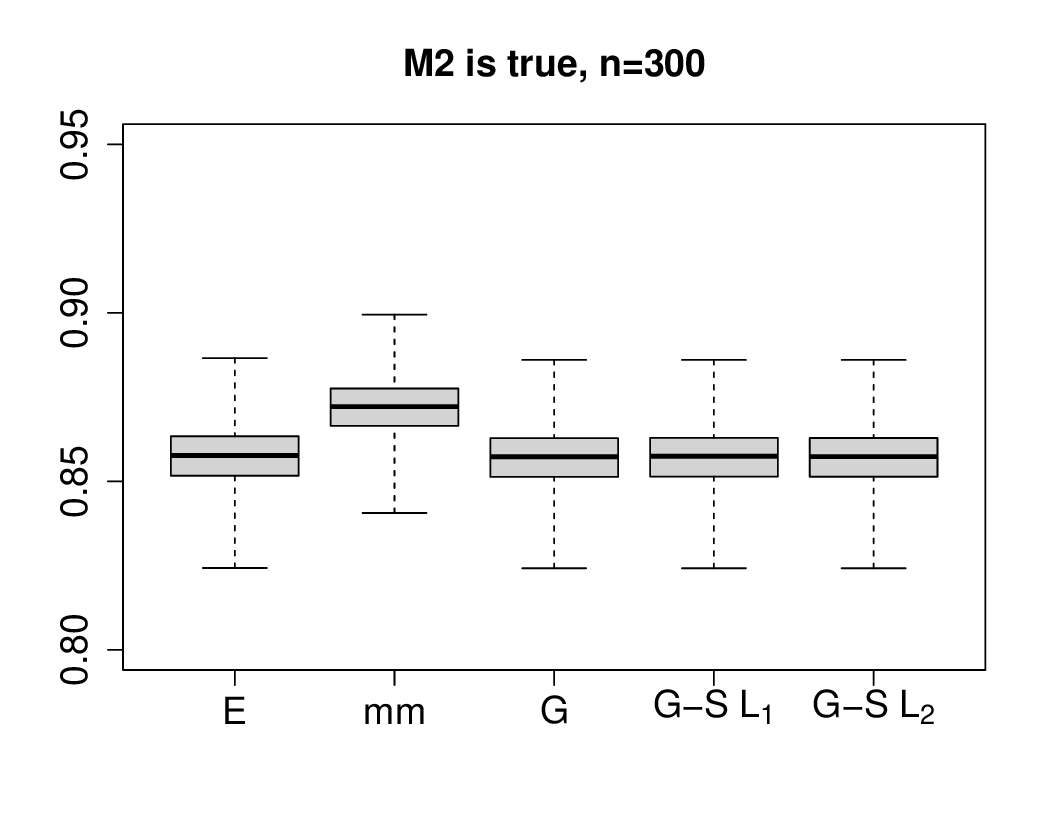}
  \end{subfigure}
  \begin{subfigure}{5.2cm}
    \centering\includegraphics[scale=0.27]{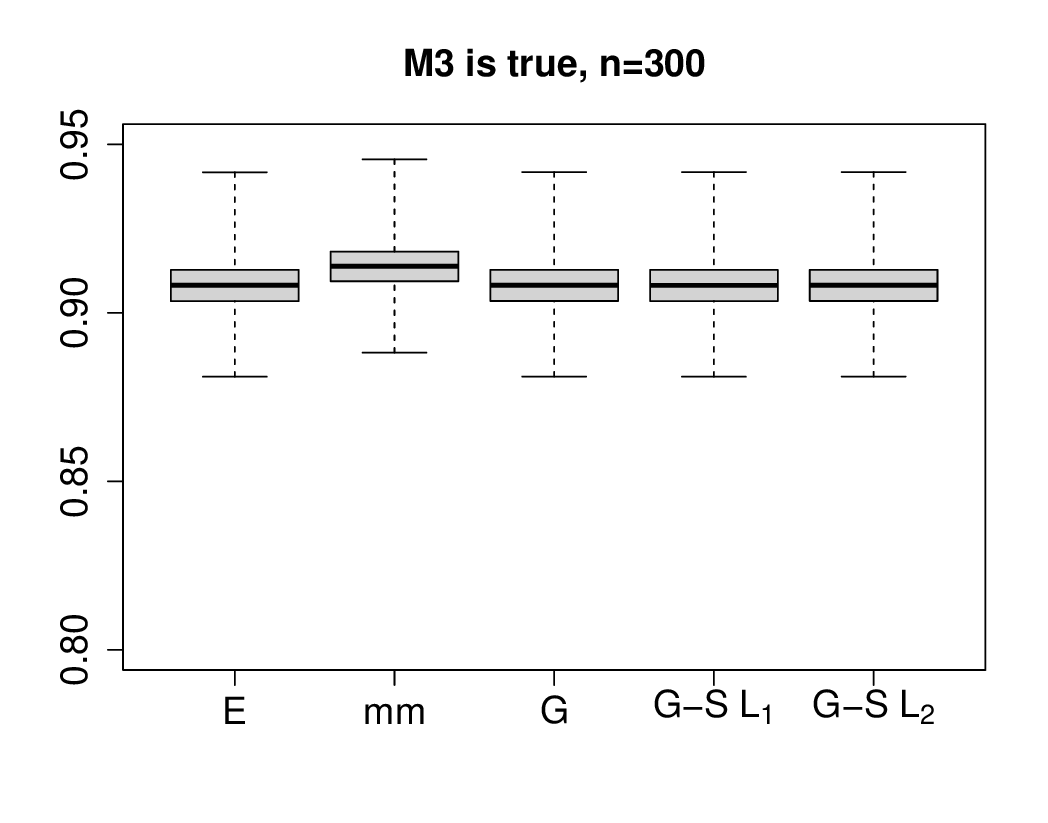}
  \end{subfigure}
   \caption{The boxplots of $1000$ expected $S_1$-scores  (first and second rows from the top) and $S_2$-scores (third and fourth rows) of the estimators: the empirical estimator $(E)$, minimax esimator $(mm)$, Grenander estimator $(G)$, Grenander-Stone estimator with $L_{1}$ penalization (G--S $L_{1}$) and Grenander-Stone estimator with $L_{2}$ penalization (G--S  $L_{2})$ for the models \textbf{M1}, \textbf{M2} and \textbf{M3}.}\label{s1ndecr_l1}
\end{figure} 

  

The performance of Stone-Grenander estimator for isotonic underlying distributions is illustrated at Figures \ref{decr_l1} and \ref{s1ndecr_l1}. From Figure \ref{decr_l1} one can see that Grenander-Stone estimator performs better with $L_{2}$ penalization than with $L_{1}$. Also, it performs better than the empirical and minimax estimators and almost as good as Grenander estimator. 
Further, in a sense of $S_{1}$ and $S_{2}$ scores, the performance of stacked estimators is visually the same for $L_{1}$ and $L_{2}$, it is better compare to minimax estimator and slightly better that the performance of the empirical estimator.
 
\begin{figure}[!htbp] 
  \begin{subfigure}{5.2cm}
    \centering\includegraphics[scale=0.27]{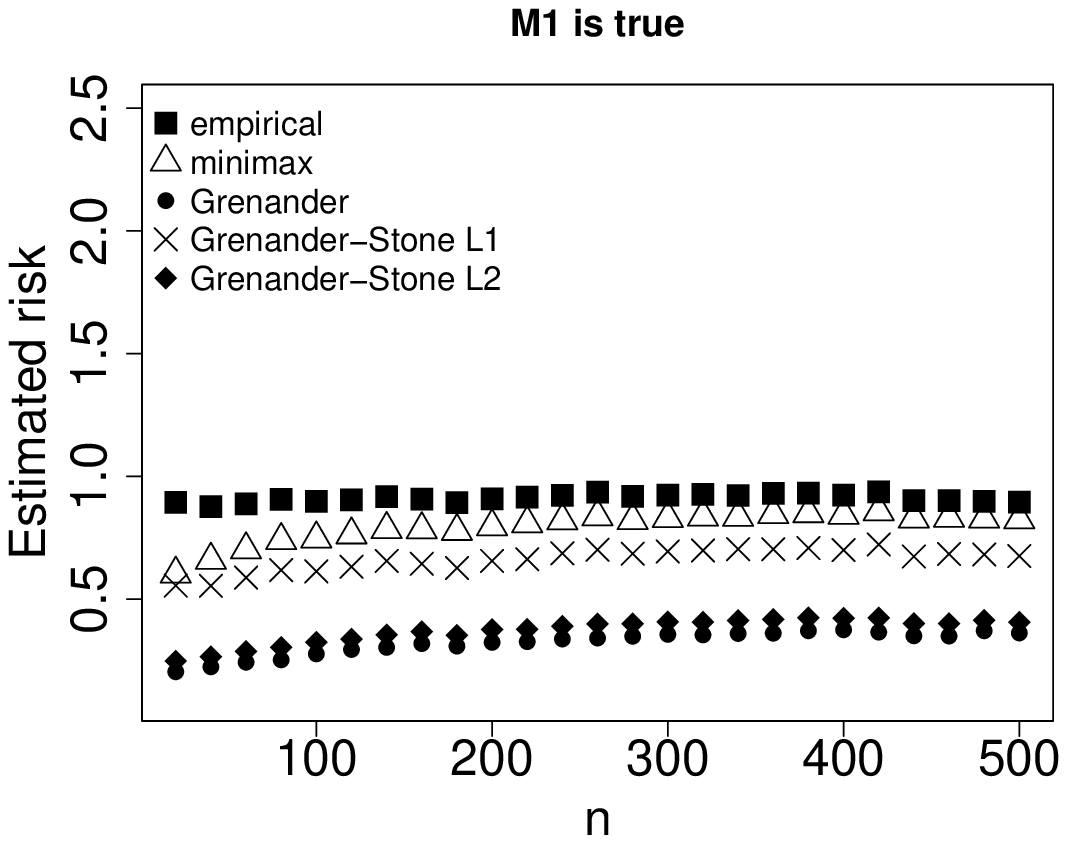}
  \end{subfigure}
  \begin{subfigure}{5.2cm}
    \centering\includegraphics[scale=0.27]{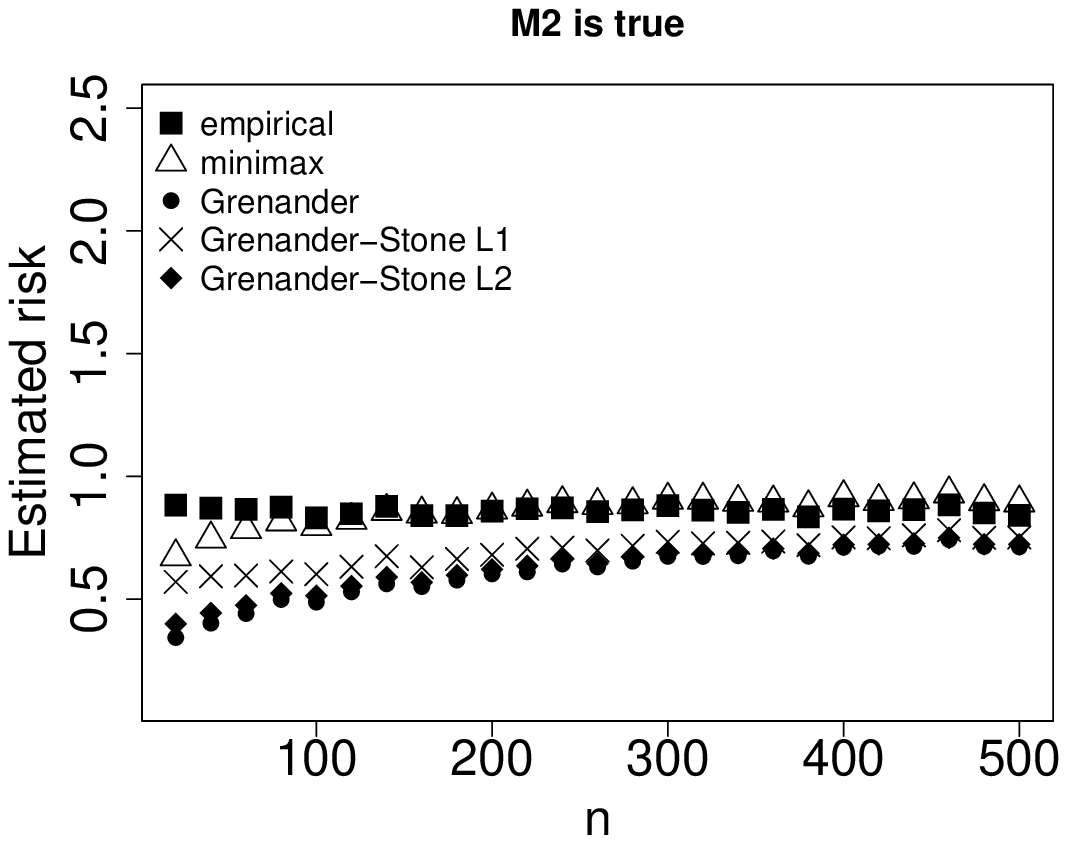}
  \end{subfigure}
  \begin{subfigure}{5.2cm}
    \centering\includegraphics[scale=0.27]{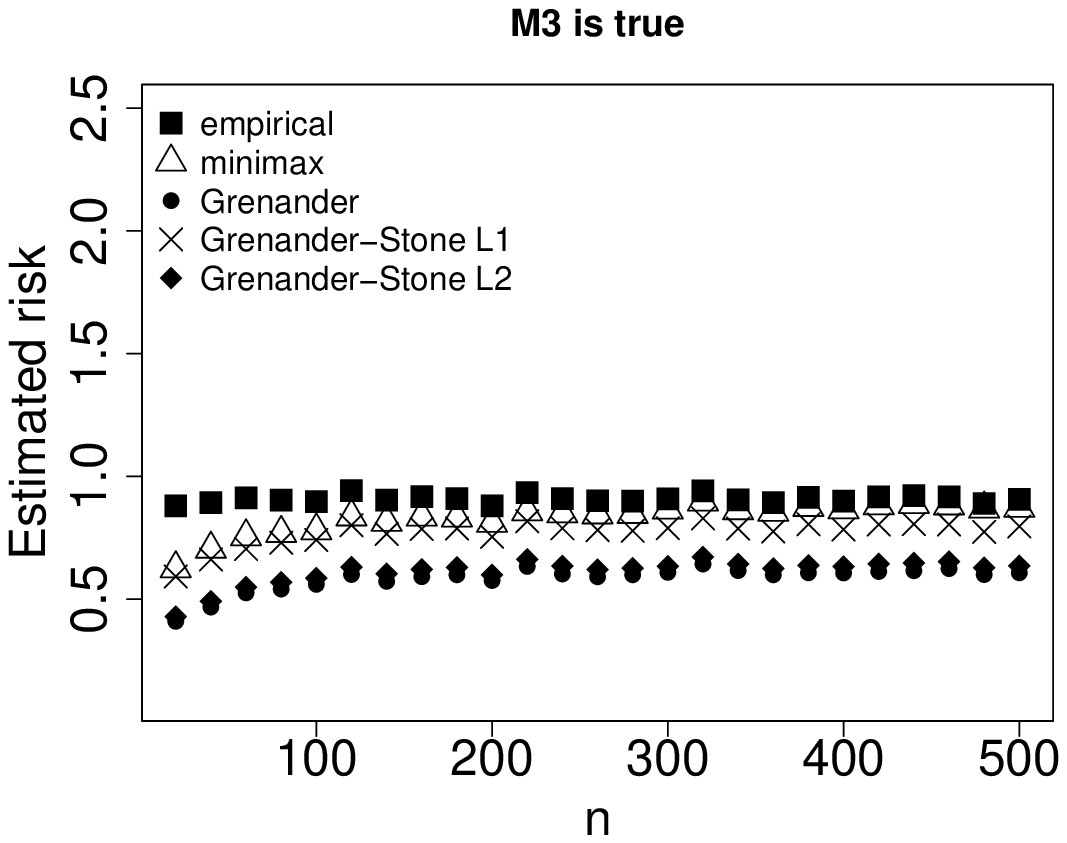}
  \end{subfigure}
   \caption{The estimates of the scaled risk for the models \textbf{M1}, \textbf{M2} and \textbf{M3}.}\label{decrRisk}
\end{figure}

Next, the superiority of Grenander-Stone estimator can be shown by ploting the estimates of scaled risk $n\mathbb{E}[||\hat{\bm{\xi}}_{n} - \bm{p}||_{2}^{2}]$ (with $\hat{\bm{\xi}}_{n}$ one of the following estimators under consideration) versus the sample size $n$, based on $1000$ Monte Carlo simulations. From Figure \ref{decrRisk} we can conclude that in the case of isotonic underlying distribution Grenander-Stone estimator with $L_{2}$ cross-validation performs almost as good as Grenander estimator and it performs significantly better than the empirical and the minimax estimators.

\subsubsection{True distribution is not isotonic}

Now let us consider the case when the underlying distribution is not isotonic:
\begin{eqnarray*}
\textbf{M4}: \bm{p} &=& NBin(7, 0.4), \\
\textbf{M5}: \bm{p} &=& \frac{3}{8}Pois(2) + \frac{5}{8}Pois(15),\\
\textbf{M6}: \bm{p} &=& q_{1}U_{2d}(1) + \dots + q_{5}U_{2d}(5),
\end{eqnarray*}
where  $\bm{q} = (-0.1, 0.4, 0.3, 0.2, 0.2)$ and $NBin$ is negative binomial distribution and $Pois$ is Poisson distribution.

\begin{figure}[!htbp] 
  \begin{subfigure}{5.2cm}
    \centering\includegraphics[scale=0.19]{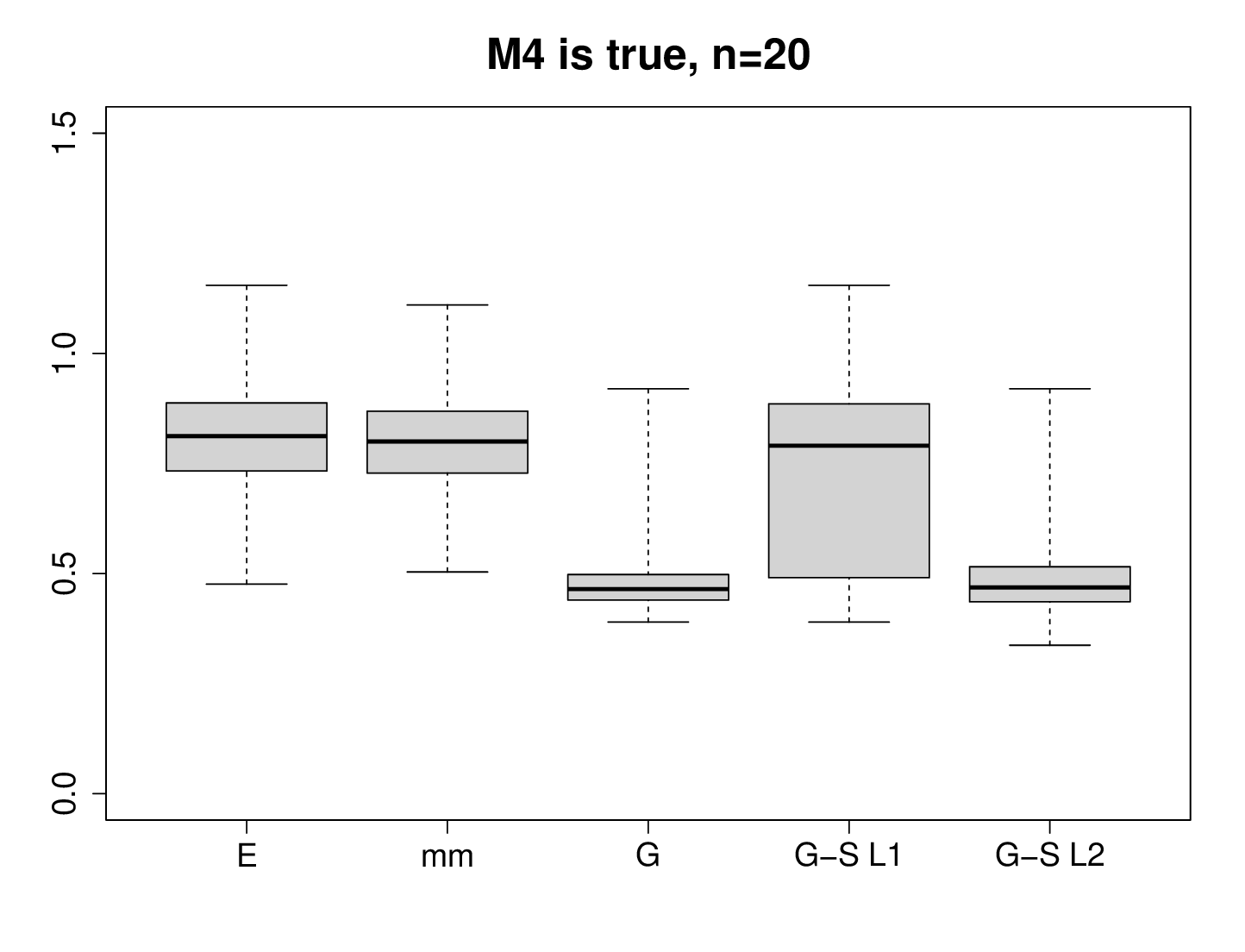}
  \end{subfigure}
  \begin{subfigure}{5.2cm}
    \centering\includegraphics[scale=0.19]{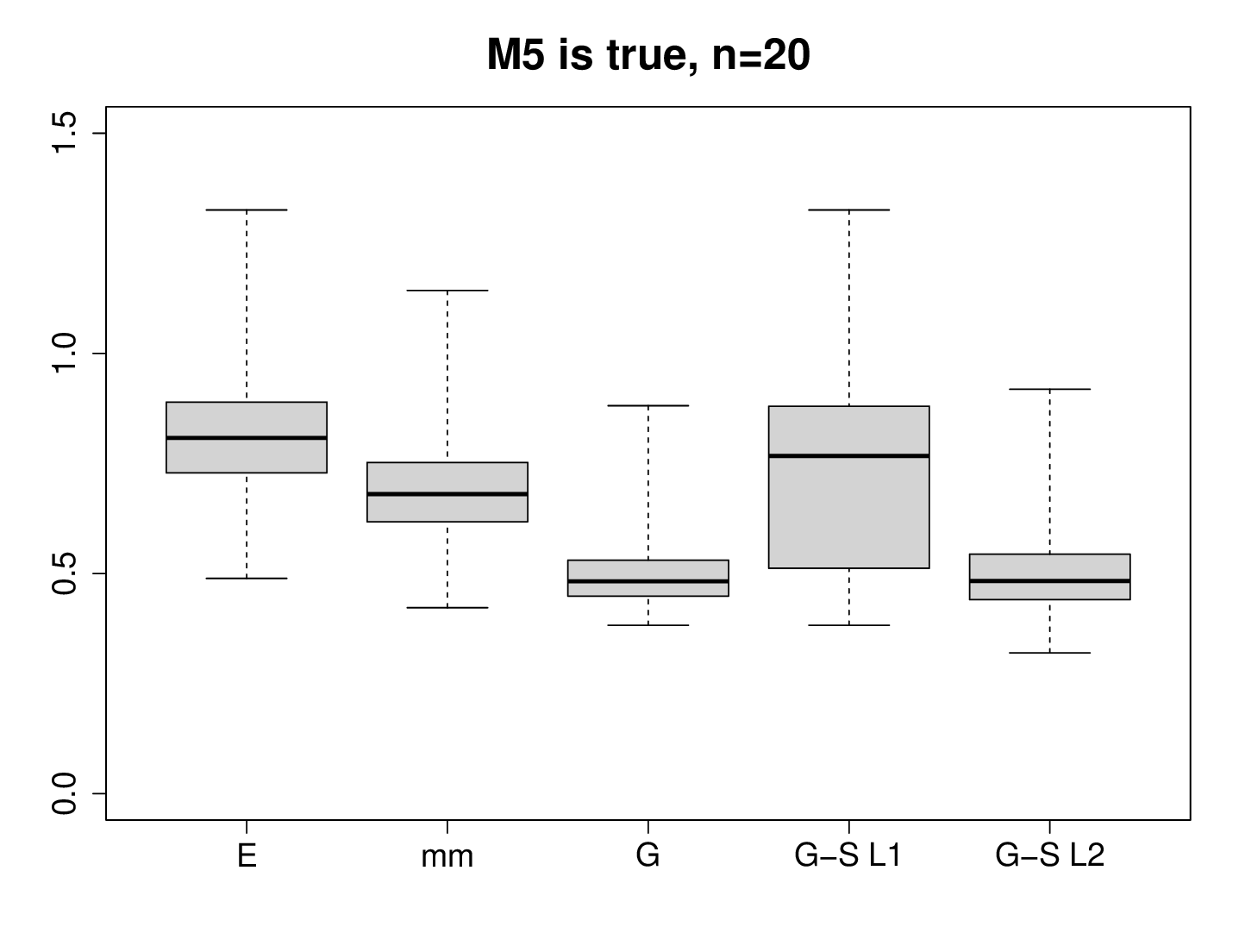}
  \end{subfigure}
  \begin{subfigure}{5.2cm}
    \centering\includegraphics[scale=0.19]{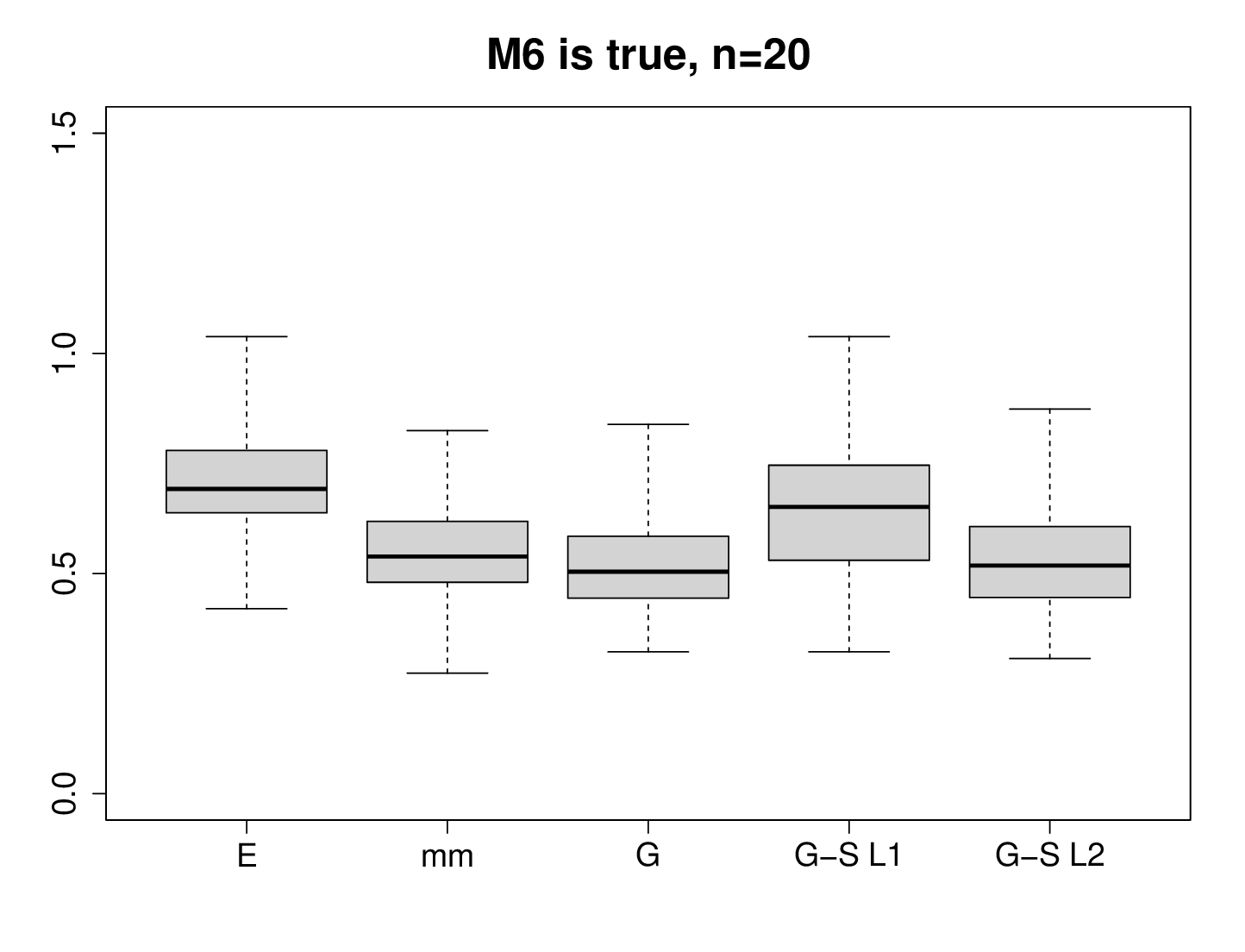}
  \end{subfigure}
  
    \begin{subfigure}{5.2cm}
    \centering\includegraphics[scale=0.19]{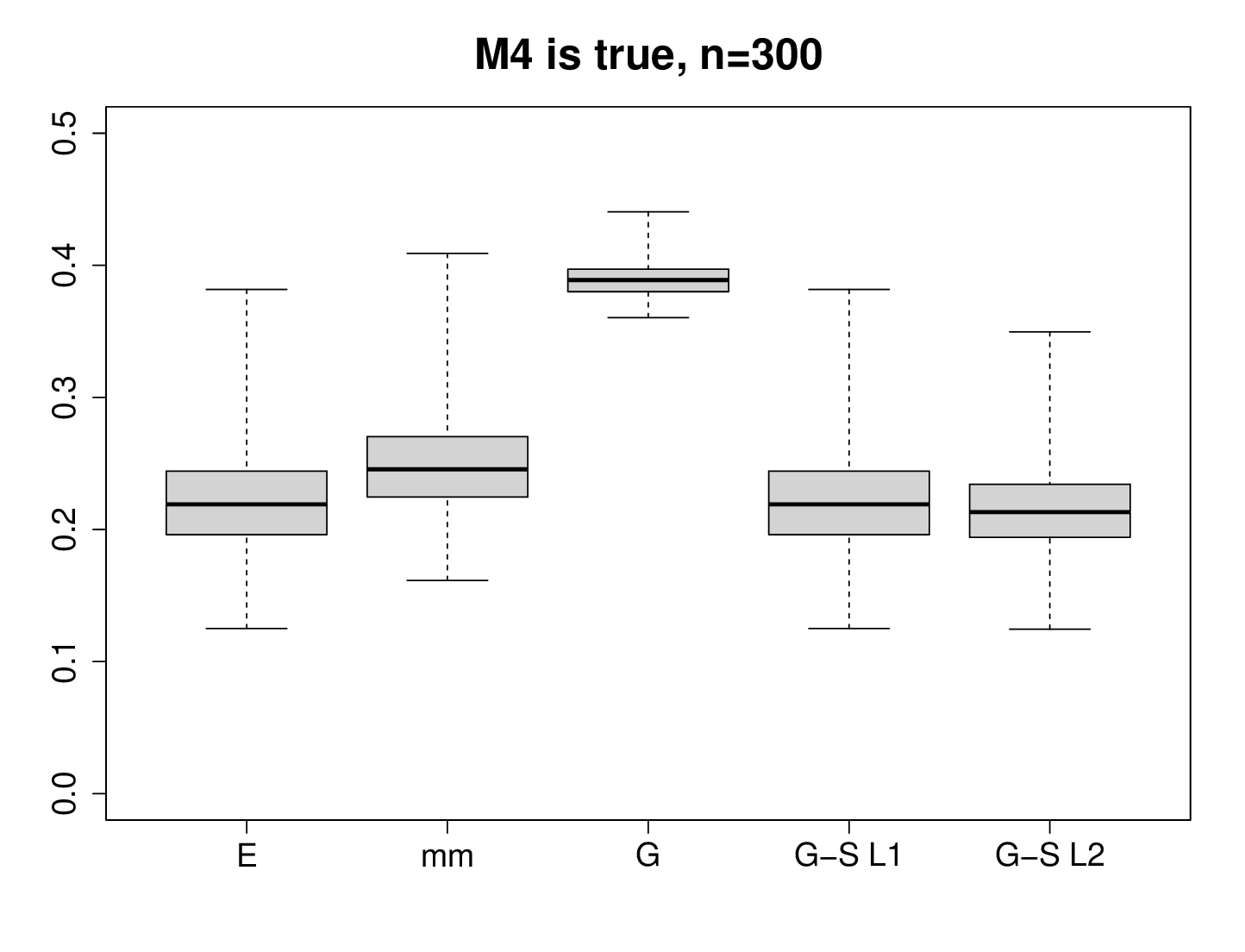}
  \end{subfigure}
  \begin{subfigure}{5.2cm}
    \centering\includegraphics[scale=0.19]{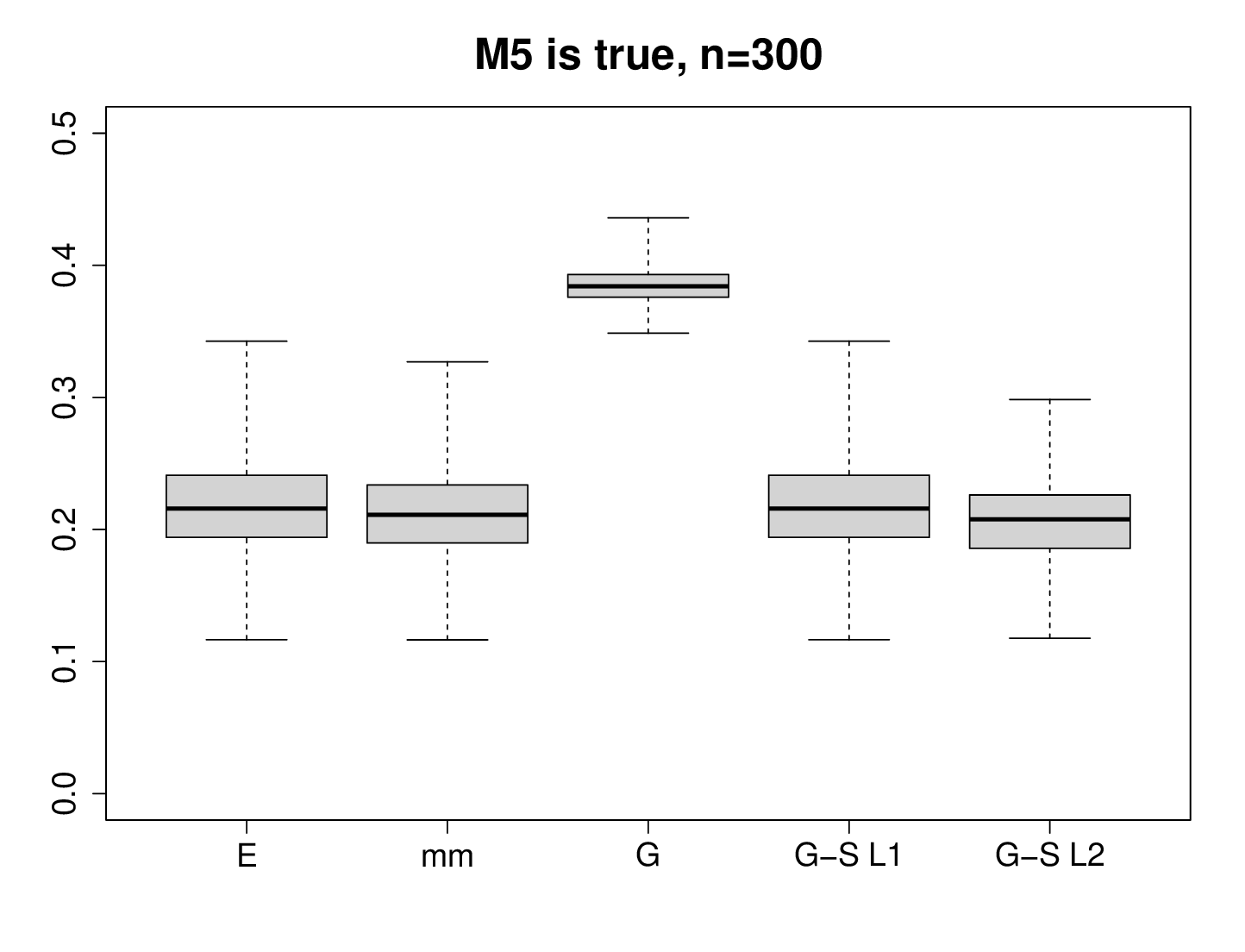}
  \end{subfigure}
  \begin{subfigure}{5.2cm}
    \centering\includegraphics[scale=0.19]{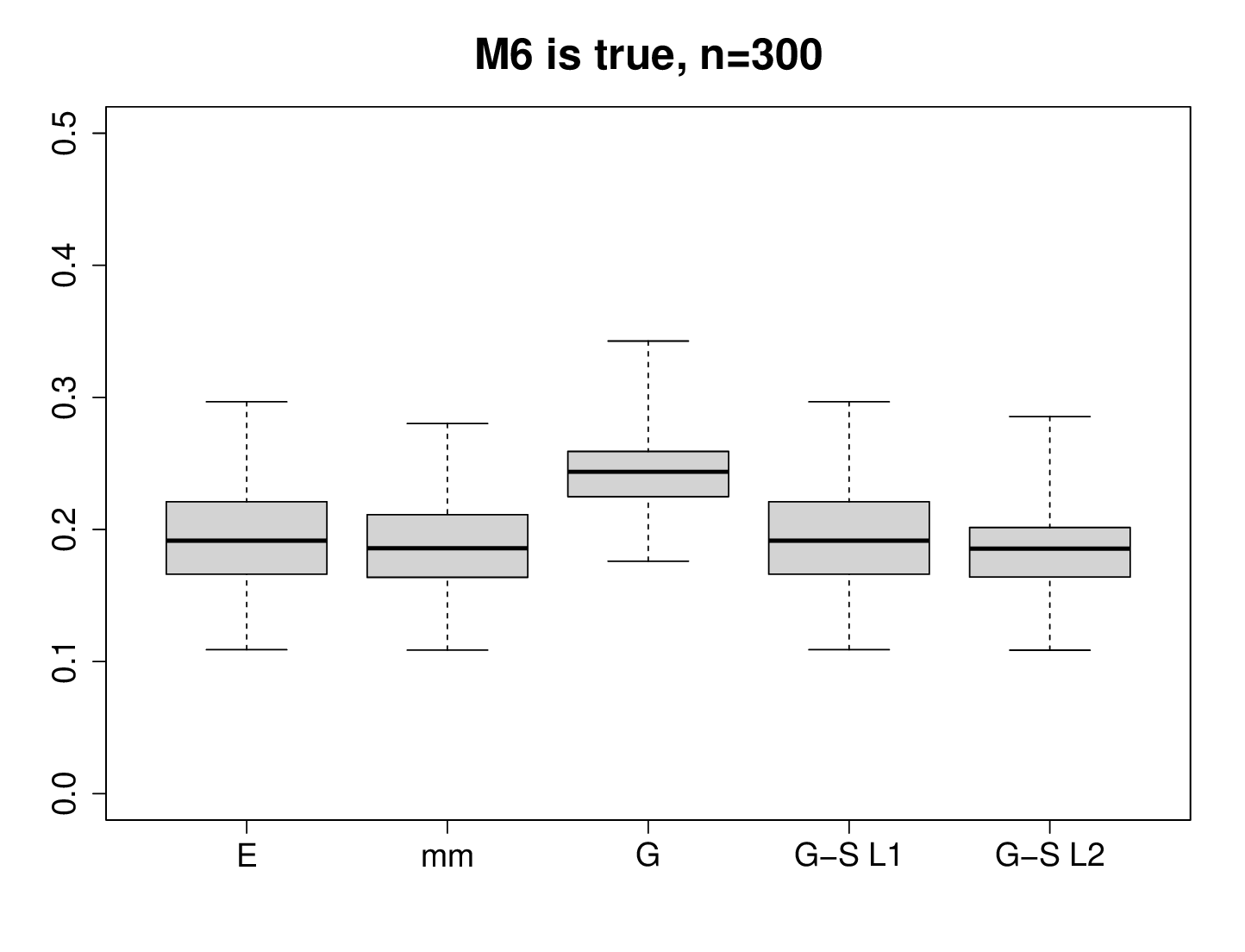}
  \end{subfigure}
  
    \begin{subfigure}{5.2cm}
    \centering\includegraphics[scale=0.19]{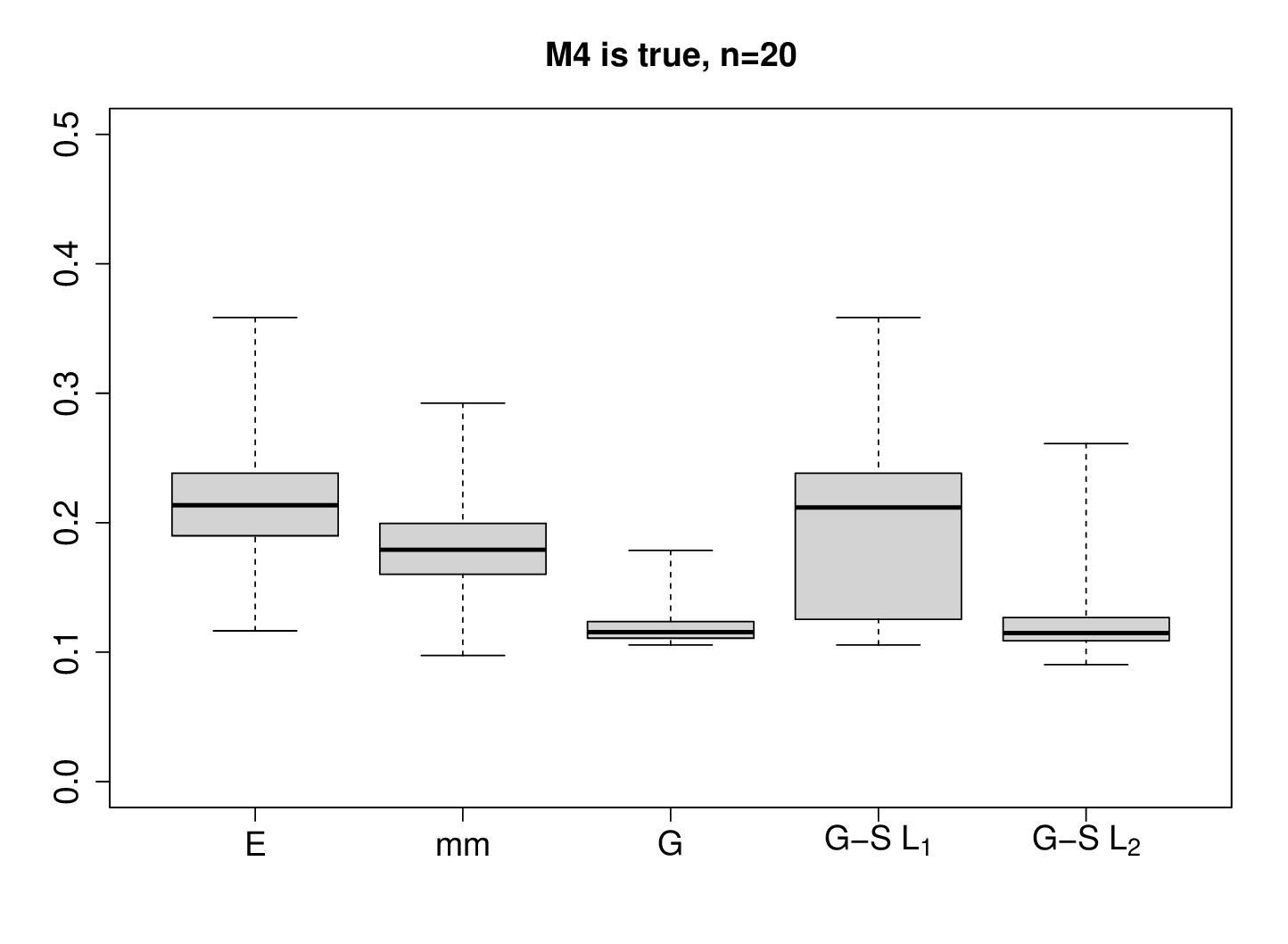}
  \end{subfigure}
  \begin{subfigure}{5.2cm}
    \centering\includegraphics[scale=0.19]{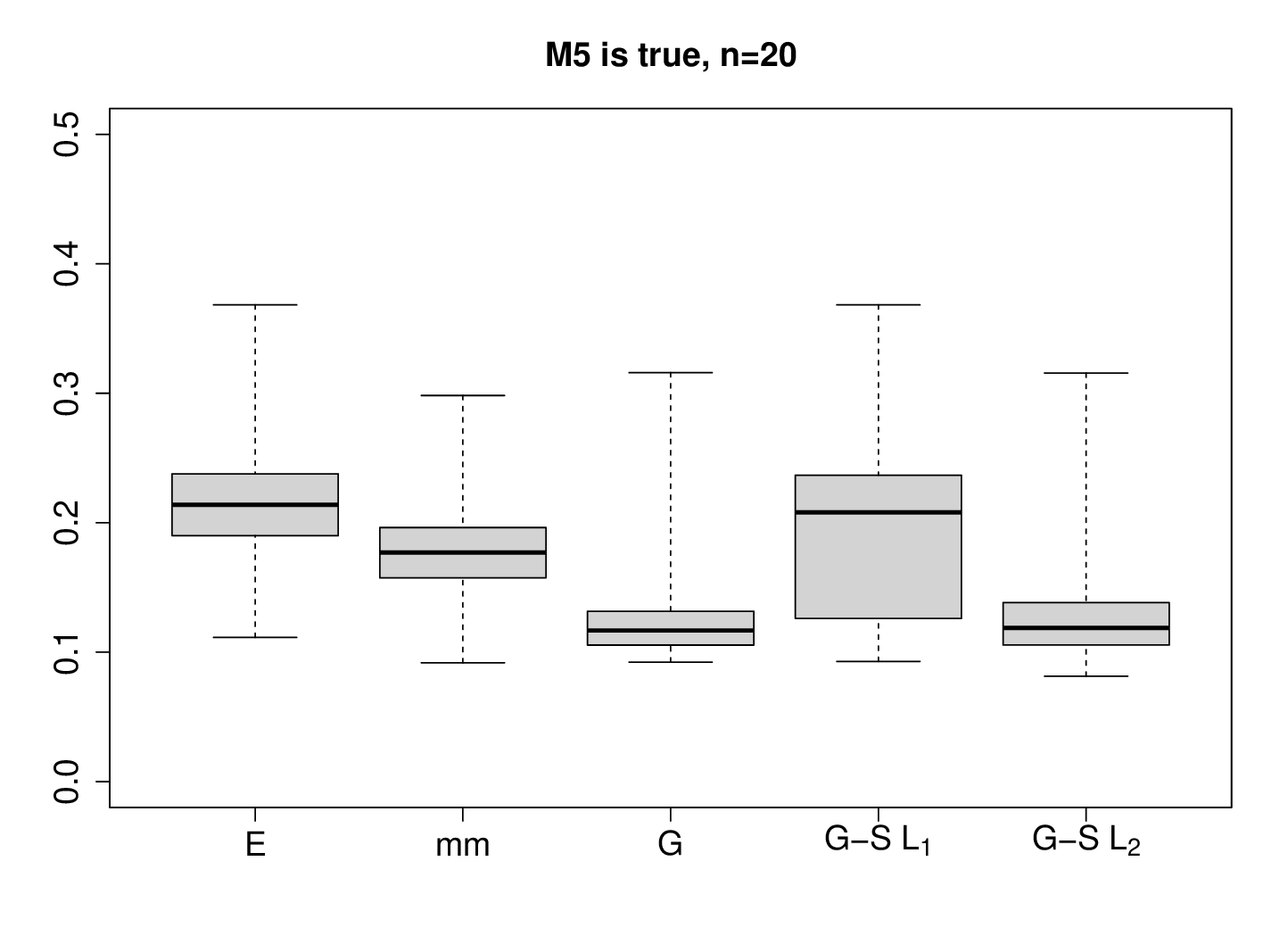}
  \end{subfigure}
  \begin{subfigure}{5.2cm}
    \centering\includegraphics[scale=0.19]{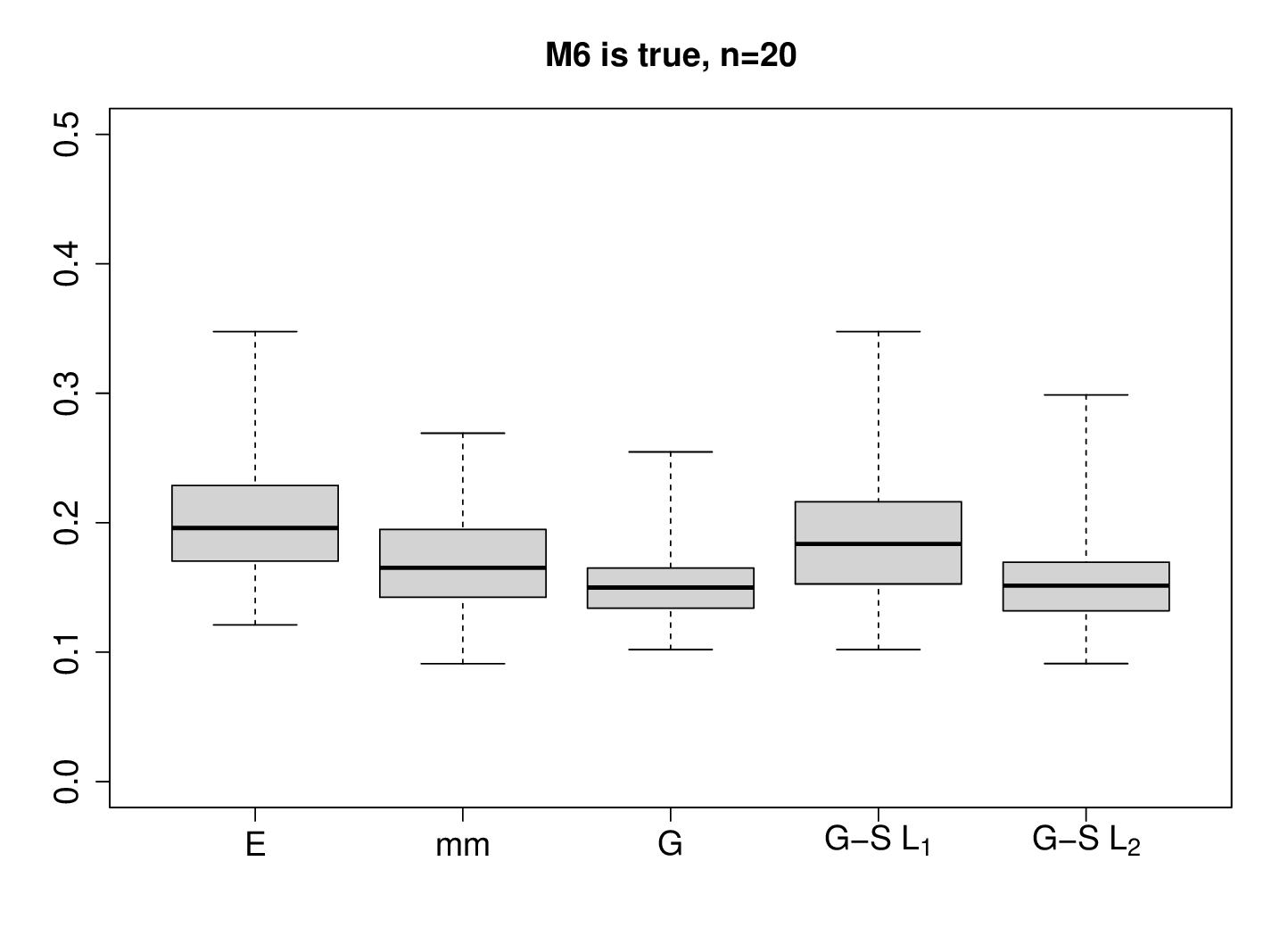}
  \end{subfigure}
  
    \begin{subfigure}{5.2cm}
    \centering\includegraphics[scale=0.19]{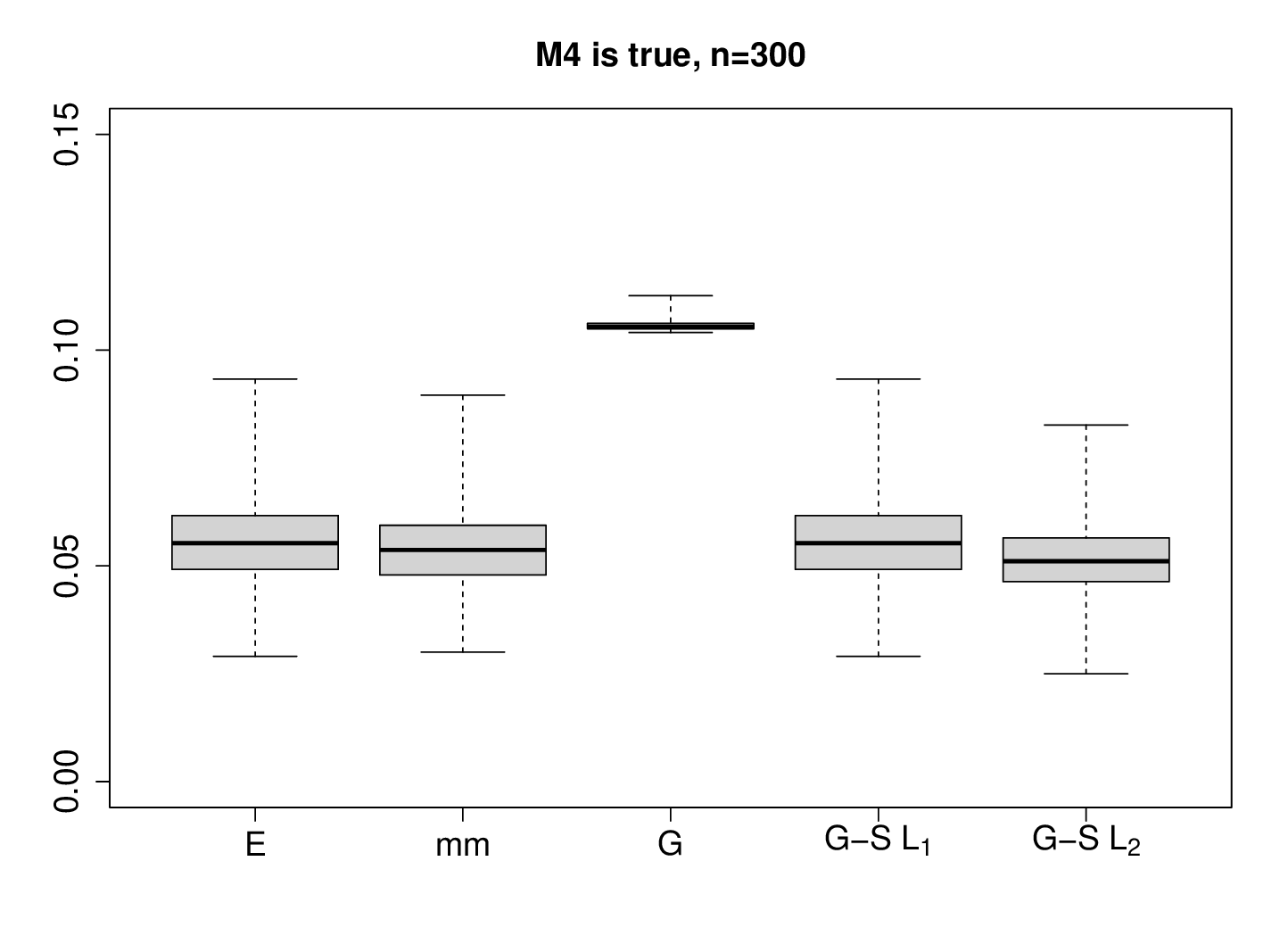}
  \end{subfigure}
  \begin{subfigure}{5.2cm}
    \centering\includegraphics[scale=0.19]{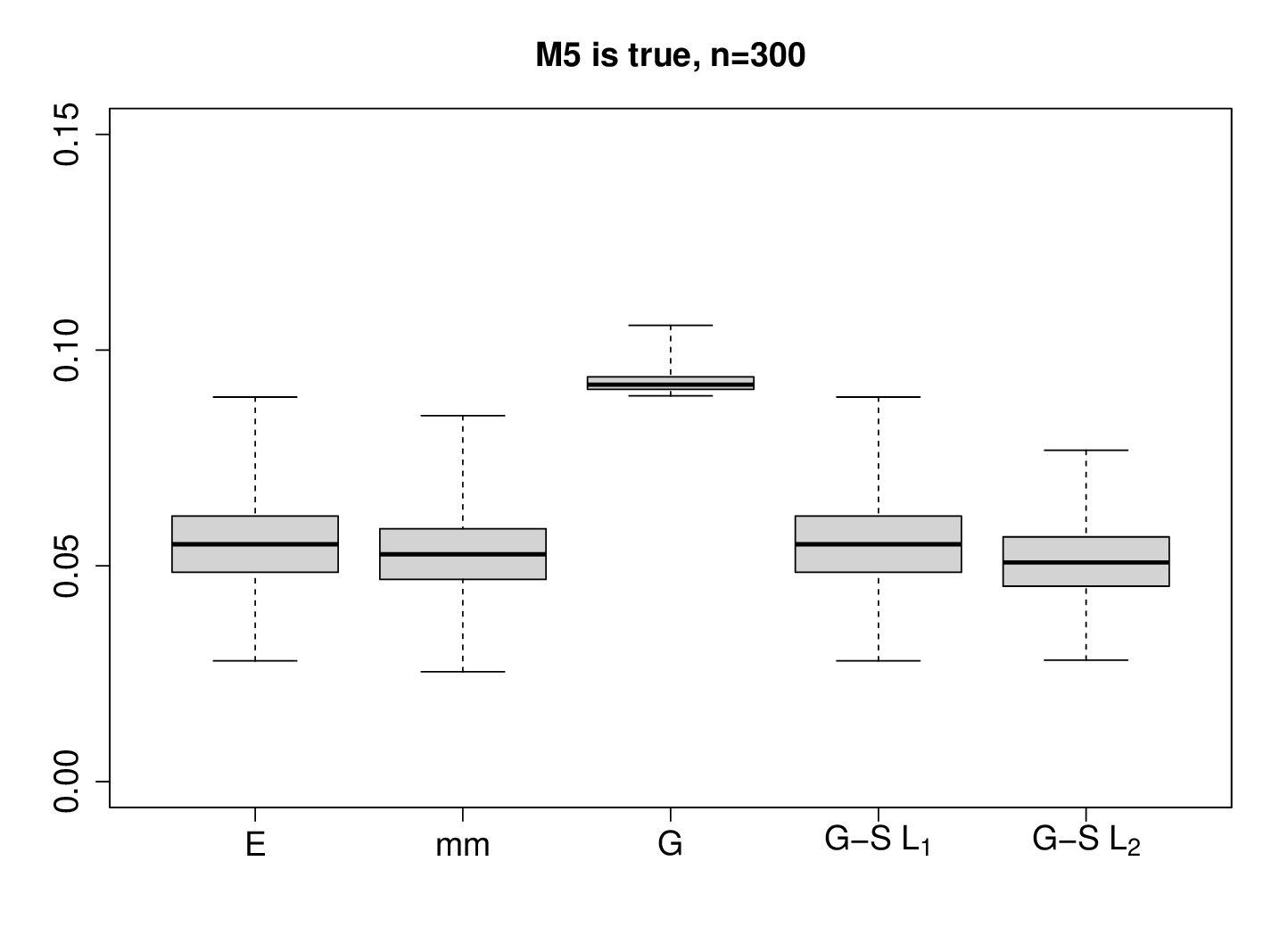}
  \end{subfigure}
  \begin{subfigure}{5.2cm}
    \centering\includegraphics[scale=0.19]{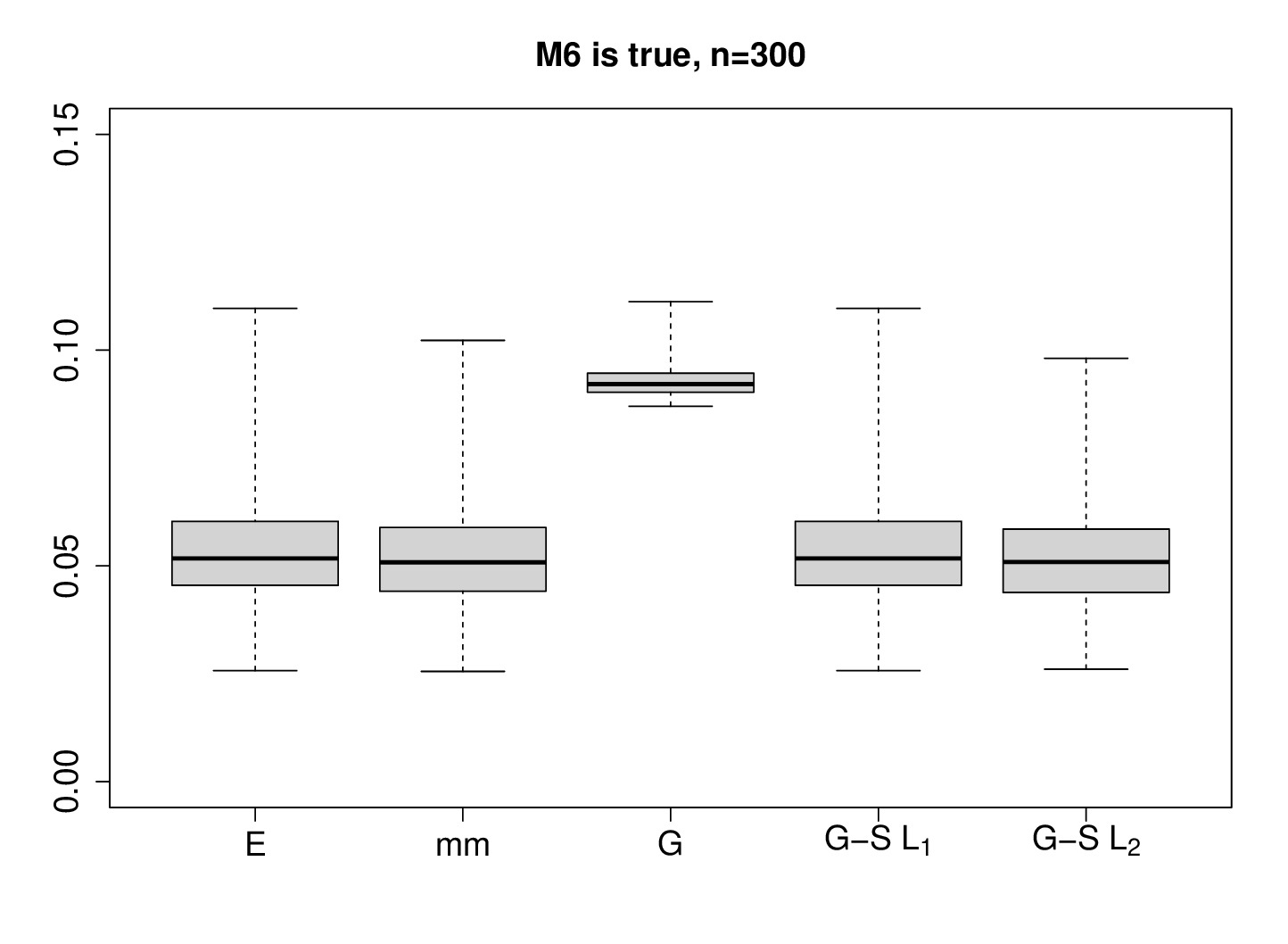}
  \end{subfigure}
   \caption{The boxplots for $\ell_1$-distances (first and second rows from the top) and $\ell_2$-distances (third and fourth rows) of the estimators: the empirical estimator $(e)$, minimax estimator $(mm)$, rearrangement estimator (r), Grenander estimator $(G)$, the stacked rearrangement estimator (sr) and the stacked Grenander estimator $(sG)$ for the models \textbf{M5}, \textbf{M6} and \textbf{M7}.}\label{nondecr_l1}
\end{figure} 

  

\begin{figure}[!htbp] 
  \begin{subfigure}{5.2cm}
    \centering\includegraphics[scale=0.27]{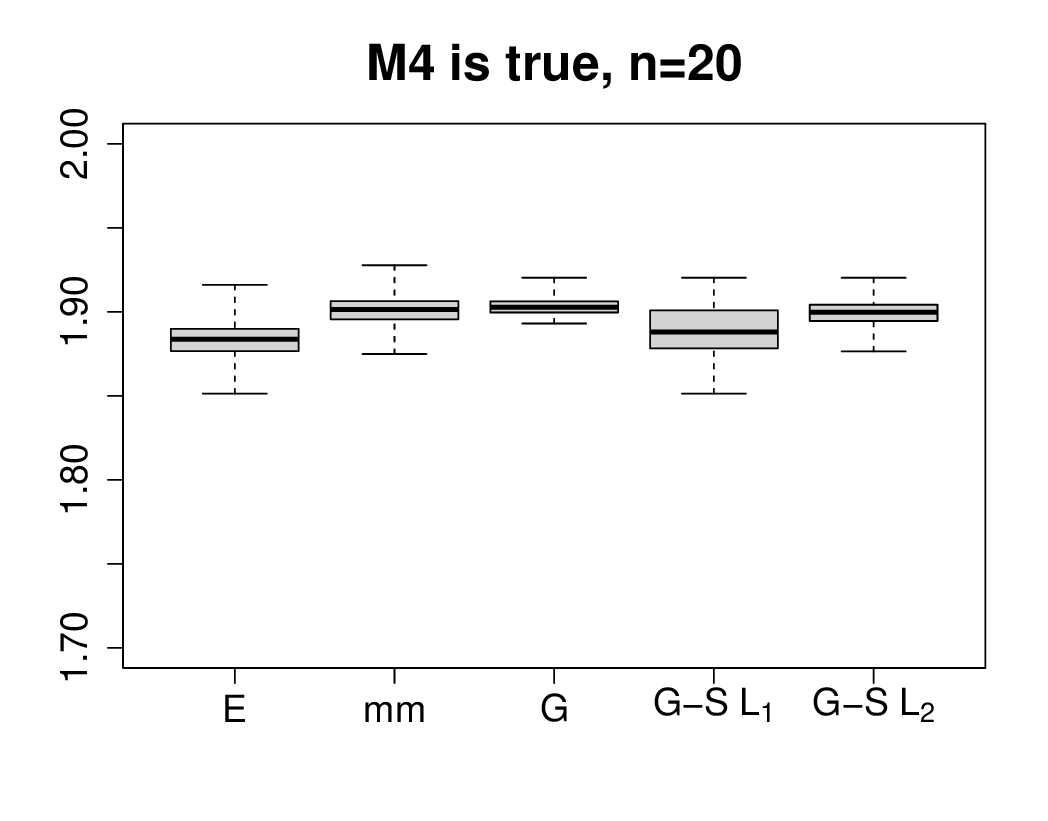}
  \end{subfigure}
  \begin{subfigure}{5.2cm}
    \centering\includegraphics[scale=0.27]{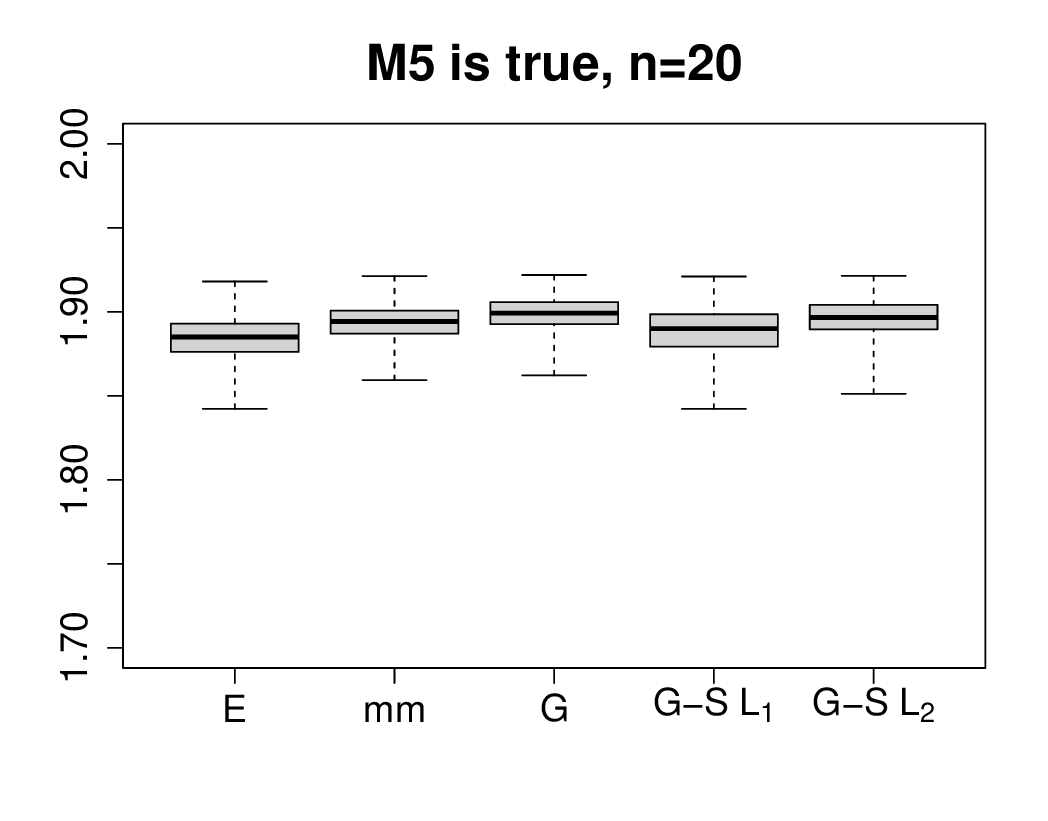}
  \end{subfigure}
  \begin{subfigure}{5.2cm}
    \centering\includegraphics[scale=0.27]{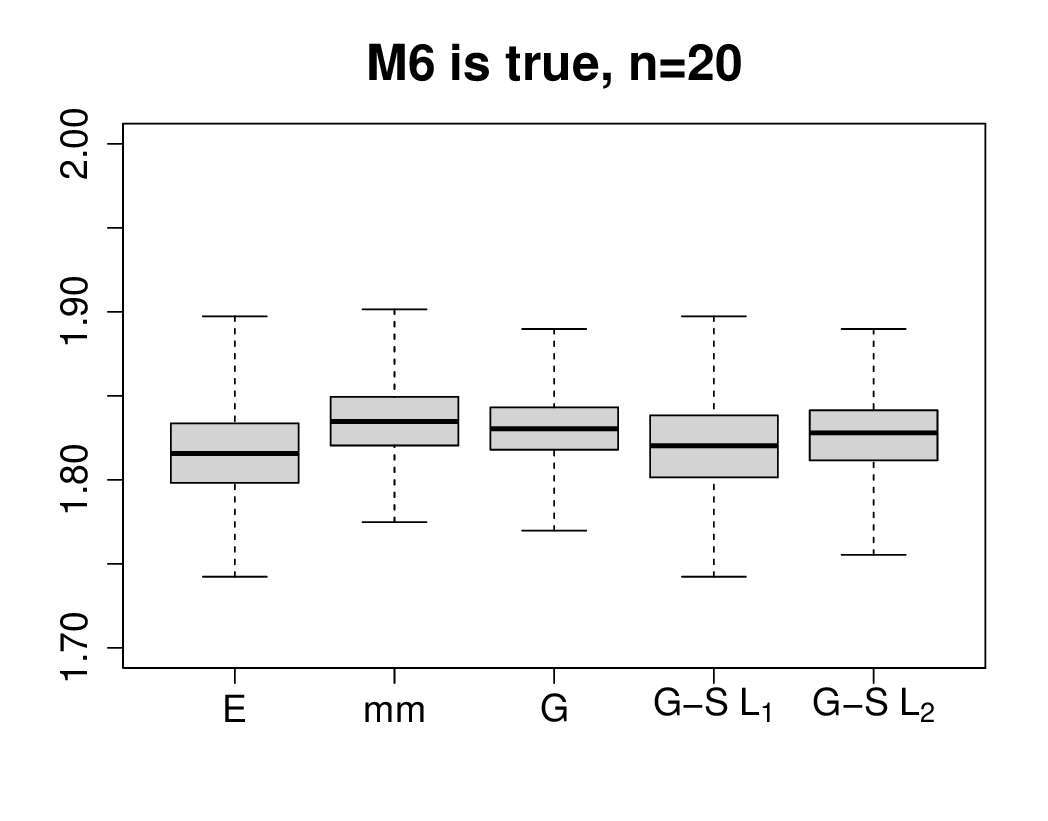}
  \end{subfigure}
  
    \begin{subfigure}{5.2cm}
    \centering\includegraphics[scale=0.27]{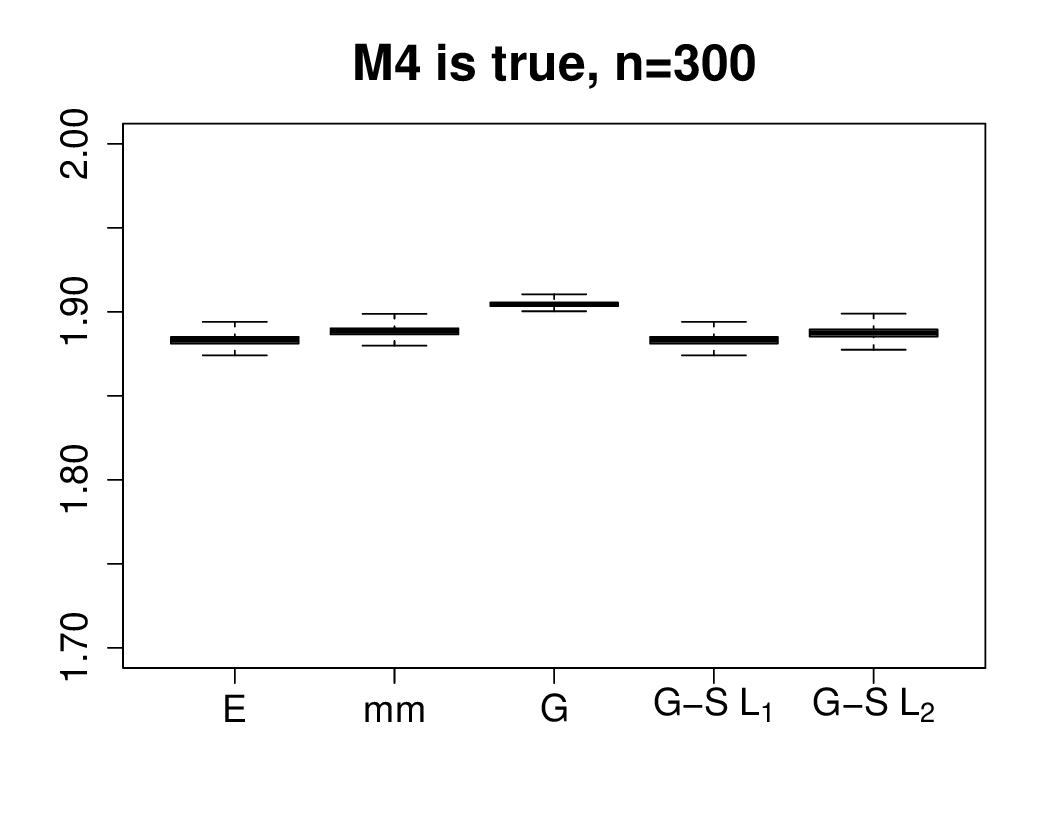}
  \end{subfigure}
  \begin{subfigure}{5.2cm}
    \centering\includegraphics[scale=0.27]{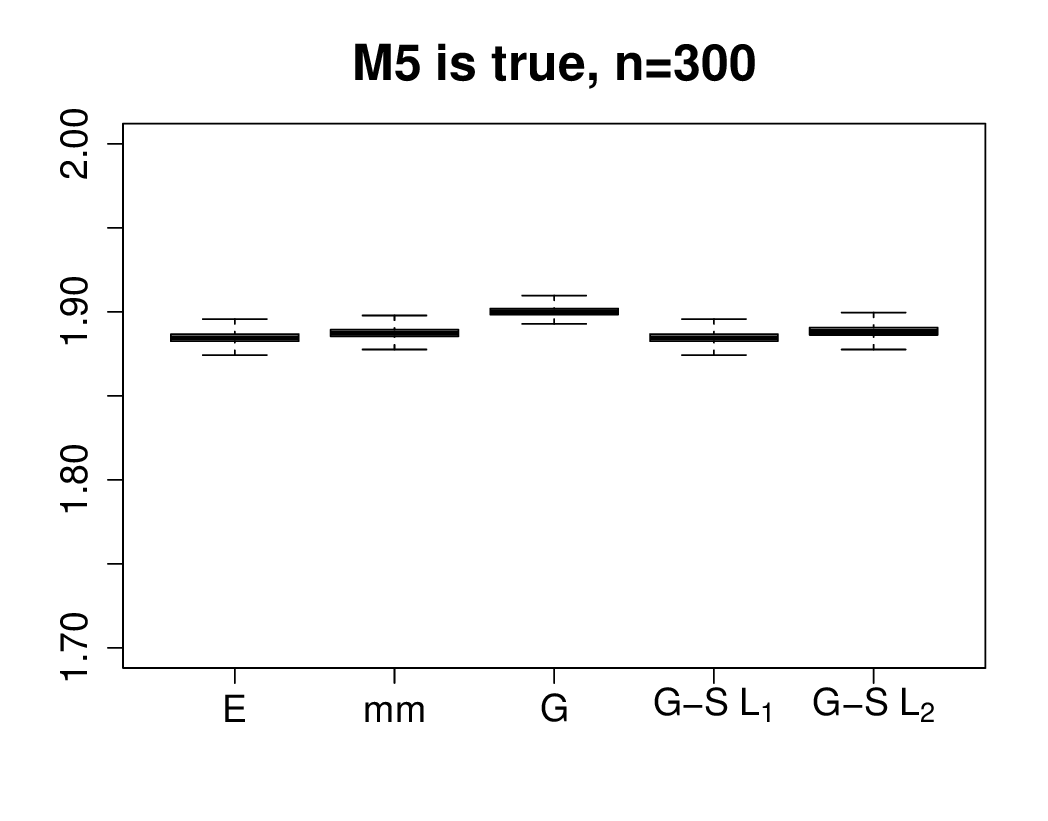}
  \end{subfigure}
  \begin{subfigure}{5.2cm}
    \centering\includegraphics[scale=0.27]{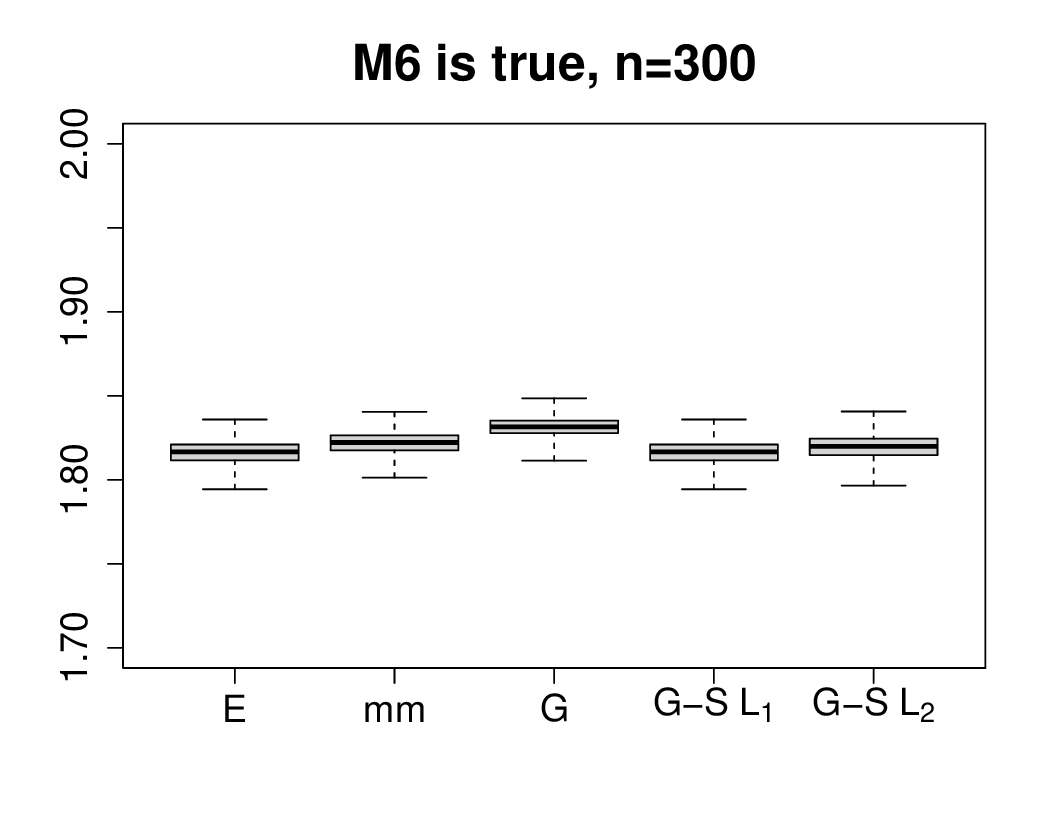}
  \end{subfigure}
  
    \begin{subfigure}{5.2cm}
    \centering\includegraphics[scale=0.27]{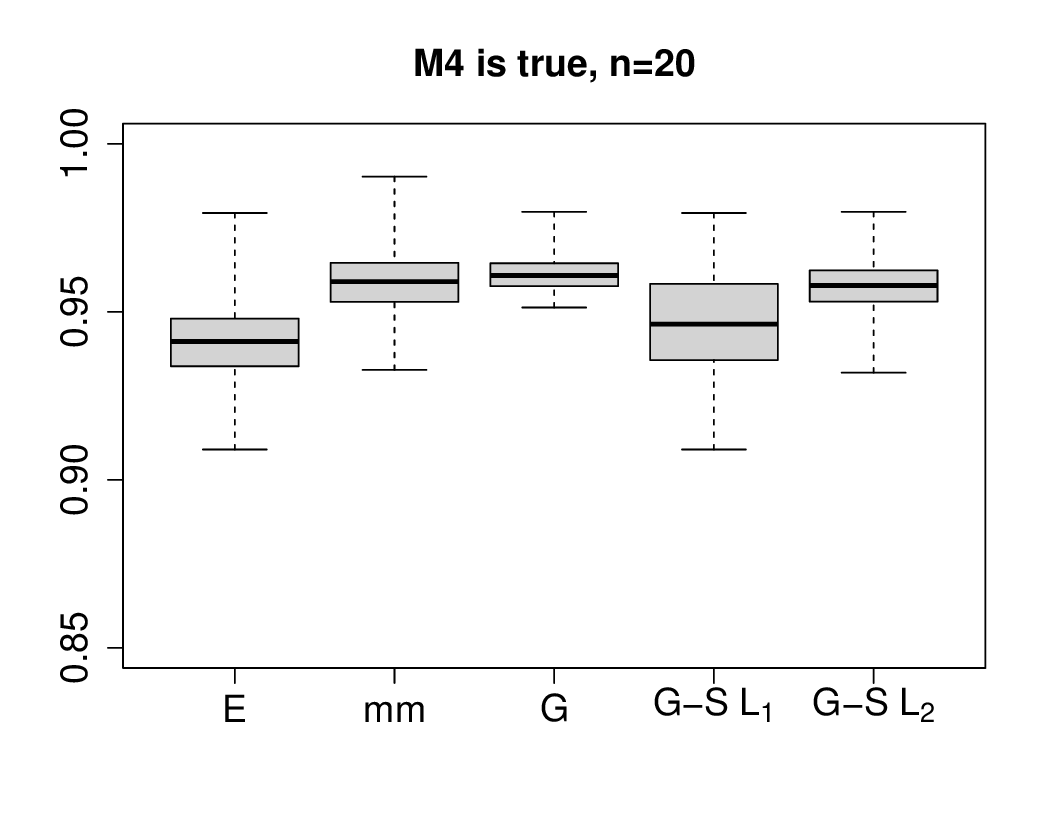}
  \end{subfigure}
  \begin{subfigure}{5.2cm}
    \centering\includegraphics[scale=0.27]{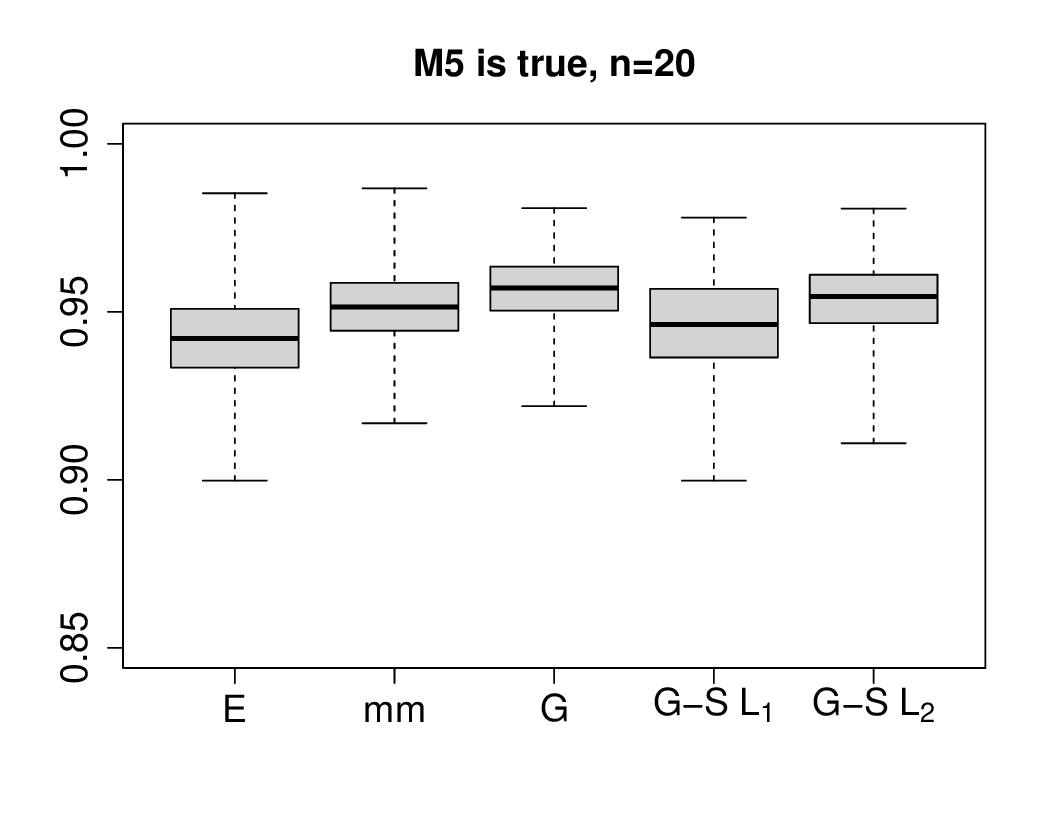}
  \end{subfigure}
  \begin{subfigure}{5.2cm}
    \centering\includegraphics[scale=0.27]{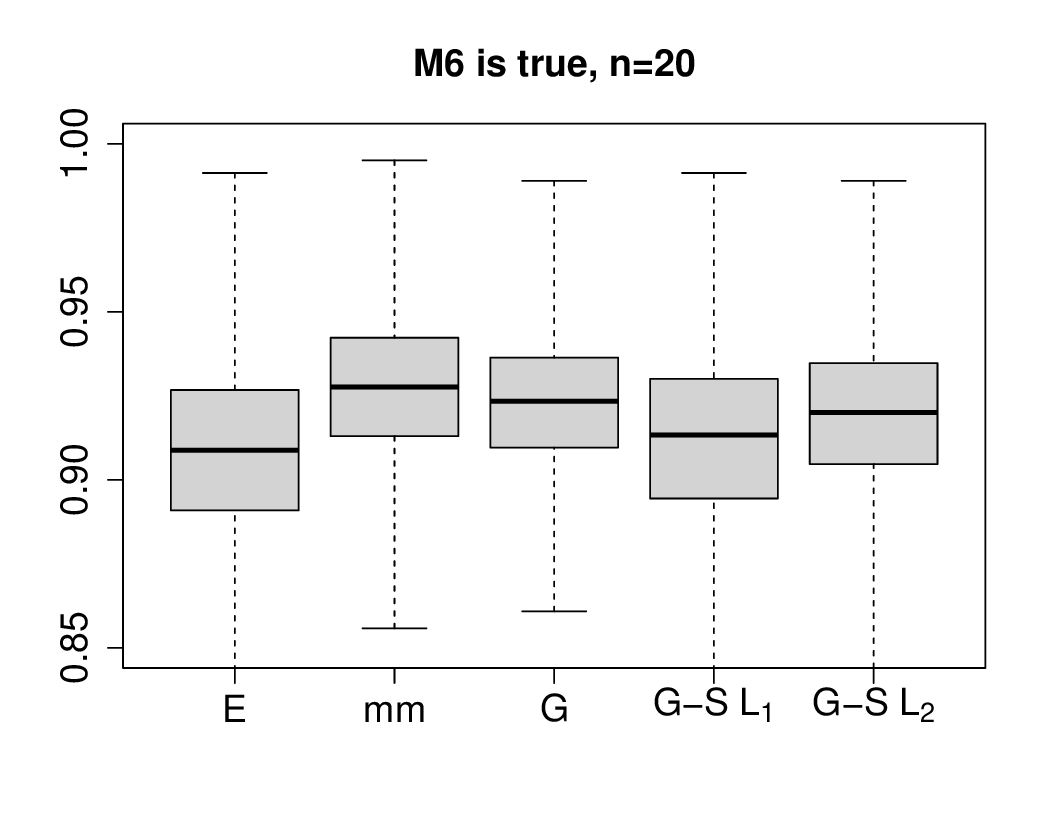}
  \end{subfigure}
  
    \begin{subfigure}{5.2cm}
    \centering\includegraphics[scale=0.27]{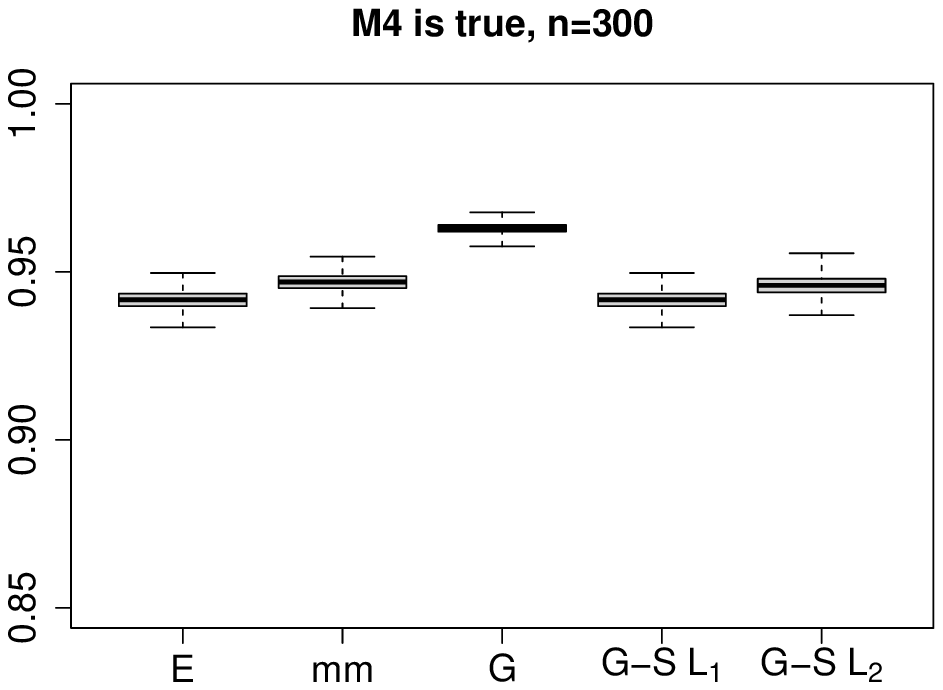}
  \end{subfigure}
  \begin{subfigure}{5.2cm}
    \centering\includegraphics[scale=0.27]{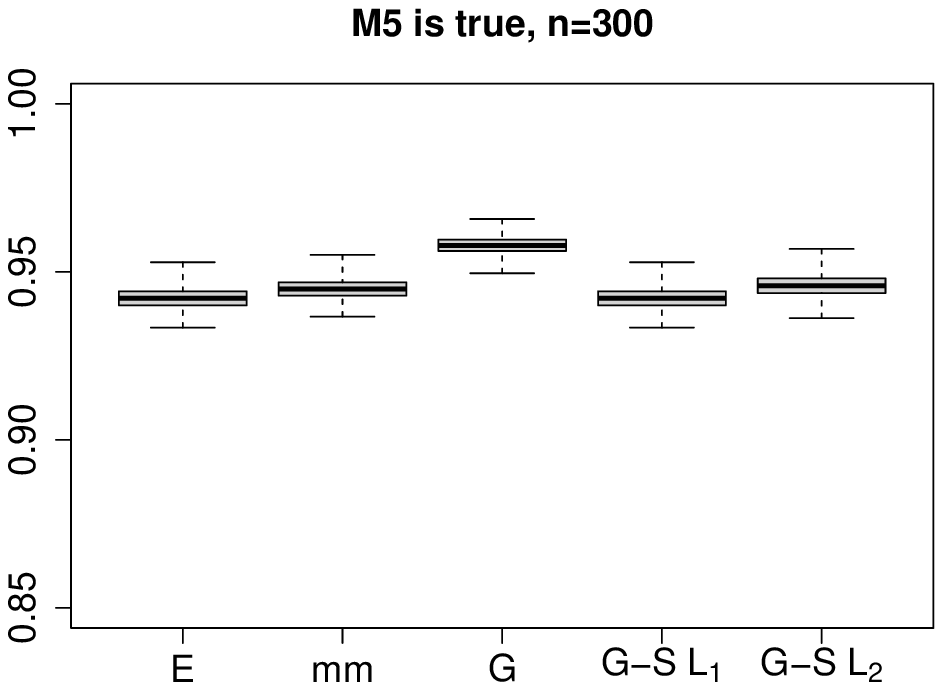}
  \end{subfigure}
  \begin{subfigure}{5.2cm}
    \centering\includegraphics[scale=0.27]{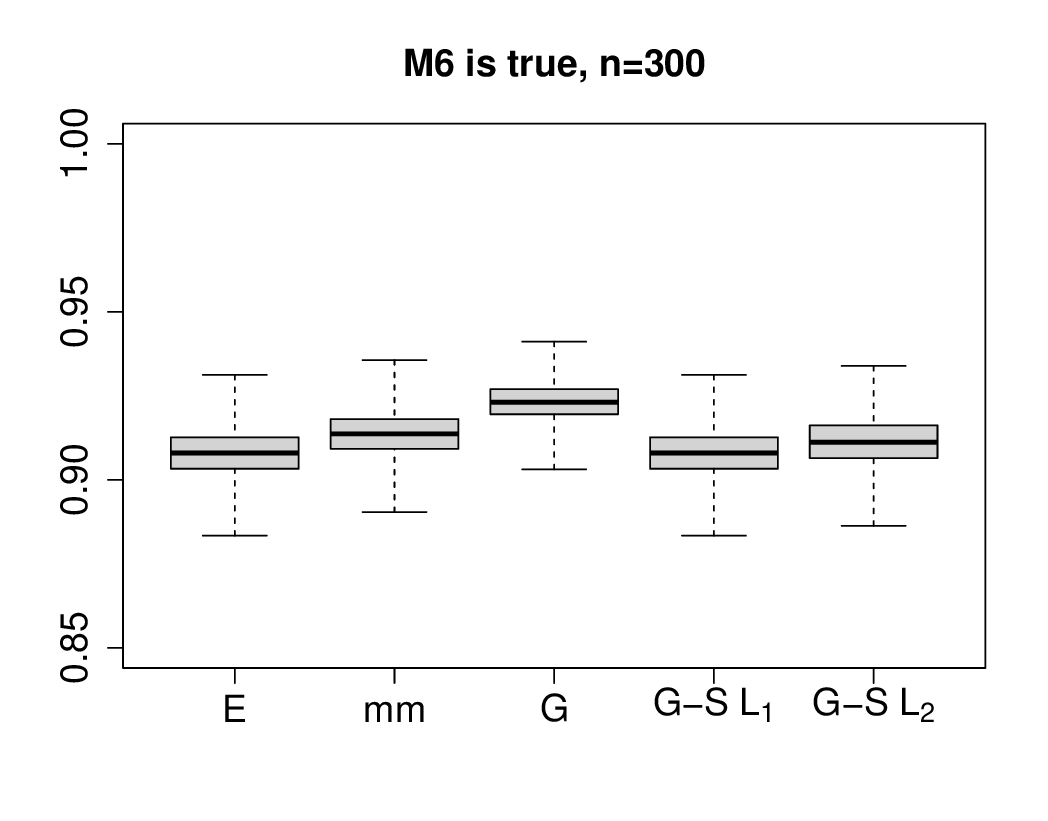}
  \end{subfigure}
   \caption{The boxplots of $1000$ expected $S_1$-scores  (first and second rows from the top) and $S_2$-scores (third and fourth rows) of the estimators: the empirical estimator $(E)$, minimax esimator $(mm)$, Grenander estimator $(G)$, Grenander-Stone estimator with $L_{1}$ penalization (G--S $L_{1}$) and Grenander-Stone estimator with $L_{2}$ penalization (G--S  $L_{2})$ for the models \textbf{M4}, \textbf{M5} and \textbf{M6}.}\label{nondecr_s1}
\end{figure}

\begin{figure}[t!] 
  \begin{subfigure}{5.2cm}
    \centering\includegraphics[scale=0.27]{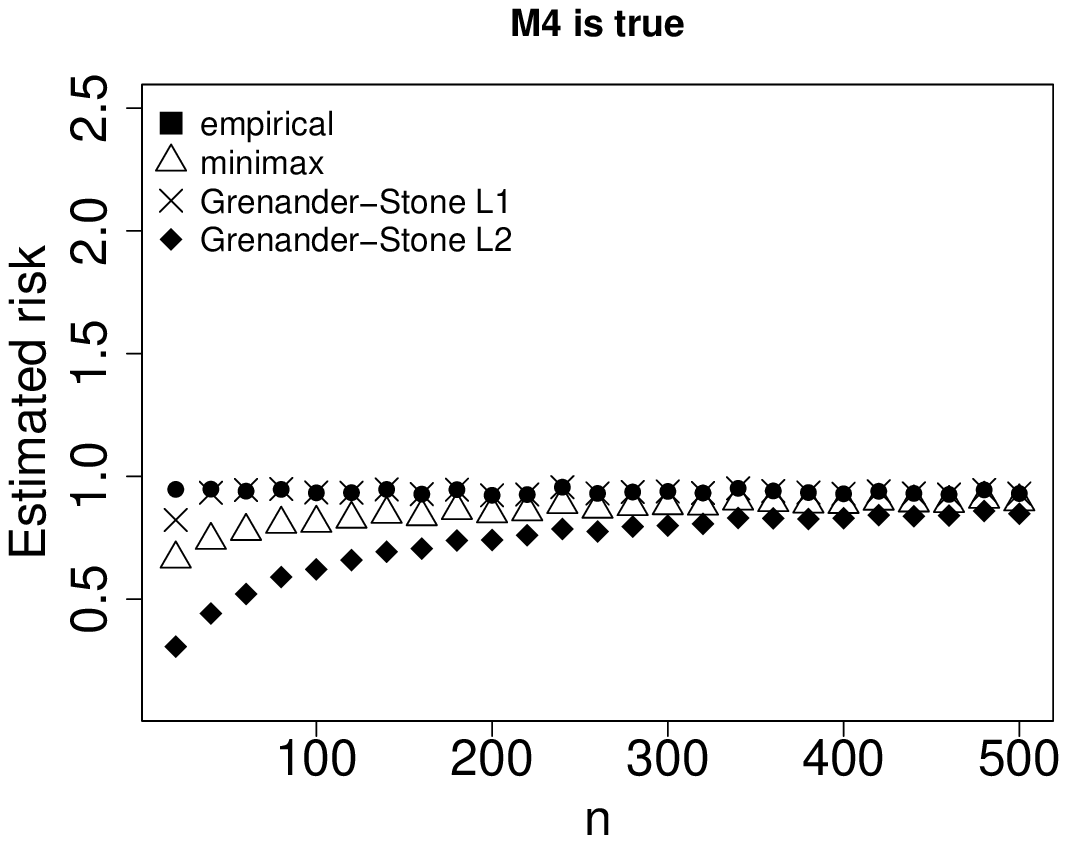}
  \end{subfigure}
  \begin{subfigure}{5.2cm}
    \centering\includegraphics[scale=0.27]{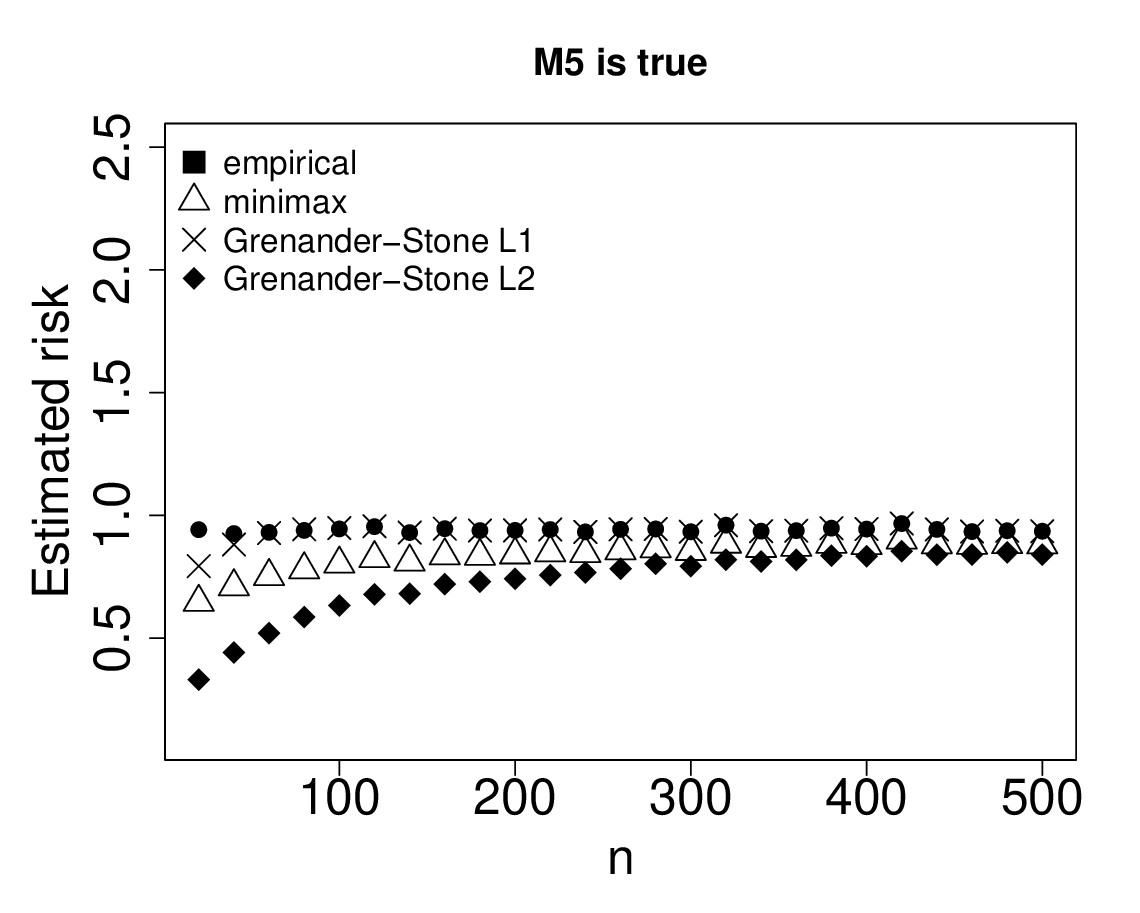}
  \end{subfigure}
  \begin{subfigure}{5.2cm}
    \centering\includegraphics[scale=0.27]{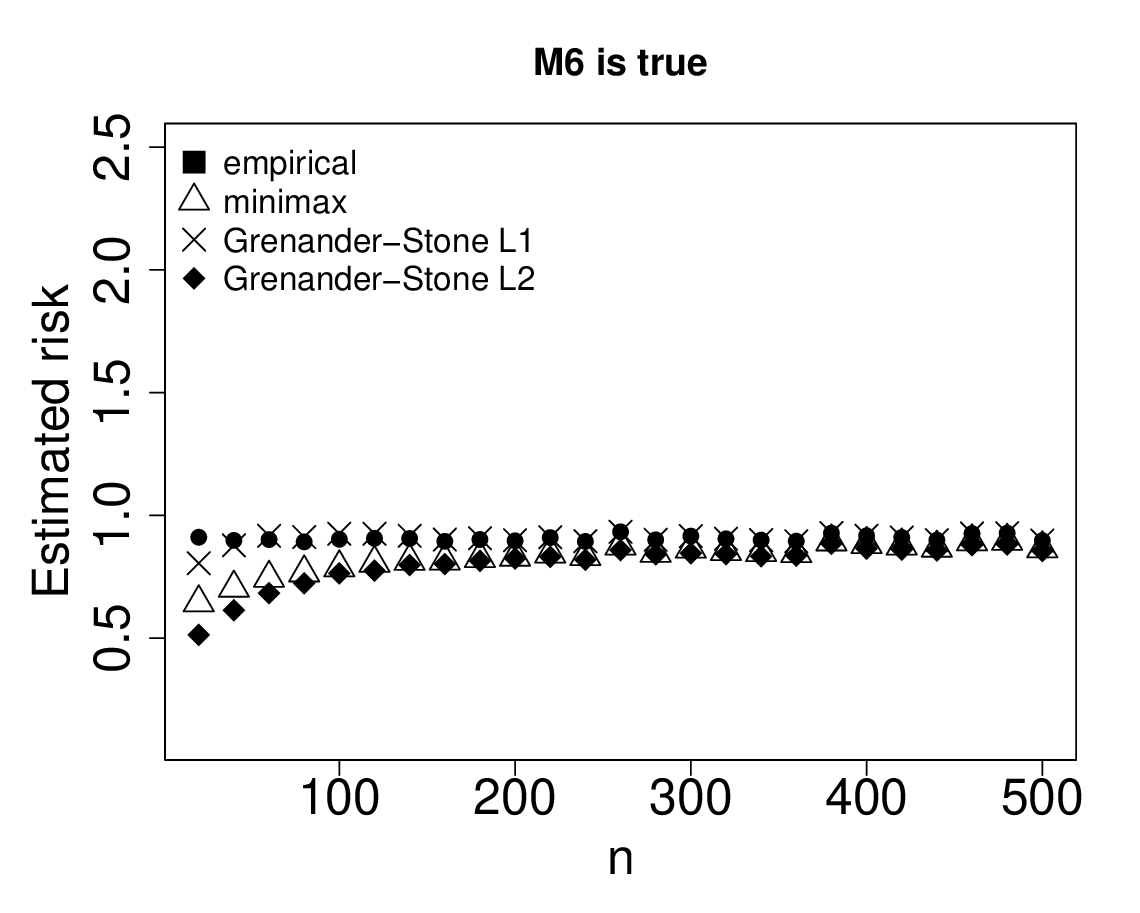}
  \end{subfigure}
   \caption{The estimates of the scaled risk for the models \textbf{M4}, \textbf{M5} and \textbf{M6}.}\label{ndecrRisk}
\end{figure}

The performance of Stone-Grenander estimator for non-monotone true distributions is presented at Figures \ref{nondecr_l1} and \ref{nondecr_s1}. From Figure \ref{nondecr_l1} it follows, first, Grenander-Stone estimator for the case of $L_{2}$ penalization is the overall winner. Second, Grenander estimator has a better performance than the empirical estimator for a small data set. The same effect holds for stacked Grenander estimator studied in \cite{pastukhov2022stacked}. This can be explained by reduction of the variance by isotonization, which compensates the bias in the case of small sized data sets. The explanation of $\ell_{2}$-risk reduction by stacking empirical estimator with some fixed probability vector is given in \cite{fienberg1973simultaneous}. Next, in a sense of forecast score, the winner is Grenander-Stone estimator with $L_{1}$ penalization, it performs in most of cases better than Grenander-Stone estimator with $L_{2}$ penalization, Grenander estimator and minimax estimator.  From Figure \ref{ndecrRisk} we can conclude that in the case of  non isotonic underlying distribution Grenander-Stone estimator with $L_{2}$ cross-validation performs better than the empirical and the minimax estimators for small sized data sets and performance becomes equal to the one of empirical estimator for large data sets.

  

\subsection{Performance of the global confidence band}
In Table \ref{emp-cov} we provide the empirical coverage probabilities for the proposed global confidence bands of the Grenander-Stone estimator, i.e. the table shows the proportion of times that
\begin{equation*}
\max\Big((\hat{\phi}_{n,j} - \frac{\hat{q}_{\alpha}}{\sqrt{n}}), 0\Big) \leq p_{j} \leq  \hat{\phi}_{n,j} + \frac{\hat{q}_{\alpha}}{\sqrt{n}}, \, \text{for all } j\in\mathbb{N}
\end{equation*}
among 1000 runs. The quantiles $\hat{q}_{\alpha}$ are estimated based on 100000 Monte-Carlo simulations. 

\begin{table}[!htbp]
\centering
\begin{tabular}{ |l|l|l|l|l|l|l|l| }
\hline
Estimator & Sample size & M1 & M2 & M3 & M4 & M5 & M6 \\ \hline
\multirow{2}{*}{Empirical estimator} 
 & n=100 & 0.958 & 0.95 &0.953 & 0.973 & 0.97 & 0.965 \\
 & n=1000 & 0.952 & 0.952 &0.944 & 0.945 & 0.957 & 0.965 \\ \hline
\multirow{2}{*}{Grenander-Stone estimator, $L_1$} 
 & n=100 & 0.959 & 0.966 & 0.966 & 0.973 & 0.97 & 0.953 \\
 & n=1000 & 0.953 & 0.955 &0.945 & 0.945 & 0.955 & 0.956 \\ \hline
\multirow{2}{*}{Grenander-Stone estimator, $L_2$} 
 & n=100 & 0.991 & 0.98 &0.974 & 0.997 & 0.997 & 0.975 \\
 & n=1000 & 0.984 & 0.955 &0.953 & 0.959 & 0.964 & 0.972 \\\hline
\end{tabular}
\caption{\label{emp-cov} Empirical coverage probabilities for the confidence bands for $\alpha = 0.05$.}
\end{table}

\section{Application to the real data}\label{app}
In this section we apply Grenander-Stone estimator to the data of the time periods from the onset of the symptoms to the notification of disease (Data set A) and time periods from notification to death (Data set B) for visceral leishmaniasis in Brazil collected from $2007$ to $2014$  \citep{maia2019premature}. 

In the case of Data set A, Grenander-Stone estimator with $L_{1}$ penalization coincides with the empirical estimator, i.e. $\beta^{(L_{1})} = 0$, while in the case of $L_{2}$ penalization the mixture parameter is $\beta^{(L_{2})} = 0.28$. For Data set B, $\beta^{(L_{1})} = 1$ and $\beta^{(L_{2})} = 1$, therefore, Grenander-Stone estimator for both $L_{1}$ and $L_{2}$ losses coincides with Grenander estimator.

\begin{figure}[h!] 
  \begin{subfigure}{7.7cm}
    \centering\includegraphics[scale=0.39]{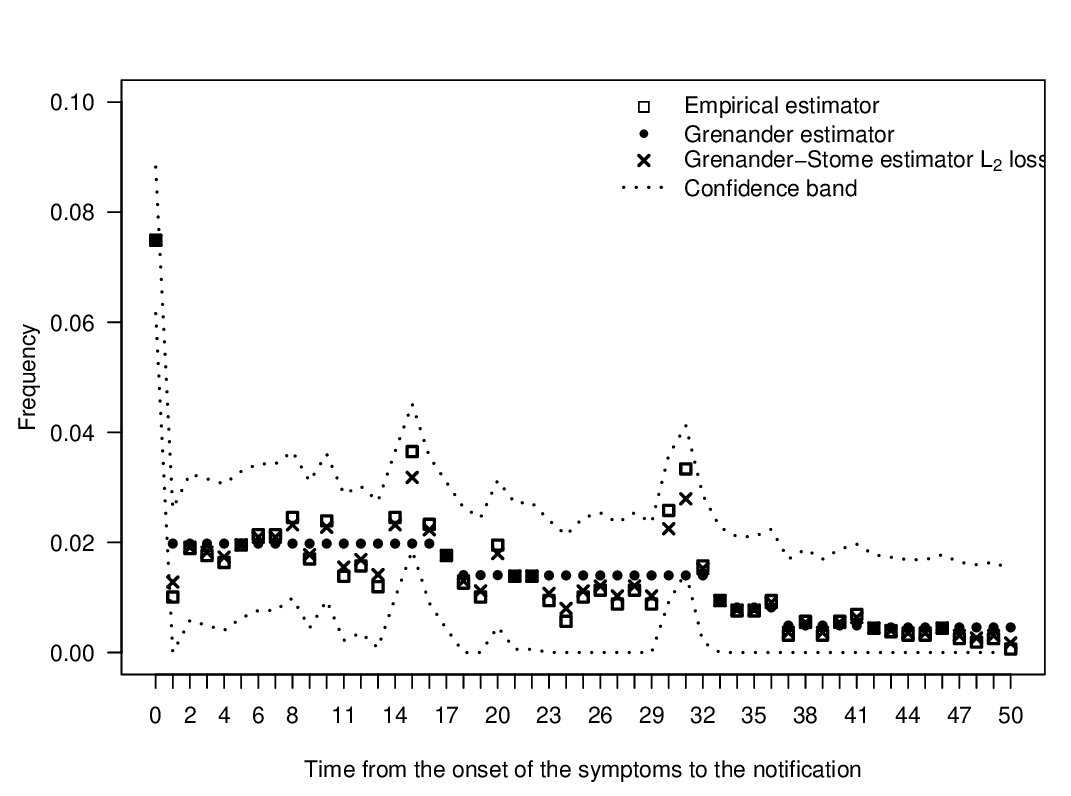}
  \end{subfigure}
  \begin{subfigure}{7.7cm}
    \centering\includegraphics[scale=0.39]{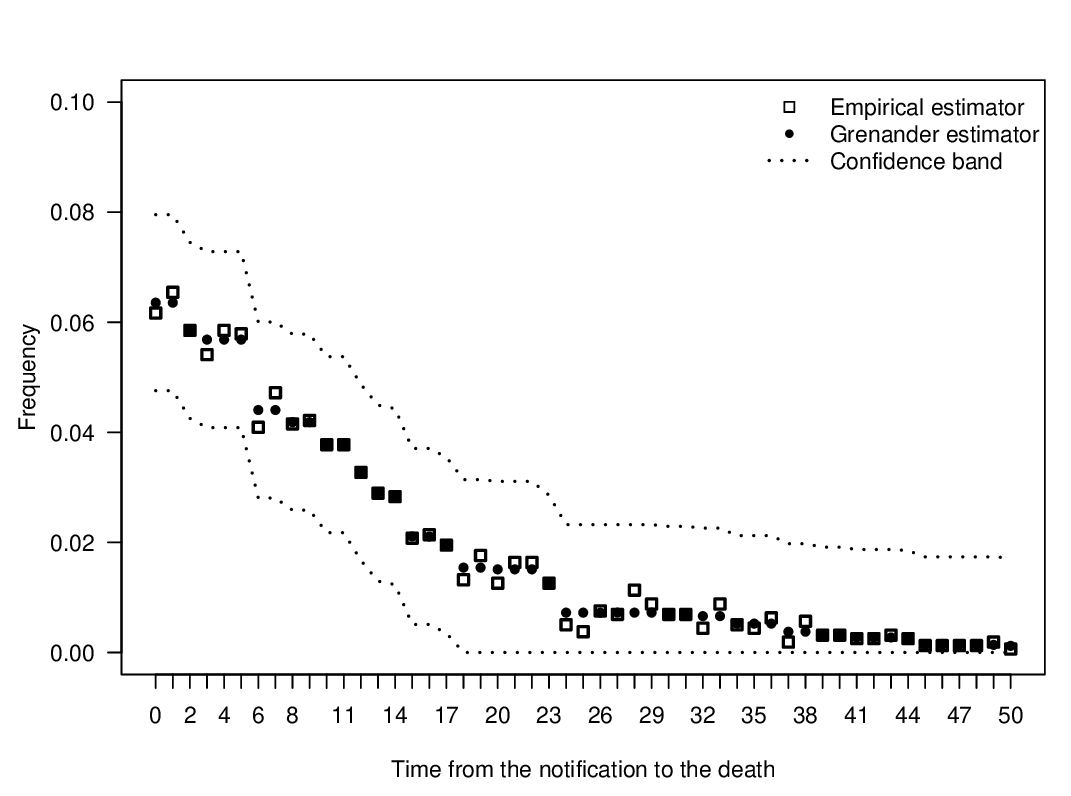}
  \end{subfigure}
    \caption{The frequency distribution of the time periods from onset of the symptoms to the notification of disease  (left panel), and from the notification to the death. The sample size is $n = 1589$, the largest order statistic is $t_{n} = 1207$, the range of the distribution shown on the plot is $j \in [0, 50]$.}\label{datafg}
\end{figure}

The estimators for Data sets A and B together with 95\% asymptotic global confidence bands are shown in the Figure \ref{datafg}.

\section{Discussion and future directions}\label{discus}
In this paper we provided and studied generalised framework of stacked constrained estimation of discrete infinitely supported distribution. The main focus was on estimation of isotonic distribution with respect to the partial order on its support.

The main conclusion of this paper is that Grenander-Stone estimator is computationally feasible and it has a good performance even for small data sets, while it is less restrictive than Grenander estimator. As noted in \cite{stone1974cross}, it is difficult to compare the performance $L_{1}$ and $L_{2}$ penalization. The case of $L_{1}$ penalization choses between the empirical estimator or  Grenander estimator, while $L_{2}$ penalization mostly provides a mixture of the estimators. 

The possible generalization of Grenander-Stone estimator is stacking of the empirical estimator with other shape constrained estimators, such as unimodal or convex maximum likelihood estimators. Second, the current approach can be applied to the histogram estimation of continuous distributions. Another direction of the research is stacking of constrained estimator based on other loss functions.

\bigskip
\begin{center}
{\large\bf SUPPLEMENTARY MATERIAL}
\end{center}

\section*{Proofs of all results}\label{appn} 
We start with properties of isotonic regression of a general vector in $\ell_{2}$ and proofs of auxiliary results which we need for the main proofs. Let $\bm{x} \in \ell_{2}$ be a vector indexed by some countable set $\mathcal{S} = \{\bm{s}_{0}, \bm{s}_{1}, \bm{s}_{2}, \dots \}$ with a partial order relation $\preceq$  defined on it. Next, let $\bm{x}^{*} \in \ell_{2}$ be the isotonic regression of $\bm{x}$ over the partially-ordered index set $\mathcal{S}$ i.e. 
\begin{eqnarray*}
\bm{x}^{*} = \underset{\bm{f} \in \mathcal{F}^{is}}{\arg \min} \sum_{\bm{s} \in \mathcal{S}}(f_{\bm{s}} - x_{\bm{s}})^{2},
\end{eqnarray*}
where $\mathcal{F}^{is} \subset \ell_{2}$ is the set of all isotonic vectors indexed by $\mathcal{S}$.

In Lemma \ref{propisot} below we provide properties of a general isotonic regression which we use in the proofs.
\begin{lemma}\label{propisot}[Properties of a general isotonic regression]
Let $\bm{x}^{*} \in \ell_{2}$ be isotonic regression of some vector $\bm{x} \in \ell_{2}$.
Then, the following holds:
\begin{enumerate}[label=(\roman*)]
\item $\bm{x}^{*}_{n}$ exists and it is unique.
\item $\sum_{j} x_{j} = \sum_{j} x^{*}_{j}$. 
\item $\bm{x}^{*}$, viewed as a mapping from $\ell_{2}$ into $\ell_{2}$, is continuous. 
\item $\bm{x}^{*}$ satisfies the same bounds as the basic estimator, i.e. $a \leq x^{*}_{j} \leq b$, for all $j = 1, 2, \dots$. 
\item $\Pi^{is}(a\bm{x} |\mathcal{F}^{is})= a\Pi^{is}(\bm{x}|\mathcal{F}^{is})$ for all $a \in \mathbb{R}_{+}$.
\item Error reduction property: for any $\bm{f} \in \mathcal{F}^{is}$ one has
\begin{equation*}\label{}
||\bm{x}^{*} - \bm{f}||_{k} \leq ||\bm{x} - \bm{f}||_{k},
\end{equation*} 
for all  $1 \leq k \leq \infty$.
\end{enumerate}  
\end{lemma}
\textbf{Proof}. Statements $(i)$, $(ii)$ and $(iii)$ follow from Theorem 8.2.1, Corollary B of Theorem 8.2.7 and Theorem 8.2.5, respectively, in \cite{robertsonorder}, statements $(iv)$, $(v)$ and $(vi)$ follow from Corollary B of Theorem 7.9, Theorems 7.5, respectively, in \cite{brunk1972statistical}, statement $(vi)$ follows from  Theorem 1.6.1 in \cite{robertsonorder}. \eop

Before we provide the solution to the isotonic regression with partially ordered support we need some more definitions. A subset $L$ of $\mathcal{S}$ is called a lower set with respect to $\preceq$ if $\bm{u} \in L$, $\bm{v} \in S$, $\bm{v} \preceq \bm{u}$ imply $\bm{v} \in L$. A subset $U$ of $\mathcal{S}$ is called a upper set with respect to $\preceq$ if $\bm{u} \in U$, $\bm{v} \in S$, $\bm{u} \preceq \bm{v}$ imply $\bm{v} \in U$. Further, let $\mathcal{L}$ and $\mathcal{U}$ denote the classes of all lower and upper sets, respectively. Both these classes are closed under arbitrarily union and intersection, and a set $U$ is an upper set if and only if $U^{c}$ is a lower set and, consequently, a set $L$ is a lower set if and only if $L^{c}$ is an upper set.  
 
A subset $B$ of $\mathcal{S}$ is a level set if there are a lower set $L$ and an upper set $U$ such that $B = L \cap U$. Note, that since $\mathcal{S}$ is both lower set and upper upper set, then each lower set and each upper set are level sets. Further, as it follows from Theorem 2.3 in \cite{brunk1972statistical}, a subset $B$ of $\mathcal{S}$ is a level set if and only if there exists a vector $\bm{x}^{*}$ indexed by $\mathcal{S}$, which is isotonic with respect to partial order $\preceq$ on $\mathcal{S}$ and there exists  $c \in \mathbb{R}$ such that $B = \{\bm{s} \in \mathcal{S}: x^{*}_{\bm{s}} = c \}$, i.e. $\bm{x}^{*}$ is constant on the index subset $B$. For a general introduction to constrained inference we refer to \cite{brunk1972statistical, robertsonorder}. 

Next, for a given vector $\bm{x}$ indexed by $\mathcal{S}$ and a subset $A$ of  $\mathcal{S}$ let 
\begin{equation*}
Av(\bm{x}, A) = \frac{\sum_{\bm{s} \in A} x_{\bm{s}}}{|A|}.
\end{equation*}

The next lemma describes the Minimum Lower Sets algorithm (cf. p. 76 in \cite{brunk1972statistical}) for the solution to the isotonic regression in the case of partially ordered support.

\begin{lemma}[Minimum Lower Sets algorithm]
For a given vector $\bm{x}$ indexed by index set $\mathcal{S}$ with partial order relation $\preceq$ defined on it the isotonic regression $\bm{x}^{*}$ is given by the following procedure.

First, we select the largest lower set $L_{1}$  such that 
\begin{equation*}
Av(\bm{x}, L_{1}) = \min \{Av(\bm{x}, L): L \in \mathcal{L} \}.
\end{equation*}
Let $B_{1} = L_{1}\cap \mathcal{S} \equiv L_{1}$, i.e. $B_{1}$ is the level set on which the isotonic regression $\bm{x}^{*}$ takes the smallest value:
\begin{equation*}
x^{*}_{\bm{s}} =  Av(\bm{x}, B_{1}) \quad \text{for} \quad \bm{s}\in B_{1}.
\end{equation*}

Second, let us consider the level sets such that $B_{2} = L \cap L_{1}^{c}$ and select largest $L_{2}$ such that  $B_{2} = L_{2}\cap L_{1}^{c}$ and
\begin{equation*}
Av(\bm{x}, B_{2}) = \min \{Av(\bm{x}, L \cap L_{1}^{c}): L \in \mathcal{L} \}.
\end{equation*}
Then, we have
\begin{equation*}
x^{*}_{\bm{s}} =  Av(\bm{x}, B_{2}) \quad \text{for} \quad \bm{s}\in B_{2}.
\end{equation*}

Further, we continue the process until the index set $\mathcal{S}$ is exhausted.
\end{lemma}
Also, analogously to one-dimensional case the algorithm described above can be written 
by one of the following max-min formulas, cf. Theorem 2.8 in \cite{brunk1972statistical}:
\begin{equation*}
\begin{aligned}
x^{*}_{\bm{s}} & = \underset{U: \bm{s} \in U \subset \mathcal{U} }{\max} \,\, \underset{L: \bm{s} \in L \subset \mathcal{L} }{\min} Av(\bm{x}, L \cap U) = 
\underset{U: \bm{s} \in U \subset \mathcal{U} }{\max}  \,\, \underset{L \subset \mathcal{L} }{\min} Av(\bm{x}, L \cap U) \\
& =  \underset{L: \bm{s} \in L \subset \mathcal{L} }{\min} \,\, \underset{U: \bm{s} \in U \subset \mathcal{U} }{\max} Av(\bm{x}, L \cap U)  =
\underset{L: \bm{s} \in L \subset \mathcal{L} }{\min}  \,\, \underset{U \subset \mathcal{U} }{\max} Av(\bm{x}, L \cap U).
\end{aligned}
\end{equation*}

\begin{lemma}\label{loovsis}
Let $\bm{x} \in \ell_{2}$ be a vector indexed by partially ordered set $\mathcal{S}$.  Next, let $\bm{y} \in \ell_{2}$ be a vector such that for some $m$ and for some $\bm{s}_{q} \in B_{m}$ we have:
\begin{equation*}\label{} 
  y_{\bm{s}_{j}}=\begin{cases}
   \tilde{y}_{\bm{s}_{j}} > x_{\bm{s}_{j}}, &\text{ if } \, j=q, \\
    x_{\bm{s}_{j}}, &\text{ otherwise}.
  \end{cases}
\end{equation*}
Then for the isotonic regressions $\bm{x}^{*}$ and $\bm{y}^{*}$ of $\bm{x}$ and $\bm{y}$, respectively, we have
\begin{equation*}
y^{*}_{\bm{s}_{j}} = x^{*}_{\bm{s}_{j}}, \, \text{ for all } \bm{s}_{j} \in B_{1}\cup\dots\cup B_{m-1}.
\end{equation*}
and 
\begin{equation*}
y^{*}_{\bm{s}_{j}} > x^{*}_{\bm{s}_{j}}, \, \text{ for } \bm{s}_{j} \in B_{m}.
\end{equation*}
\end{lemma}
\textbf{Proof}
First, let us $\{B_{1}, B_{2}, \dots \}$ denote the level sets of $\bm{x}^{*}$ and assume that $q \in B_{1}$. Then, for all $\bm{s} \in B_{1}$ we have 
\begin{equation*}
x^{*}_{\bm{s}} =  Av(\bm{x}, B_{1}).
\end{equation*}

Second, let us $\{ B^{'}_{1}, B^{'}_{2}, \dots \}$ denote the level sets of $\bm{y}^{*}$. Then, since $y_{\bm{s}_{q}} \geq x_{\bm{s}_{q}}$, we have 
\begin{equation*}
Av(\bm{x}, B_{1}) \leq Av(\bm{x}, B^{'}_{1}) \leq Av(\bm{y}, B^{'}_{1}).
\end{equation*}
Therefore, if $q \in \tilde{B}^{'} \neq B^{'}_{1}$, then
\begin{equation*}
y^{*}_{\bm{s}} =  Av(\bm{y}, \tilde{B}^{'}) \geq Av(\bm{y}, B^{'}_{1}) > Av(\bm{x}, B^{'}_{1}) \geq Av(\bm{x}, B_{1}) = x^{*}_{\bm{s}}.
\end{equation*}

Next, assume that $q \in B_{2}$. Then, first, note that in this case $ B^{'}_{1} = B_{1}$. Next, again, we have that for all $\bm{s} \in B_{2}$:
\begin{equation*}
x^{*}_{\bm{s}} =  Av(\bm{x}, B_{2}),
\end{equation*}
and
$q \in \tilde{B}^{'} \neq B^{'}_{2}$, then
\begin{equation*}
y^{*}_{\bm{s}} =  Av(\bm{y}, \tilde{B}^{'}) \geq Av(\bm{y}, B^{'}_{2}) > Av(\bm{x}, B^{'}_{2}) \geq Av(\bm{x}, B_{2}) = x^{*}_{\bm{s}}.
\end{equation*}
We continue the process until the index set $\mathcal{S}$ is exhausted.
\eop

\textbf{Proof of Theorem~\ref{thmbetacv}.}
We consider the cases of $L_{1}$ and $L_{2}$ separately. 

\textbf{The case of $L_{1}$ loss}. Recall the definitions of $\hat{\bm{\phi}}^{\backslash (i)}_{n}$ and $\hat{\bm{\phi}}^{\backslash [j]}_{n}$ at (\ref{cvestimi}) and (\ref{cvestimj}). Next, when $z_{i} = j \ (j=0, \dots, t_{n})$, then there exists $\bm{\delta}^{(i)} = \bm{\sigma}^{[j]}$, for some $i \ (i = 1, \dots, n)$, where $\sigma^{[j]}$ defined in (\ref{indgm}). Following the derivation in \cite{stone1974cross}, for the $L_{1}$ loss the leave-one-out cross-validation criterion $C_{1}(\beta)$ is given by
\begin{equation}\label{}
\begin{aligned}
CV_{1}(\beta) &={} \frac{1}{n}\sum_{i=1}^{n}L_{1}(\bm{\delta}^{(i)}, \hat{\bm{\phi}}^{\backslash (i)}_{n}) = \frac{1}{n}\sum_{i=1}^{n}||\bm{\delta}^{(i)} - \hat{\bm{\phi}}^{\backslash (i)}_{n}||_{1} =\frac{1}{n}\sum_{j=0}^{t_{n}}x_{j}||\bm{\sigma}^{[j]} - \hat{\bm{\phi}}^{\backslash [j]}_{n}||_{1} \\ 
&=\frac{1}{n}\sum_{j=0}^{t_{n}}x_{j}\big[ \big\{ 1 - \beta \hat{g}^{\backslash[j]}_{n, j} - (1-\beta)\frac{x_{j} - 1}{n-1} \big\}   +  \sum_{k \neq j}^{t_{n}}  \big\{ \beta \hat{g}^{\backslash[j]}_{n, k} - (1-\beta)\frac{x_{k}}{n-1} \big\} \big].\nonumber
\end{aligned}
\end{equation}
After the simplification we obtain
\begin{equation}\label{}
CV_{1}(\beta) = \frac{2}{n(n-1)}\beta \Big[ \sum_{j=0}^{t_{n}}\big\{ (x_{j} - n \hat{g}^{\backslash[j]}_{n, j})^{2} + (n+1)\hat{g}^{\backslash[j]}_{n, j}(x_{j} - n \hat{g}^{\backslash[j]}_{n, j}) + n(\hat{g}^{\backslash[j]}_{n, j})^{2} \big\} - n \Big] + c_{1},\nonumber
\end{equation}
where the term $c_{1}$ does not depend on $\beta$. Therefore, the result of Theorem \ref{thmbetacv} for $L_{1}$ loss follows.

\textbf{The case of $L_{2}$ loss}. For the leave-one-out cross-validation criterion in the case of $L_{2}$ loss we have
\begin{equation}\label{}
\begin{aligned}
CV_{2}(\beta) &={} \frac{1}{n}\sum_{i=1}^{n}L_{2}(\bm{\delta}^{(i)}, \hat{\bm{\phi}}^{\backslash (i)}_{n}) = \frac{1}{n}\sum_{i=1}^{n}||\bm{\delta}^{(i)} - \hat{\bm{\phi}}^{\backslash (i)}_{n}||_{2}^{2} =\frac{1}{n}\sum_{j=0}^{t_{n}}x_{j}||\bm{\sigma}^{[j]} - \hat{\bm{\phi}}^{\backslash [j]}_{n}||_{2}^{2}\\ 
&=\frac{1}{n}\sum_{j=0}^{t_{n}}x_{j}\big[ \big\{1 - \beta \hat{g}^{\backslash[j]}_{n, j} - (1-\beta)\frac{x_{j} - 1}{n-1} \big\}^{2}   +  \sum_{k \neq j}^{t_{n}}  \big\{ \beta \hat{g}^{\backslash[j]}_{n, k} + (1-\beta)\frac{x_{k}}{n-1} \big\}^{2} \big]\nonumber.
\end{aligned}
\end{equation}
Then, after simplification we get 
\begin{equation*}\label{}
CV_{2}(\beta) = \frac{1}{n(n-1)^{2}} (a_{n} \, \beta^{2}  - 2b_{n} \, \beta ) + c_{2},
\end{equation*}
where  the term $c_{2}$ does not depend on $\beta$, and 
\begin{equation*}
\begin{aligned}
a_{n} &={}  \sum_{j=0}^{t_{n}}x_{j}\big[ \big( x_{j} - 1 - \hat{g}^{\backslash[j]}_{n, j}(n-1)\big)^{2} + \sum_{k \neq j}^{t_{n}}\big(x_{k} - \hat{g}^{\backslash[j]}_{n, k}(n-1)\big)^{2}\big] \\
&= n - \sum_{j=0}^{t_{n}}\big[2x_{j}(x_{j} - \hat{g}^{\backslash[j]}_{n, j}(n-1)) - x_{j} \sum_{k=0}^{t_{n}}(x_{k} - (n-1)\hat{g}^{\backslash[j]}_{n, k})^{2} \big],
\end{aligned}
\end{equation*}
and 
\begin{equation*}
\begin{aligned}
b_{n} &={}  \sum_{j=0}^{t_{n}} x_{j}\big[ \Big\{ \sum_{k \neq j}^{t_{n}} x_{k}(x_{k} - \hat{g}^{\backslash[j]}_{n, k}(n-1)) \Big\} - (n-x_{j})(x_{j} - \hat{g}^{\backslash[j]}_{n, j}(n-1) -1) \big] \\
&=n^2 - \sum_{j=0}^{t_{n}}x_{j}^{2} + (n-1) \sum_{j=0}^{t_{n}} x_{j} \sum_{k=0}^{t_{n}} x_{k} ( \hat{g}^{\backslash[k]}_{n, k} - \hat{g}^{\backslash[j]}_{n, k}).  
\end{aligned}
\end{equation*}
First, from the derivation of $a_{n}$ it follows that $a_{n} \geq 0$. Therefore, when $a_{n} > 0$ the maximum of $CV_{2}(\beta)$ is achieved when
\begin{equation*}\label{}
  \beta^{(L_{2})}=\begin{cases}
    \frac{b_{n}}{a_{n}}, & \text{if } \, 0 \leq b_{n} \leq a_{n}, \\
    0, & \text{if } \,b_{n} < 0, \\
    1, & \text{otherwise}.
  \end{cases}
\end{equation*}

Next, unlike the case of a constant $\bm{\lambda}$ in Stone estimator $\hat{\bm{p}}^{S}_{n}$, in our case both coefficients $a_{n}$ and $b_{n}$ can be zero. Indeed,  $a_{n} =0$ when
\begin{equation*}
\big[ \big( x_{j} - 1 - \hat{g}^{\backslash[j]}_{n, j}(n-1)\big)^{2} + \sum_{k \neq j}^{t_{n}}\big(x_{k} - \hat{g}^{\backslash[j]}_{n, k}(n-1)\big)^{2}\big] = 0,
\end{equation*}
for all $j = 0, \dots, t_{n}$. Therefore, $a_{n} = 0$ iff for all $j = 0, \dots, t_{n}$
\begin{equation*}\label{condast0}
  \begin{cases}
    \hat{g}^{\backslash[j]}_{n, k} = \frac{x_{k} - 1}{n-1}, & \text{if} \, k = j, \\
    \hat{g}^{\backslash[j]}_{n, k} =  \frac{x_{k}}{n-1}, & \text{if} \, k \neq j,
  \end{cases}
\end{equation*}
which is equivalent to 
\begin{equation}\label{cond_thtj}
\bm{g}^{\backslash[j]}_{n} \equiv \Pi \big(\frac{\bm{x} - \bm{\sigma}^{[j]}}{n-1}\big)  = \frac{\bm{x} - \bm{\sigma}^{[j]}}{n-1},
\end{equation}
for all $j = 0, \dots, t_{n}$. If the condition (\ref{cond_thtj}) is satisfied, then also $b = 0$. 
Next, if $x$ is strictly decreasing, then $\hat{\bm{\phi}}_{n} = \hat{\bm{p}}_{n}$, and for consistency of notation we can set $\beta = 0$.
Finally, we obtain the following result for $\beta^{(L_{2})}_{n}$:
\begin{equation*}\label{}
  \beta^{(L_{2})}_{n} = \begin{cases}
    \frac{b_{n}}{a_{n}}, & \text{if } \, a_{n} \neq 0 \text{ and } 0 \leq b_{n} \leq a_{n}, \\
    0, & \text{if } \, a_{n} = 0 \text{ or } b_{n} < 0,\\
    1, & \text{otherwise}.
  \end{cases}
\end{equation*}
\eop

\textbf{Proof of Theorem~\ref{consrtgnr}.}
First, since the map $\Pi$ is a continuous map from $\ell_{2}$ to $\ell_{2}$, and the empirical estimator $\hat{\bm{p}}_{n}$ is strongly consistent, then from continuous mapping theorem it follows that 
\begin{equation*}
\Pi(\hat{\bm{p}}_{n})\stackrel{a.s.}{\to} \bm{g},
\end{equation*}
for some $\bm{g} = \Pi(\bm{p})$, in $\ell_{2}$-norm. Next, the almost sure convergence of $\Pi(\hat{\bm{p}}_{n})$ to $\bm{g}$ in $\ell_{k} \ (1 \leq k \leq \infty)$ follows from  Lemma C.2 in the supporting material of \cite{balabdaoui2013asymptotics}.

Next, we assume that $\bm{g} = \Pi(\bm{p}) = \bm{p}$. Then, from the subadditivity of norms it follows that
\begin{equation}\label{ineqnorm}
\begin{aligned}
||\hat{\bm{\phi}}_{n} - \bm{p}||_{k} {}  \leq \hat{\beta}^{(L_{d})}_{n}||\hat{\bm{g}}_{n} - \bm{p}||_{k} + (1-\hat{\beta}^{(L_{d})}_{n})||\hat{\bm{p}}_{n} - \bm{p}||_{k},
\end{aligned}
\end{equation} 
for $1 \leq k \leq \infty$ and $d = 1, 2 $. Therefore, the result of the theorem follows from statement $(v)$ of the Assumption \ref{asmt} and Theorem 2.4 and Corollary 4.1 in \cite{jankowski2009estimation}.

Further, we assume that $\bm{g} \neq \bm{p}$ and consider the cases of $L_{1}$ and $L_{2}$ losses for the cross-validation separately. 

\textbf{The case of $L_{1}$ loss}.
First, in the similar way as in Theorem 2 in \cite{pastukhov2022stacked} one can show that
\begin{equation*}\label{}
\hat{\bm{\gamma}}_{n}\stackrel{a.s.}{\to} \bm{g},
\end{equation*}
in $\ell_{k}$-norm for $1 \leq k \leq \infty$.

Recall the definition of $\beta^{(L_{1})}_{n}$ at Theorem \ref{thmbetacv}. Let us rewrite $\beta^{(L_{1})}_{n}$ as
\begin{equation*}\label{}
  \beta^{(L_{1})}_{n}= 
  \begin{cases}
    0, & \text{if }  B_{n} \geq n, \\
    1, & \text{otherwise},
  \end{cases}
\end{equation*}
where
\begin{equation*}\label{}
\begin{aligned}
B_{n} &= 
\sum_{j=0}^{t_{n}}(x_{j} - n \hat{\gamma}_{n, j})^{2} + \sum_{j=0}^{t_{n}} (n+1)\hat{\gamma}_{n, j}(x_{j} - n \hat{\gamma}_{n, j}) + \sum_{j=0}^{t_{n}} n(\hat{\gamma}_{n, j})^{2}\\
&\equiv B_{n}^{(1)} + B_{n}^{(2)} + B_{n}^{(3)}.
\end{aligned}
\end{equation*}
Next, from continuous mapping theorem, almost sure convergence of $\hat{\bm{\gamma}}_{n}$, and assumption $(iii)$ on the theorem it follows that 
\begin{equation*}\label{}
\begin{aligned} 
B_{n}^{(1)}/n^{2} &\stackrel{a.s.}{\to} ||(\bm{p} -  \bm{g})||_{2}^{2} > 0, \\
 B_{n}^{(2)}/n^{2} &\stackrel{a.s.}{\to} \langle \bm{g}, (\bm{p} -  \bm{g}) \rangle = 0, \\ 
  B_{n}^{(3)}/n^{2} &\stackrel{a.s.}{\to} 0.
\end{aligned}
\end{equation*}
Therefore, we have shown that 
\begin{equation*}\label{}
B_{n}/n^{2} \stackrel{a.s.}{\to} ||(\bm{p} -  \bm{g})||_{2}^{2},
\end{equation*} 
which means that there exists sufficiently large random $n_{0}$ such that $B_{n} > n$ for all $n > n_{0}$ almost surely, and, therefore, $\beta^{(L_{1})}_{n} = 0$
for all $n > n_{0}$ almost surely, which provides that
\begin{equation*}\label{}
\hat{\bm{\phi}}_{n} = \hat{\bm{p}}_{n} \text{ for all } n > n_{0} \text{  almost surely}.
\end{equation*} 
The strong consistency of $\hat{\bm{\phi}}_{n}$ follows from strong consistency of the empirical estimator $\hat{\bm{p}}_{n}$.

Finally, the result for the rate of convergence of $\hat{\bm{\phi}}_{n}$ follows from the rate of convergence of the empirical estimator i.e. as proved in Corollaries 4.1 and 4.2. in \cite{jankowski2009estimation}
\begin{equation*}
n^{1/2}||\hat{\bm{p}}_{n} - \bm{p}||_{k} = O_{p}(1) \quad (2 \leq k \leq \infty),
\end{equation*}
and, if in addition $\sum_{j=0}^{\infty}\sqrt{p}_{j} < \infty$ holds, then
\begin{equation*}
n^{1/2}||\hat{\bm{p}}_{n} - \bm{p}||_{1} = O_{p}(1).
\end{equation*}

\textbf{The case of $L_{2}$ loss}. Recall the definitions of $a_{n}$ and $b_{n}$ at Theorem \ref{thmbetacv}. Let us start with the asymptotic properties of $a_{n}$. Let us rewrite $a_{n}$ as
\begin{equation*}
a_{n} = \mathfrak{a}^{(1)}_{n} + \mathfrak{a}^{(2)}_{n},
\end{equation*}
where
\begin{equation*}
\begin{aligned}
\mathfrak{a}^{(1)}_{n} &= n - 2\sum_{j=0}^{t_{n}}x_{j}(x_{j} - \hat{g}^{\backslash[j]}_{n, j}(n-1)), \\
\mathfrak{a}^{(2)}_{n} &= \sum_{j=0}^{t_{n}}x_{j}\sum_{k=0}^{t_{n}}(x_{k} - (n-1)\hat{g}^{\backslash[j]}_{n, k})^{2}.
\end{aligned}
\end{equation*}
First, from continuous mapping theorem, almost sure convergence of $\hat{\bm{\gamma}}_{n}$ and the statement $(iii)$ in the Asumption \ref{asmt} we have  
\begin{equation*}
\mathfrak{a}^{(1)}_{n}/n^{2} \stackrel{a.s.}{\to}  -2\langle \bm{p}, (\bm{p} -  \bm{g}) \rangle < 0.
\end{equation*}
Second, since $\Pi$ is continuous map, then for any fixed $j$ we have $\hat{\bm{g}}^{\backslash[j]}_{n}  \stackrel{a.s.}{\to} \bm{g}$. Furthermore, since $\bm{g} \neq \bm{p}$, if $\bm{p}$ is not decreasing, then there exists sufficiently large random $n_{0}$ such that 
\begin{equation*}
\mathfrak{a}^{(2)}_{n}/n^{3} > 0  \text{ for all } n > n_{0} \text{ almost surely}. 
\end{equation*}
Therefore, there exists sufficiently large random $n_{1}$ such that
\begin{equation*}
a_{n}/n^{3} > 0, \text{ for all } n > n_{1} \text{ almost surely}.
\end{equation*}

Next, we study the asymptotic properties of $b_{n}$. Let us rewrite $b_{n}$ as
\begin{equation*}
b_{n} = \mathfrak{b}^{(1)}_{n} + \mathfrak{b}^{(2)}_{n},
\end{equation*}
where
\begin{equation*}
\begin{aligned}
\mathfrak{b}^{(1)}_{n} &= n^2 - \sum_{j=0}^{t_{n}}x_{j}^{2},\\
 \mathfrak{b}^{(2)}_{n} &= (n-1) \sum_{j=0}^{t_{n}} x_{j} \sum_{k=0}^{t_{n}} x_{k} ( \hat{g}^{\backslash[k]}_{n, k} - \hat{g}^{\backslash[j]}_{n, k}).
\end{aligned}
\end{equation*}
First, note that 
\begin{equation*}
\mathfrak{b}^{(1)}_{n}/n^{2.5} \stackrel{a.s.}{\to} 0. 
\end{equation*}
Second, since the map $\Pi$ is contractive (statement $(ii)$ of Assumption \ref{asmt}) with respect to $\ell_{\infty}$-norm, one has
\begin{equation*}
\sup_{k} | \hat{g}^{\backslash[k]}_{n, k} - \hat{g}^{\backslash[j]}_{n, k}| \leq \sup_{k} | \hat{p}^{\backslash[k]}_{n, k} - \hat{p}^{\backslash[j]}_{n, k}| = \frac{1}{n-1}.
\end{equation*}
Therefore,  for $\mathfrak{b}^{(2)}_{n}$ we have
\begin{equation*}
|\mathfrak{b}^{(2)}_{n}| \leq (n-1) \sum_{j=0}^{t_{n}} x_{j} \sum_{k=0}^{t_{n}} x_{k} \frac{1}{n-1} = n^{2},
\end{equation*}
which yields $b_{n}/n^{2.5} \stackrel{a.s.}{\to} 0$, and since $\beta^{(L_{2})}_{n} \geq 0$, it follows that 
\begin{equation}\label{nbtnconv}
n^{1/2} \beta^{(L_{2})}_{n} \stackrel{a.s.}{\to} 0.
\end{equation}

Finally, the inequality (\ref{ineqnorm}) and the strong consistency of the empirical estimator $\hat{\bm{p}}_{n}$ yields strong consistency of Grenander-Stone estimator for the case of $L_{2}$ loss penalization. The $n^{1/2}$-rate of convergence follows from (\ref{nbtnconv}) and the $n^{1/2}$-rate of convergence of the empirical estimator. 
\eop

\textbf{Proof of Lemma~\ref{propL1}.}
Let us consider some arbitrary $\bm{h} \in \mathcal{P}$ and let $t \leq \infty$ be the size of support of $\bm{h}$. First, note that 
\begin{equation*}\label{}
S_{1}(\bm{f}, \bm{p}) = 2 - \sum_{k=0}^{t} p_{k} (\beta h_{k} + (1-\beta)p_{k}) = 2 -\beta \sum_{k=0}^{t} p_{k}(h_{k} - p_{k}) + ||\bm{p}||_{2}^{2}.
\end{equation*} 
Next, since $\sum_{k=0}^{t} p_{k}(h_{k} - p_{k}) \geq 0$, then $S_{1}(\bm{f}, \bm{p})$ is minimised by $\beta =0$, which yields $\bm{f} = \bm{p}$. \eop

\textbf{Proof of Theorem~\ref{GScons}.}
From Theorem  \ref{consrtgnr} it follows that in order to prove the theorem, we have to prove all the statements of Assumption \ref{asmt}.

First, from the statements $(i)$, $(ii)$ and $(iv)$ of Lemma \ref{propisot} it follows that for any $\bm{x} \in \ell_{2}$, if $\bm{x} \in \mathcal{P}$, then for the isotonic regression $\bm{x}^{*}$ of $\bm{x}$ one has $\bm{x}^{*} \in \mathcal{P}$. Therefore, we have 
\begin{equation*}\label{}
\Pi^{is}(\hat{\bm{p}}_{n}  \mid \mathcal{P}^{is}) = \Pi^{is}(\hat{\bm{p}}_{n}  \mid \mathcal{F}^{is}).
\end{equation*} 

Now, we prove the statements.

Statement $(i)$. It follows from statement $(iii)$ of Lemma \ref{propisot}.

Statement $(ii)$ follows from contraction property of isotonic regression with respect to $\ell_{\infty}$-norm, which was proved in Lemma 1 in \cite{yang2019contraction} in one dimensional case, and Lemma 9 in the Supplement to \cite{pananjady2022isotonic} for the case of a general partial order.
 
Statement $(iii)$ follows from the property of $\ell_{2}$ projection onto the cone.

Statement $(iv)$. Note that 
\begin{equation*}\label{}
\begin{aligned}
  \hat{\gamma}_{n,j}=\begin{cases}
    \Pi^{is}(\hat{\bm{p}}^{\backslash [j]}_{n}|\mathcal{F}^{is}), & \text{if } \, x_{j} > 0, \\
    0, & \text{otherwise},
      \end{cases} 
\end{aligned}
\end{equation*}
Next, for a fixed $j$ let $\bm{x} = (x_{1}, \dots, x_{j} -1, x_{t_{n}})$ and $\bm{y} = (x_{1}, \dots, x_{j}, x_{t_{n}})$, where $x_{i}$ are the sample frequences. Then, from Lemma \ref{loovsis} if follows that for isotonic regressions of $\bm{x}$ and $\bm{y}$ one has $x^{*}_{j} < y^{*}_{j}$. Next, note that
\begin{equation*}\label{}
(n-1)\frac{x^{*}_{j}}{n-1} < n \frac{y^{*}_{j}}{n},
\end{equation*}
and, using the property $(v)$ of isotonic regression from Lemma \ref{propisot} we prove the statement.

Statement $(v)$ follows from the inequality (\ref{ineqnorm}), the error reduction property of isotonic regression in $(vi)$ of Lemma \ref{propisot} and  from Theorem 2.1 in \cite{jankowski2009estimation}. \eop
 
\textbf{Proof of Theorem~\ref{Marshal}.}
First, note that from the max-min formulas in (\ref{minmaxis}) it follows that
\begin{equation}\label{suppis}
\sup\{j: \hat{g}_{n,j}>0 \} = \sup\{j: \hat{p}_{n,j}>0 \} \equiv t_{n}
\end{equation}
for any $n$. Next, let $D_{j}$ be a cumulative sum function of any decreasing distribution on $\mathbb{N}$. From (\ref{suppis}) it follows that for a given sample $\mathcal{Z}_{n}$ we have
\begin{equation*}\label{}
\sup_{j\in[t_{n}+1, \infty]} |\hat{F}^{G}_{n, j}  - D_{j}| =
\sup_{j\in[t_{n}+1, \infty]}|\hat{F}_{n, j}  - D_{j}| = 
\sup_{j\in[t_{n}+1, \infty]}|1 - D_{j}|,
\end{equation*}
where  $\hat{F}_{n,j} = \sum_{k=0}^{j}\hat{p}_{n, j} \ (j \in \mathbb{N})$, and $\hat{F}^{G}_{n,j} = \sum_{k=0}^{j}\hat{g}_{n, j} \ (j \in \mathbb{N})$.

Second, let $\hat{F}_{n,j}^{GS} = \sum_{k=0}^{j}\hat{\phi}_{n, j} \ (j \in \mathbb{N})$, and note that from the subadditivity of norms it follows that
\begin{equation*}\label{}
\sup_{j\in \mathbb{N}}|\hat{F}_{n,j}^{GS}  - D_{j}| \leq  \beta^{L_{d}}_{n}\sup_{j\in \mathbb{N}}|\hat{F}^{G}_{n, j}  - D_{j}| + (1 - \beta^{L_{d}}_{n})\sup_{j\in \mathbb{N}}|\hat{F}_{n, j}  - D_{j}|.
\end{equation*} 

Third, from Lemma in Section 2.2 in \cite{brunk1972statistical} it follows that for a cumulative sum function $D_{j}$ of any decreasing probability vector we have
\begin{equation*}\label{}
\sup_{j\in[0, t_{n} ]}|\hat{F}^{G}_{n, j}  - D_{j}| \leq \sup_{j\in[0, t_{n} ]}|\hat{F}_{n, j}  - D_{j}|,
\end{equation*} 
which completes the proof.
\eop

\textbf{Proof of Theorem~\ref{convqntl}.}
From Proposition B.7 at \cite{balabdaoui2016maximum} it follows that we have to prove:\\
(a) $\hat{\bm{\phi}}_{n}$ converges to $\bm{p}$ point-wise almost surely,\\
(b) $\lim_{m\to\infty} \lim_{n\to\infty} \sum_{j > m} \hat{\phi}_{n, j} = 0$.

First, we assume that the estimator is Grenander-Stone estimator if one dimensional case with the support $\mathcal{S} = \{0,1,2, \dots \}$. Then, the Condition (a) follows from Theorem \ref{consrtgnr}. Let us prove the Condition (b). Let $\hat{F}_{n, j} \ (j \in \mathbb{N})$, be the empirical distribution function, let 
\begin{equation*}
\hat{F}_{n, j}^{G} = \sum_{k=0}^{j}\hat{g}_{n, j} \ (j \in \mathbb{N}),
\end{equation*} 
and 
\begin{equation*}
\hat{F}_{n, j}^{GS} = \sum_{k=0}^{j}\hat{\phi}_{n, j} \ (j \in \mathbb{N}).
\end{equation*} 

Next, from the properties of a simple order isotonic regression it follows that $\hat{F}_{n, j}^{GS}$ is given by the least concave majorant of $\hat{F}_{n, j}$  \citep{brunk1972statistical, robertson1980algorithms}. Therefore, we have $\hat{F}_{n, j}^{G} \geq \hat{F}_{n, j}$ for any $j$ and $n$. Therefore, 
\begin{equation*}
\sum_{j > m} \hat{\phi}_{n, j} = 1 - \hat{F}_{n, m}^{GS} = 1 - \beta^{(L_{d})}_{n}\hat{F}_{n, j}^{G}  - (1-\beta^{(L_{d})}_{n})\hat{F}_{n, j} \leq 1 - \hat{F}_{n, m},
\end{equation*}
for $d = 1, 2$, and from the properties of the empirical distribution function it follows 
\begin{equation*}
\lim_{m\to\infty} \lim_{n\to\infty} (1 - \hat{F}_{n, m}) = 0,
\end{equation*}
almost surely.

Now we can show that the Conditions (a) and (b) hold for any strongly consistent estimator $\hat{\bm{\phi}}_{n}$ of $\bm{p}$. Note, that since the sum is bounded linear functional, the it is continuous. Therefore, for any strongly consistent estimator $\hat{\bm{\phi}}_{n}$ we have 
\begin{equation*}
\sum_{j > m} \hat{\phi}_{n, j} \stackrel{a.s.}{\to} \sum_{j > m} p_{j},
\end{equation*}
and, finally, 
\begin{equation*}
\lim_{m\to\infty} \sum_{j > m} p_{j} = 0.
\end{equation*}

\bibliographystyle{agsm}

\bibliography{gsest}{9}
\end{document}